\input ./style/arxiv-general.cfg
\documentclass[aap,MSNbibl,seceqn,dvips]{arximspdf}
\makeatletter
   \@ifpackageloaded{graphicx}{}{\usepackage{graphicx}}
\makeatother

\doi{10.1214/14-AAP1038} 
\volume{25}
\issue{4}
\pubyear{2015}
\firstpage{1868}
\lastpage{1935}
\docsubty{FLA}

\makeatletter

\newcommand{\rrvert}{\vert}
\newcommand{\llvert}{\vert}
\def\m{\mu} 

\def\mbb{\mathbb}
\def\mbf{\mathbf}
\def\mrm{\mathrm}
\def\mc{\mathcal}
\def\unl{\underline}
\def\ovl{\overline}

\def\P{\mathbb{P}}
\def\E{\mathbb{E}}
\def\R{\mathbb{R}}
\def\FF{\mathcal{F}}

\def\T{\tilde}
\def\WT{\widetilde}

\def\a{\alpha}
\def\b{\beta}
\def\d{\delta}
\def\th{\theta}
\def\s{\sigma}
\def\ps{\psi}
\def\F{\Phi}
\def\td{\mathrm{d}}
\def\Eq{\Leftrightarrow}

\newtheorem{Thm}{Theorem}[section]
\newtheorem{Cor}[Thm]{Corollary}
\newtheorem{Lemma}[Thm]{Lemma}
\newtheorem{lemma}[Thm]{Lemma}
\newtheorem{Prop}[Thm]{Proposition}
\newproclaim{As}{Assumption}
\newproclaim{Ex}[Thm]{Example}
\newproclaim{Def}[Thm]{Definition}
\newproclaim{Rem}[Thm]{Remark}

\newcommand{\eqref}[1]{(\ref{#1})}

\makeatother

\begin{document}
\begin{frontmatter}

\title{On Gerber--Shiu functions and optimal dividend distribution for
a L\'{e}vy risk process in the presence of a penalty function}
\runtitle{Optimal dividend distribution under a penalty}

\begin{aug}
\author[A]{\fnms{F.}~\snm{Avram}\ead[label=e1]{Florin.Avram@univ-Pau.fr}},
\author[B]{\fnms{Z.}~\snm{Palmowski}\ead[label=e2]{zpalma@math.uni.wroc.pl}}\thanksref{T1}
\and
\author[C]{\fnms{M.~R.}~\snm{Pistorius}\corref{}\thanksref{T2}\ead[label=e3]{m.pistorius@imperial.ac.uk}}
\runauthor{F. Avram, Z. Palmowski and M.~R. Pistorius}
\affiliation{University of Pau, University of Wroc\l aw and
Imperial College London}
\thankstext{T1}{Supported by the Ministry of Science and
Higher Education of Poland under the Grant DEC-2011/01/B/HS4/00982
(2012--2013).}
\thankstext{T2}{Supported in part by EPSRC
Grant EP/D039053, EPSRC Mathematics Platform Grant
EP/I019111/1 and NWO-STAR.}
\address[A]{F. Avram\\
Departement de Mathematiques\\
Universit\'e de Pau\\
France 64000\\
\printead{e1}}
\address[B]{Z. Palmowski\\
Mathematical Institute\\
University of Wroclaw\\
pl. Grunwaldzki 2/4\\
50-384 Wroclaw\\
Poland\\
\printead{e2}}
\address[C]{M.~R. Pistorius\\
Department of Mathematics\\
Imperial College London\\
South Kensington Campus\\
London SW7 2AZ\\
United Kingdom\\
\printead{e3}}
\end{aug}

\received{\smonth{12} \syear{2012}}
\revised{\smonth{5} \syear{2014}}

%
\begin{abstract}
This paper concerns an optimal dividend distribution problem for an
insurance company whose risk process evolves as a spectrally
negative L\'{e}vy process (in the absence of dividend payments).
The management of the company is assumed to control timing and
size of dividend payments. The objective is to maximize the sum of
the expected cumulative discounted dividend payments received until the
moment of ruin and a penalty payment at the moment of ruin, which
is an increasing function of the size of the shortfall at ruin; in
addition, there may be a fixed cost for taking out dividends.
A complete solution is presented to the corresponding stochastic
control problem. It is established that
the value-function is the unique \emph{stochastic solution} and
the pointwise smallest \emph{stochastic supersolution}
of the associated HJB equation.
Furthermore, a necessary and sufficient condition is identified for optimality
of a single dividend-band strategy, in terms of a particular
Gerber--Shiu function.
A number of concrete examples are analyzed.
\end{abstract}

%

%
\end{frontmatter}

\section{Optimal control of L\'{e}vy risk models. The spectrally
negative L\'{e}vy
risk model}\label{sec:intro}
Recall the classical Cram\'er--Lundberg model
%
\begin{equation}
\label{cramermod} X_t - X_0= {\eta} t- S_t,\qquad
S_t=\sum_{k=1}^{N_t}
C_k - \lambda m t,
\end{equation}
which is used in collective risk theory (e.g., Gerber
\cite{Gerber}) to describe the surplus $X=\{X_t, t\in\mbb R_+\}$
of an insurance company. Here, $X_0\ge0$ is the initial level of reserves,
$C_k$ are i.i.d. positive random
variables representing the claims made, $ N=\{N_t, t\in\mbb R_+\}$
is an independent Poisson process with intensity $\lambda$
modeling the times at which the claims occur, and $p t$, with
$p:=\eta+ \lambda m$, represents the premium income up to
time $t$, with profit rate $\eta>0$ and mean $m<\infty$ of $C_1$.

In later years, model (\ref{cramermod}) was generalized
to the ``perturbed model''
%
\begin{equation}
X_t-X_0:= \s B_t + \eta t -
S_t,\label{SNmod}
\end{equation}
where $B_t$ denotes an
independent standard Brownian motion, which models small scale
fluctuations of the risk process.

Since the jumps of $X$ are all negative, the moment generating
function $\E[\mathrm{e}^{\th(X_t-X_0)}]$ exists for all $\th\ge0$ and
$t\in
\mbb R_+$,
and is log-linear in $t$, defining thus a function $\psi(\th)$
satisfying $\E[\mathrm{e}^{\th( X_t-X_0)}] =
\mathrm{e}^{t\psi(\th)}$ with
%
\begin{equation}
\label{eq:psi}\psi(\th)= \frac{\s^2}{2} \th^2 + \eta\th+ \int
_{\mbb R_+\setminus\{0\}} \bigl( \mathrm{e}^{- \th x}- 1 + \theta x\bigr)
\nu(\td x),
\end{equation}
where $\nu(\td x)=\lambda F_C(\td x)$, $x \in\mathbb R_+$,
with $F_C$ the distribution function of $C_1$, is the ``L\'evy
measure'' of the
compound Poisson process $S_t$, and $\eta=\psi'(0)$ is the mean of $X_1-X_0$.

The cumulant exponent $\psi(\th)$
is well defined, at least on the positive half-line, where it is
strictly convex with the property that
$\lim_{\th\to\infty}\ps(\th)=+\infty$. Moreover, $\psi$ is strictly
increasing on $[\F(0),\infty)$, where $\F(0)$ is the largest root of
$\ps(\th)=0$. The right-inverse function of $\ps$ is denoted by
$\F\dvtx[0,\infty)\to[\F(0),\infty)$.

An important generalization is to replace the process $S$
in \eqref{SNmod} by a general subordinator
[a nondecreasing L\'{e}vy process, with L\'evy measure $\nu(\td
x), x \in\R_+$, which may have infinite mass]. Under this model,
the ``small fluctuations'' can arise either continuously, due to
the Brownian motion, or due to the infinite jump-activity.

Taking $S$ to be a pure jump-martingale with i.i.d. increments and
negative jumps with L\'{e}vy measure $\nu(\td x)$, one arrives
thus to a general integrable spectrally negative L\'{e}vy process
$X = \{X_t, t\in\mbb R_+\}$, that is, a stochastic process that
has stationary independent increments, no positive jumps and
c\`{a}dl\`{a}g paths, such that $X_t$ integrable for any $t\in\mbb R_+$,
defined on some filtered probability space $(\Omega,\FF,{\mbf F}, \P)$,
where ${\mbf F}=\{\FF_t\}_{t\in\mbb R_+}$ is the natural
filtration generated by $X$ satisfying the usual conditions of
right-continuity and
completeness; see Bertoin \cite{bert96}, Kyprianou \cite{Kbook},
Sato \cite{Sato}. The assumption that $X_t-X_0$ has finite mean for any
fixed $t>0$
is equivalent to the requirement that
the L\'{e}vy measure $\nu$ satisfies the integrability condition
%
\begin{equation}
\label{eq:nu1i} \nu_{1,\infty}:= \int_{[1,\infty)} x \nu(\td x) <
\infty.
\end{equation}
{To avoid degeneracies, the case that $X$
has monotone paths is excluded.}
The (possibly random) initial value $X_0$ is assumed to be nonnegative.
Conditioning the probability measure $\P$ on the value of $X_0$ gives
rise to
the family of probability measures $\{\P_x, x\in\mbb R_+\}$ that
satisfy $\P_x[X_0=x]=1$.

An alternative characterization of spectrally negative L\'{e}vy
processes is via the ``$q$-harmonic homogeneous scale function''
$W^{(q)}$, a nondecreasing function defined on the real line that
is 0 on $(-\infty,0)$, continuous on $\mathbb R_+$, with Laplace
transform $\mc L W^{(q)}$ given by
%
\begin{equation}
\label{Wq} \mc L W^{(q)}(\theta) = \bigl(\ps(\th) - q
\bigr)^{-1},\qquad \th > \F(q).
\end{equation}

Despite the diversity of possible path behaviors displayed by
spectrally negative L\'evy processes,
a wide variety of results may be elegantly expressed in a unifying
manner via the homogeneous
scale function
$W^{(q)}$, bypassing thus ``probabilistic complexity'' via
unified analytic methods. This paper further illustrates this
aspect, by unveiling the way the scale function intervenes in a
quite complex control problem.

\textit{De Finetti's dividend problem.} Under the assumption that the
increments of the surplus process have positive mean, the L\'{e}vy
risk model has the unrealistic property that it converges to
infinity with probability one.

In answer to this
objection, De Finetti \cite{Fin} introduced the risk process with
dividends
%
\begin{equation}
\label{eq:U} U^\pi_t = X_t -
D^\pi_t,\qquad t \ge0,
\end{equation}
where $\pi$ is an
``admissible'' dividend control policy and $D^\pi_t$ denotes the
cumulative amount of dividends that has been transferred to a
beneficiary up to time $t$, and where $U^\pi_{0-}=X_0\ge0$ is the
initial capital.

Writing $\tau^\pi=\inf\{t\in\mbb R_+\dvtx U^\pi_t < 0\}$ for the time
at which ruin occurs, the objective is to maximize the expected
cumulative dividend payments until the time of ruin
\[
v_*(x):=\sup_{\pi\in\Pi} \E_x \biggl[ \int
_{[0, {\tau^\pi})} \mathrm{e}^{-qt}\,\td D^\pi_t
\biggr],
\]
with
$\E_x[\cdot]=\E[\cdot|X_0=x]$ and where $\Pi$ denotes the set of all
admissible strategies, and $q>0$ is the discount rate.

Note that ruin may be either exogeneous or endogeneous (i.e.,
caused by a claim or by a dividend payment). A dividend strategy
is admissible if ruin is always exogeneous, or more precisely, an
admissible dividend strategy $D^\pi=\{D^\pi_t, t\in\mbb R_+\}$ is a
right-continuous $\mbf
F$-adapted stochastic process that will satisfy that, at any time
preceding the epoch of ruin, a dividend payment
is smaller than the size of the
available reserves, that is, for any $t\leq\tau^\pi$,
%
\begin{equation}
\label{exoruinliq} \cases{ \Delta D^\pi_t:=
D^\pi_{t}-D^\pi_{t-} \leq
\bigl(X_t - D^{\pi}_{t-} \bigr)\vee0, $\mbox{ and}$& \vspace*{2pt}
\cr
D^{\pi(c)}_t - D^{\pi(c)}_u
\leq p(t-u),  \qquad \mbox{$\forall u\in[0,t)$ if $\nu_{0,1}<
\infty$,}}
\end{equation}
where $D^{\pi(c)}$ denotes the continuous part of $D^{\pi}$,
$\nu_{0,1}:=\int_{(0,1)} x\nu(\td x)$
and $p:=\eta+ \nu_{0,1} + \nu_{1,\infty}$.
In the second line in \eqref{exoruinliq} it is stated that if
the jump-part of $X$ is of bounded variation, it is not admissible to pay
dividends at a rate larger than the premium rate $p$ at any time
$t$ that there are no reserves (i.e., $U^\pi_t=0$), as this would
lead to immediate ruin.

\textit{Single barrier policies.} Recall first the simplest case when
there are no transaction costs.
One possible dividends distribution policy is the ``barrier
policy'' $\pi_b$ of transferring all surpluses above a given level
$b$, which results in the
value
\[
v_b(x): =v_{\pi_b}(x) = \E_x \biggl[\int
_{[0, {\tau_b})} \mathrm {e}^{-qt}\,\td D_t^b
\biggr] = \frac{W^{(q)}(x)}{W^{(q)\prime}(b)},\qquad x\in[0,b]
\]
and $v_b(x)=x-b+v_b(b)$ for $x>b$, where $\tau_b=\inf\{t\ge0\dvtx X_t<
D^b_t\}$, and
$D^b = D^{\pi_b}$ is a local time-type strategy, given
explicitly in terms of $X$ by $D_{0-}^b=0$ and
%
\begin{equation}
\label{Deebee} D_t^b = \sup_{s\leq t}
(X_s - b )^+,\qquad t\in\mbb R_+,
\end{equation}
with $x^+=\max\{x,0\}$. As this equation shows, a nonzero optimal barrier
must be an inflection point of the scale function if the latter
is smooth.

\textit{Multiple bands policies.} However, single barrier
strategies might not be optimal; cf. Gerber \cite{Ger69,Ger72}. The
optimal strategy
may be a ``multi-bands strategy,''
involving several ``continuation bands'' $[a_i,b_i), i=0,1,\ldots,$
with upper
reflecting boundaries $b_i$, separated by ``lump-sum dividend taking bands''
$[b_i,a_{i+1}), i=0,1,\ldots, $ of
jumping to the next reflecting barrier below $b_i$, by paying all
the excess as a lump-sum payment; see also Hallin \cite{Hal}, who
formulated a system
of time dependent
integro-differential equations associated to
multi-bands policies. Azcue and Muler \cite{AM}
established the optimality of multi-bands strategies under the Cram\'
{e}r--Lundberg model
in the presence of proportional and excess-of-loss reinsurance,
adopting a viscosity approach.
A direct approach was developed in Schmidli \cite{Schmbook} where
a recursive algorithm was provided to find, in terms of solutions to
certain integro-differential
equations, the value function of the optimal dividend problem
under the Cram\'{e}r--Lundberg model in the absence of a
penalty. Recently, Albrecher and Thonhauser \cite{AT08} proved the optimality
of bands strategies, in the case that the reserves attract a fixed
interest rate.

Gerber showed also
that for exponential claims (and with no constraints
on the dividends rate), the optimal policy involved only one
barrier (and one continuation band); however, constructing examples
where more than
one band was necessary remained an open problem for a long time.

\textit{Optimality conditions for single barrier strategies.} The interest
in bands strategies was reawakened by Azcue and Muler \cite{AM},
who produced the first example (with Gamma claims) in
which a single constant barrier is not optimal. Let
%
\begin{equation}
\label{bstar} b^* = \sup \bigl\{b>0\dvtx W^{(q)\prime}(b) \leq
W^{(q)\prime}(x)\mbox{ for all $x$} \bigr\}
\end{equation}
denote the last global minimum of the derivative of the $q$-scale function.

Avram et al. \cite{APP} showed that
%
\begin{equation}
\label{APP}(\Gamma v_{b^*} - qv_{b^*}) (x)\leq0\qquad
\mbox{for all $x>b^*$},
\end{equation}
where $\Gamma$ denotes the infinitesimal generator of $X$, is a
sufficient optimality condition for the single barrier strategy under
a general
spectrally negative L\'{e}vy model. In fact, conditions
\eqref{bstar}--\eqref{APP} is both necessary and sufficient, as
follows by examining the variational inequality characterizing the
problem; see Loeffen \cite{Loeffen1}, Lemmas 1, 2.

A simpler sufficient condition for the optimality of single band
policies was
obtained by Loeffen \cite{Loeffen1,Loeffen2} (with and without
transaction costs), who showed that it is enough to check
that the last local minimum of the $q$-scale function is also a
global minimum. Even more direct optimality conditions in terms of the
L\'{e}vy measure $\nu$
were provided by
Kyprianou et al. \cite{krs}, and Loeffen and Renaud \cite{LR}, who
showed, respectively, that log-convexity of the density and of the
survival functions suffice (the second condition is more general). Note
that the second result allowed also for an affine penalty
function with slope less than unity, and that both results imply
complete monotonicity of the L\'{e}vy density, and constitute therefore
powerful generalizations of Gerber's unicity result \cite{Ger69,Ger72}.

It turns out that $b^*$ in \eqref{bstar} is always the right endpoint
of the first continuation band. As already demonstrated in the rather
terse example in Azcue and Muler (\cite{AM}, page 274), left and right
endpoints of subsequent bands can in principle be determined
recursively (the former by ensuring the ``smoothness'' of the value
function, and the latter similarly with $b^*$, by selecting last global
maxima of updated value
functions, adjusted by using the values of previous bands as stopping
penalties). A characterization of points of nondifferentiability was
provided in
Schmidli \cite{Schmbook}. However, an explicit smoothness
condition \eqref{eq:astar2} in terms of scale functions seems not to
have been reported previously.

Quite paradoxically, it is possible that beyond the lump-sum dividend
taking band following the first continuation band, waiting for higher
barriers $b_i, i \geq2$, may become again optimal. The level $a_2$
where the second continuation band starts may be determined by
examining the family of functions $G_2^{(a)}(b)$ defined in \eqref
{eq:astar2}, which are computed from a second Gerber--Shiu function,
which uses the first value functions as stopping penalties, and so on, leading
ultimately to all the optimal band levels; see Section~\ref{sec:multib}.

\textit{Fixed transaction costs.} It is interesting to consider also
the effect of adding fixed transaction cost $K>0$ that are not
transferred to the beneficiaries when dividends are being paid.
The objective of the beneficiaries becomes then to maximize
$v_{\pi,K}(x)$, that is, $v_*(x) = \sup_{\pi\in\Pi} v_{\pi,K}(x)$ with
\[
v_{\pi,K}(x) = \E_x \biggl[\int_{[0, {\tau^\pi})}
\mathrm{e}^{-qt}\,\td D^\pi_t - K \int
_{[0, {\tau^\pi})}\mathrm{e}^{-qt}\,\td N^\pi_t
\biggr],
\]
where $N^\pi=\{N^\pi_t, t\in\mbb R_+\}$ is the stochastic process
that counts
the number of jumps of $D^\pi$ in the interval $[0,t]$,
%
\begin{equation}
\label{eq:count} N^\pi_t = \#\bigl\{s\in[0,t]\dvtx\Delta
D^\pi_s>0\bigr\},\qquad t\in\mbb R_+.
\end{equation}
The introduction of a fixed transaction cost $K>0$ has the usual
effect of
changing the optimal reflection boundaries $b$
into strips $[b_{-}, b_{+}]$, so that when $U_t=b_+$,
a~lump-sum dividend $b_+-b_{-}$ is
paid, and the reserves process is
diminished to the lower ``entrance'' point $b_{-}$. To emphasize
this disappearance of reflection barriers, the term \emph{band} will be
used throughout
when $K >0$, and also when more than one barrier is present.

The typical optimal dividend strategy consists of ``lump-sum
payments'' (see, e.g., Alvarez and Virtanen \cite{AlvVir}
and Thonhauser and Albrecher \cite{TA2}), with $\pi$ of the form $
\pi
=\{(J_k,
T_k), k\in\mbb N\}$, where $0\leq T_1\leq T_2 \leq\cdots$ is an
increasing sequence of $\mbf F$-stopping times representing the
times at which dividend payments are made, and $J_i \geq K$ is a
sequence of positive $\mc F_{T_i}$-measurable random variables
representing the sizes of the dividend payments. Then
\[
D^\pi_t = \sum_{k=1}^{N^{\pi}_t}
J_k,
\]
where $N^\pi_t = \#\{k\dvtx T_k\leq t\}$ is the number of times that
dividends have been paid by time $t$.

For single bands policies for example, the dividend distribution
consists of the fixed amount $J_i=b_{+}-b_{-}$.

\textit{Balancing dividends and ruin penalties}.
Several alternative objectives have been proposed recently, involving
final penalties $w(x)$ at ruin (see Dickson and Waters~\cite{DW},
Gerber et al. \cite{GerberLY} and Zajic \cite{Zajic}), or continuous
payoffs until ruin; see Albrecher and Thonhauser \cite{AT07}, Cai et
al. \cite{cai2009expectation}.
For example, the case where the insurance company is bailed out by the
beneficiaries every time that there is a shortfall in the reserves was
investigated in Avram et al. \cite{APP} and in Kulenko \&
Schmidli~\cite{kulschm}.
This paper continues the investigation of the impact of a general
penalty and fixed transaction costs on the optimal dividends policy.
The considered objective is to maximize the expected cumulative
discounted dividend payments until the moment of ruin less the penalty, which
is an increasing function of the shortfall at the moment of ruin, by
controlling the timing and
size of dividend payments.
This problem is phrased as an optimal control problem, which will be solved
by constructing explicitly a solution of the associated
Hamilton--Jacobi--Bellman (HJB) equation, in terms of scale functions
of the L\'{e}vy process~$X$.

\textit{Stochastic solutions}. Given results concerning the
smoothness of scale functions (see, e.g., Kyprianou et al. \cite{KypScale}),
it is not to be expected that the candidate value-function is a
classical solution
of the HJB equation. In fact, it will turn out that
the candidate value function is
continuous but not $C^1$ on $\mbb R_{+}\setminus\{0\}$ if $X$ has
bounded variation, and is $C^1$ but not $C^2$ on $\mbb R_{+}\setminus
\{
0\}$, if $X$ has unbounded variation.
To verify optimality of the candidate optimal value-function under weak
regularity conditions,
a probabilistic approach is adopted in this paper. It is established
that the value-function is the unique
\emph{stochastic solution} of the HJB equation corresponding to the
optimal control problem
under consideration. The notion of stochastic solution may informally
be considered as a probabilistic counterpart of the analytical notion
of viscosity solution: While viscosity sub- and supersolutions are
defined in terms of pointwise approximations by smooth solutions to the
variational inequalities associated to the HJB equation, stochastic
super- and subsolutions are phrased in terms of super- and
submartingale properties of related stochastic processes.
The version of the notion of stochastic solution deployed here
is an adaptation of Stroock and Varadhan's \cite{SV} classical notion,
which was originally introduced
in the setting of linear parabolic PDEs, to the current setting; see
Definition~\ref{def:sss}. A~stochastic version of Perron's method using
the stochastic solution concept was recently developed in Bayraktar and
S\^{\i}rbu \cite{BS} for the case of linear parabolic PDEs.

The viscosity solution method is a classical approach that has been
used extensively in the study of existence and uniqueness of solutions
to HJB equations; cf. Bardi and Capuzzo-Dolcetta \cite{BC} and Fleming
and Soner \cite{FM} for general treatments. The HJB equation \eqref
{eq:HJB} corresponding to the stochastic control problem considered in
the current paper is a nonlinear integro-differential equation with
constant coefficients and with a gradient constraint, which is of first
or second order depending on whether or not a Gaussian component is
present in the dynamics of $X$.
Due to the negative jumps of $X$ and the boundary condition on the
negative half-axis (the specified penalty at the epoch of ruin), one is
led to the notion of constraint viscosity solutions which, in the
context of different optimization problems, has been developed for
first order equations by Sayah \cite{Sa1} and Soner \cite{So1}, and
for second order equations in Alvarez and Tourin \cite{AvT}, Benth et
al. \cite{Benth} and Pham \cite{Pham}.
In, for example, Azcue and Muler \cite{AM,AM2} and Albrecher and
Thonhauser \cite{AT08},
dividend optimization problems are studied under the Cram\'
{e}r--Lundberg model
using the viscosity solution method.

By deploying probabilistic tools from among others martingale theory,
analogues are derived of key results from viscosity solution theory. In
particular, existence and uniqueness of a stochastic solution to the
HJB equation is shown (Theorem~\ref{thm:repg2}),
where the uniqueness is established deploying a {comparison principle}
(Proposition~\ref{prop:stcomp}).
A {(local) verification theorem} (Theorem~\ref{cor:repg}) is derived
as tool for verifying optimality of a constructed candidate value-function,
as direct consequence of a dual representation of the value function as
pointwise minimum of stochastic supersolutions (Proposition~\ref{thm:repg}).

\textit{Gerber--Shiu functions.}
A key point in the presented approach is the decomposition of the
candidate value
function preceding and within a continuation band $[a,b]$
%
\begin{equation}
\label{eq:dec} v_{a,b}(x) = \cases{ f(x), \quad& $x < a,$\vspace*{2pt}
\cr
F(x) + W^{(q)}(x)G(a,b),\quad &$x\in[a,b]$,}
\end{equation}
into a \textit{nonhomogeneous} solution $F(x)$, which will be called a
\emph{Gerber--Shiu function}, and the product of the homogeneous scale
function $W^{(q)}(x)$ and a ``barrier-influence'' function $G(a,b)$
defined in \eqref{eq:Gw}, which needs to be maximized at $b$ and be
smooth at $a$.

Note that the function $G$ in the decomposition \eqref{eq:dec} is only
determined up to a constant, but becomes fixed once $F$ has been
selected; see \eqref{eq:Gfabb}.

To ensure smoothness at $a$, it seems then natural to use a ``smooth
Gerber--Shiu function'' $F_f(x)$
associated to a given penalty $f(x), x \in(-\infty,a)$. Informally,
$F_f(x)$ is the ``smooth
nonhomogeneous solution'' of the Dirichlet problem on $\{x \geq a\}$
with boundary condition $f(x), x \in(-\infty,a)$. More
precisely, it is defined in
Definitions \ref{d:GS} and \ref{Def:Fw} in Section~\ref{sec:RLevy}
by subtracting
a multiple of
the homogeneous scale function $W^{(q)}(x)$ out of the
solutions of
either the two-sided, or the reflected exit problem,
such that the remaining part is \emph{continuous} on $\mbb R$ if $f$
is continuous, and \emph{continuously differentiable} on $\mbb R$ if
$f$ is continuously
differentiable on $\mbb R_-$ and $X$ has unbounded variation. This
results in the explicit formula \eqref{eq:Fwr}.

For exponential penalties $w(x)=\mathrm{e}^{x v}$, the Gerber--Shiu
function takes a simple form \eqref{eq:Zv}, which may be used also as a
generating function for the expected payoffs associated to
polynomial penalties $x^k, k=0,1,\ldots.$

Decomposition \eqref{eq:dec} with $F_f(x)$ chosen to fit the imposed
penalty $f(x)=w(x)$ already determines the value function on the first
continuation band (and the value function in the lump-sum dividend
taking bands surrounding it); see
Proposition~\ref{prop:lp2} and Theorem~\ref{thmos}. It also yields a
{necessary} and {sufficient}
criterion for optimality of two-dividend barrier policies with one
barrier at zero, which is analogous to \eqref{APP}; %
see Theorem~\ref{thm:cm3}.

\textit{Contents}.
The remainder of the paper is organized as
follows. Sections \ref{sec:prob} and \ref{sec:DP} are devoted to the
formulation of
the dividend-penalty and the corresponding HJB equation.
In Section~\ref{sec:mart} the definition of stochastic solution is
given in this context, and a
verification result is established.
Section~\ref{sec:RLevy} is concerned with Gerber--Shiu functions, and
Sections \ref{aux} and \ref{aux2} are devoted to single and two-band
strategies. Section~\ref{sec:paste} is devoted to a key auxiliary
result (Lemma~\ref{lem:paste}).
Conditions for optimalty of single and two-band strategies and a
construction of the candidate value-function in terms of scale
functions are given in Sections \ref{ssec:singo}, \ref{sec:doublop}
and \ref{sec:multib}.
The optimal value function is shown to be the unique stochastic solution
of the HJB equation in Section~\ref{sec:ex}. Some examples are analyzed
in Section~\ref{sec:exa}.
Some of the proofs are deferred to the \hyperref[app]{Appendix}.

\section{The dividend-penalty control problem}\label{sec:prob}
Assume that the beneficiaries control the timing and size
of dividend payments made by the company, and are
liable to pay at the moment $\tau^\pi$ of ruin the penalty
$-w(U^\pi_{\tau^\pi})$, which may be used to cover (part of) the
claim that led to insolvency, where $w$ is a penalty.

%
\begin{Def}\label{def:R}
(i) For any $a\in\mbb R$, denote by $\mc R_a$
the set of c\`adl\`ag\setcounter{footnote}{2}\footnote{c\`adl\`ag${} = {}$right-continuous with
left-limits.} functions
$w\dvtx(-\infty,a]\to\mbb R$ that are left-continuous at $a$,
admit a finite first
left-derivative $w_-'(a)$ at $a$ and satisfy the
integrability condition
%
\begin{equation}
\label{eq:cw2} \sup_{y>1} \int_{[y,\infty)}
\sup_{u\in[y-1,y]}\bigl|w(a+u-z)\bigr|\nu (\td z) < \infty.
\end{equation}

(ii) A \emph{penalty} $w\dvtx\R_-\to\R_-$, with $\mbb R_-=(-\infty,0]$,
is a function from the set $\mc R_0$ that is increasing.
The collection of penalties is denoted by $\mc P$.
\end{Def}

The beneficiaries seek to maximize the sum of the expected
discounted cumulative dividend payments and an expected penalty payment
by paying out dividends according to an admissible policy.
The present value of the
penalty payment discounted at rate $q>0$, considered as function of the
level of initial reserves, is called the \emph{Gerber--Shiu penalty function}
associated to the {penalty} $w$,
and is given by
\[
{\mathcal W}_w^\pi(x):= \E_x \bigl[
\mathrm{e}^{-q\tau^\pi} w \bigl(U^\pi_{\tau^\pi} \bigr) \bigr],\qquad
x\in\mbb R_+.
\]
For any penalty $w\in\mc P$, it holds that, for any level of initial
capital $x\in\mbb R_+$, $\mc W_w^\pi(x)$ is bounded uniformly over
$\pi\in\Pi$; see Lemma~\ref{lem:est}.

The objective of the beneficiaries of the insurance company is
described by
the following stochastic
control problem:
%
\begin{equation}\quad
\label{optdiv} v_*(x) = \sup_{\pi\in\Pi} v_\pi(x),\qquad
v_\pi(x):= {\mathcal W}_w^\pi(x)+
\E_x \biggl[\int_{[0,{\tau^\pi})} \mathrm{e}^{-qt}
\mu _K(\td t) \biggr],
\end{equation}
for $x\in\mbb R_+$,
where $\Pi$ denotes the set of admissible dividend policies $\pi$
and $\mu_K$ is the (signed) random measure on $(\mbb R_+,
\mc B(\mbb R_+))$ defined by
%
\begin{equation}
\label{eq:muk} \mu_K^\pi\bigl([0,t]\bigr) =
D^\pi_t - K N^\pi_t,
\end{equation}
with $N^\pi_t$ and $D^\pi_t$ equal to the counting process defined in
\eqref{eq:count} and the cumulative
amount of dividends that has been paid out by time $t$, respectively.
It is assumed
throughout that $w$ is a penalty
($w\in\mc P$) and that there is positive net income, $\eta:=
\mathbb E[X_1]>0$.
A solution to the stochastic control
problem in \eqref{optdiv} consists of a pair $(u,\pi^*)$ of a
function $u\dvtx\mbb R_+\to\mbb R$ and a policy $\pi^*\in\Pi$ satisfying
$v_*(x)= u(x) = v_{\pi^*}(x)$ for all $x\in\mbb R_+$.

\section{Dynamic programming and HJB equation}\label{sec:DP}
The analysis of the stochastic optimal control problem \eqref{optdiv}
starts from the observation that
the value function~$v_*$ satisfies a dynamic programming equation.

\begin{Prop}\label{prop:mart}
\textup{(i)} Extending $v_*$ to the negative half-axis by $v_*(x)=w(x)$
for $x<0$, we have for any
$\tau\in\mc T$, the set of $\mbf F$-stopping times, $v_*(x) = \sup_{\pi
\in\Pi} v_{\pi,\tau}(x)$ where
%
\begin{equation}
\label{eq:vstar} v_{\pi,\tau}(x):= \E_x \biggl[
\mathrm{e}^{-q(\tau\wedge\tau
^\pi)} v_*\bigl(U^\pi_{\tau
\wedge\tau^\pi}\bigr) + \int
_{[0,\tau\wedge\tau^\pi]} \mathrm{e}^{-qs} \mu_K^\pi(
\td s) \biggr]. %
\end{equation}

\textup{(ii)} For any fixed $\pi\in\Pi$,
the process $V^\pi=\{V^\pi_t,t\in\mbb R_+\}$ given by
%
\begin{equation}
\label{eq:Vpi} V^\pi_t:= \mathrm{e}^{-q(t\wedge\tau^\pi)}v_*
\bigl(U^\pi_{t\wedge
\tau^\pi
}\bigr) + \int_{[0,t\wedge\tau^\pi]}
\mathrm{e}^{-qs} \mu^\pi_{K}(\td s)
\end{equation}
is an $\mbf F$-supermartingale.
\end{Prop}

\begin{Rem}
Note that the integration domains $[0,\tau\wedge\tau^\pi]$ and
$[0,t\wedge\tau^\pi]$
in \eqref{eq:vstar} and \eqref{eq:Vpi}
are consistent with the domain $[0,\tau^\pi)$ in \eqref{optdiv}
as $\mu_K(\{\tau^\pi\})$ is equal to $0$ for any policy $\pi\in\Pi$.
\end{Rem}

The proof of Proposition~\ref{prop:mart}(i)
follows by straightforward adaptation of classical arguments
(see, e.g., \cite{AM}, pages 276--277), while that of
Proposition~\ref
{prop:mart}(ii)
is deferred to  Appendix \ref{app:dp}.

The next step is to identify the HJB equation in the current setting.
As the beneficiaries may decide to pay out part of the reserves
immediately as lump-sum dividend
the value function $v_*$ satisfies in addition to the dynamic
programming equation
the following gradient condition (see Lemma~\ref{lem:est}):
%
\begin{equation}
\label{eq:vlbub} v_*(x) - v_*(y) \geq (x-y - K) \qquad\mbox{for al $x,y> 0$ with $x >
y$, }
\end{equation}
or equivalently,
%
\begin{eqnarray}\label{qgamma}
 \mathtt d_{v_*}(x)&\ge&1 \qquad\mbox{for all $x>0$, with for any
function $g\dvtx\mbb R\to\mbb R$},
%
\nonumber
\\[-8pt]
\\[-8pt]
\nonumber
 {\mathtt
d}_g(x) &=& \inf_{y\in(0,x)}\frac{g(x) - g(x-y) + K}{y},\qquad x>0.
\end{eqnarray}
Note that in the case $K=0$ and when $v_*|_{\mbb R_+\setminus\{0\}}$
is in $C^1(\mbb R_+\setminus\{0\})$ the gradient constraint in \eqref
{eq:vlbub}
is equivalent to the condition
\[
v_*'(x)\ge1 \qquad\mbox{for all $x>0$}.
\]

Rather than to pay out dividends immediately, the beneficiaries may decide
to postpone such payments to a future epoch. Provided the value
function $v_*$ were sufficiently regular,
it would hold at level $x$ of the reserves that $\E_x[\mathrm
{e}^{-q(t\wedge
T^-_0)} v_*(X_{t\wedge T^-_0})] = v_*(x) + t(\Gamma v_*(x) - q v_*(x))
+ o(t)$ for $t\searrow0$,
where $T_0^-=\inf\{t\ge0\dvtx X_t<0\}$, and $\Gamma$ denotes the
infinitesimal generator of the Feller semi-group of
$X$ which acts
on $f\in C^2_c(\mbb R_+)$ as (cf. Sato \cite{Sato}, Theorem~31.5)
%
\begin{equation}\qquad
\label{eq:gamma11} \Gamma f(x) = \frac{\s^2}{2}f''(x)
+ \eta f'(x) + \int_{\mbb
R_+\setminus\{0\}}\bigl[f(x-y) - f(x) +
yf'(x)\bigr]\nu(\td y),
\end{equation}
for $x\in\mbb R_+$, where $f'$ denotes the derivative of $f$
and $\eta=\psi'(0)$. Heuristically, this suggests that $v_*$ satisfies
$\Gamma v_*(x) - q v_*(x)\leq0$ at any $x>0$, and that it is not
optimal to postpone a dividend payment
at level $x$ in case $\Gamma v_*(x) - q v_*(x)< 0$.

As far as the boundary condition at $x=0$ is concerned, it follows
from \eqref{optdiv} that $v_*(0)=w(0)$
if and only if ruin is immediate with zero initial capital
(i.e., $\tau^\pi=0$ $\P_0$ a.s.),
which is precisely the case if $X$ has paths of unbounded variation.
Thus the boundary condition at $x=0$ is imposed precisely
if the Gaussian coefficient $\sigma^2$ is strictly positive or
the L\'{e}vy measure $\nu$ does not finitely integrate $x$ around 0
($\nu_{0,1} = \infty$). In particular,
in the case of the Cram\'{e}r--Lundberg model or when $X$
has paths of finite variation, $v_*(0)$ is in general different from $w(0)$.

By the above discussion, one is led to the following form of the HJB
equation associated to the optimal control problem \eqref{optdiv},
expressed in a unified manner for general cost $K\ge0$:
%
\begin{equation}
\label{eq:HJB} \max \bigl\{\Gamma g(x) - q g(x), 1 - {\tt d}_g(x)
\bigr\} = 0, \qquad x>0,
\end{equation}
subject to the boundary condition
%
\begin{equation}
\label{eq:bc} \cases{ g(x) = w(x), &\quad $\mbox{for all $x<0$, and}$ \vspace*{2pt}
\cr
g(0) = w(0), & \quad $\mbox{in the case $\bigl\{\s^2>0$ or $
\nu_{0,1} = \infty \bigr\}$}$,}
\end{equation}
where the function ${\tt d}_g$ is defined in \eqref{qgamma}.
%
\subsection{Properties of the value function}

For later reference a number of properties of the value function are
collected below.

%
\begin{Lemma}\label{lem:est} \textup{(i)} The function $x\mapsto v_*(x)$ is
continuous on $\mbb R_+$, and $v_*$ satisfies equation \eqref{eq:vlbub}.\vspace*{-6pt}
\begin{longlist}[(iii)]
\item[(ii)] For any $q>0$, $x\in\mbb R_+$ and $w\in\mc P$, there exists a
$C\in
\mbb R_{+}\setminus\{0\}$
such that the following bound holds true:
\begin{eqnarray*}
&&\E_x \biggl[\sup_{t\in\mbb R_+, \pi\in\Pi} \biggl\{\mathrm
{e}^{-qt} U^\pi_{t} \mbf 1_{\{t<\tau^\pi\}} +
\int_0^t \mathrm{e}^{-qs}\,\td
D^\pi_s + \int_0^t
\mathrm{e}^{-qs}(\ovl X_s - \unl X_s)\,\td s
\biggr\} \biggr]
\nonumber
\\[-8pt]
\\[-8pt]
\nonumber
&&\qquad{}
+ \sup_{y\in\mbb
R_+}\sup_{\pi\in\Pi}\E_y
\bigl[\mathrm{e}^{-q\tau}\bigl|w\bigl(U^\pi _\tau\bigr)\bigr|
\bigr] < C,
\end{eqnarray*}
with $\ovl X_t = \sup_{s\leq t}X_s$ and
$\unl X_t = \inf_{s\leq t}X_s$ denoting the supremum and infimum of
$X_s$ over the $s\in[0,t]$.

\item[(iii)] $v_*$ is dominated by an affine function: for any $x\in\mbb R_+$,
$v_*(0) - K \leq
v_*(x) - x\leq\frac{1}{\Phi(q)}$, and
the process $V^\pi=\{V^{\pi}_t, t\in\mbb R_+\}$
defined in \eqref{eq:Vpi} is a uniformly integrable (UI)
${\mbf F}$-supermartingale.
\end{longlist}
\end{Lemma}

The proof of part (i) is deferred to  Appendix \ref{sec:lem:est}.

\begin{pf*}{Proof of Lemma~\ref{lem:est}(ii)}
The following bounds hold true:
%
\begin{equation}
\label{eq:sas} \sup_{t\in\mbb R_+} \mathrm{e}^{-qt}U^\pi_{t}
\mbf1_{\{t<\tau^\pi
\}} \leq \sup_{t\in\mbb R_+} \mathrm{e}^{-qt}
X_t \leq\sup_{t\in\mbb
R_+}\int_t^\infty
q\mathrm{e}^{-qs}\ovl X_s\,\td s.
\end{equation}
Since the running supremum $\ovl X_{\mbf e_q}$ at an independent
exponential random time $\mbf e_q$ with mean $q^{-1}$ under $\P_0$
follows an exponential distribution
with parameter $\Phi(q)$ (e.g., Bertoin \cite{bert96}, Corollary
VII.2), the expectation under $\P_x$ of the
expression on the
right-hand side of \eqref{eq:sas} is bounded by $x + 1/\Phi(q)$.

The compensation formula applied to the Poisson point process
$(\Delta X_t, t\in\mbb R_+)$, the monotonicity of $w$ and the
fact that $w(0)$ is nonpositive yield that the following
inequalities holds true, for any $x\in\mbb R_+$:
\begin{eqnarray*}\qquad
\label{eq:wu}
\E_x \bigl[\mathrm{e}^{-q\tau^\pi}w
\bigl(U^\pi_{\tau^\pi}\bigr) \bigr] &\geq& w(-1) +
\E_x \bigl[\mathrm{e}^{-q\tau^\pi}w\bigl(U^\pi_{\tau^\pi}
\bigr)\mbf 1_{\{U^\pi_{\tau^\pi} < -1\}} \bigr]
\\
&=& w(-1) + \int_0^\infty \int
_0^\infty w(y-z)\mbf1_{\{y-z<-1\}}\nu(\td z)
\tilde R_x^q(\td y),
\end{eqnarray*}
where $\tilde R_x^q(\td y)$ denote the $q$-potential measure of
$U^\pi$ under $\P_x$,
$\tilde R_x^q(\td y) = \int_0^\infty\mathrm{e}^{-q
t}\P_x(U^\pi_t\in\td y, t<\tau^\pi)$. The right-hand side of \eqref{eq:wu} is
bounded below, as $w$ satisfies the integrability condition \eqref{eq:cw2}
(as $w\in\mc P$).
\end{pf*}

\begin{pf*}{Proof of Lemma~\ref{lem:est}(iii)}
In the case $K=0$
integration by parts, the nonnegativity of $w$
and condition \eqref{exoruinliq} of ``no exogenous ruin''
imply that
\begin{eqnarray*}
v_{\pi}(x) &\leq& \E_x \biggl[\int_{[0,{\tau^\pi})}
\mathrm {e}^{-qt}\,\td D^\pi_t \biggr] =
\E_x \biggl[\int_0^{\tau^\pi}q
\mathrm{e}^{-qs} D^\pi_{s}\,\td s +
\mathrm{e}^{-q\tau^\pi}D^\pi_{\tau^\pi} \biggr]
\\
&\leq& \E_x \biggl[\int_0^{\tau^\pi}q
\mathrm{e}^{-qs} X_{s}\,\td s + \mathrm{e}^{-q\tau^\pi}X_{\tau^\pi-}
\biggr] \leq\E_x \biggl[\int_0^\infty
q \mathrm{e}^{-qs} \ovl X_s\,\td s \biggr],
\end{eqnarray*}
which is equal to $x + \frac{1}{\Phi(q)}$
since, as noted before, $\ovl X_{\mbf e_q}\sim\mrm{Exp}(\Phi(q))$ under
$\P_0$.
In the case $K>0$, then
the above bound remains valid since the value $v_*(x)$
decreases if the transaction cost $K$ increases.
The lower bound for the value-function
follows from part (i) (with $x=0$).
The uniform integrability of $V^\pi$ is a consequence of the
fact that $V^\pi$ is dominated by an integrable random variable, in view
of the bounds in parts (ii).
\end{pf*}

\subsection{Generator and boundary condition}\label{sec:gener}
From the HJB equation \eqref{eq:HJB} one would expect that, on any interval
$I$ on which the restriction $v_*|_I$ has unit derivative,
the function $\Gamma v_* - qv_*$ is nonpositive. Below this function
is expressed explicitly in terms of the characteristic triplet of $X$.
More generally, in the next result the form is specified of the
generator applied to
the functions %
$\WT\ell^w_{a,b}\dvtx\mbb R\to\mbb R$, $a,b\in\mbb R$, given by
%
\begin{eqnarray}
\WT\ell^w_{a,b}(z) = \cases{\ell_{a,b}(z), & \quad $z\ge
a$,\vspace*{2pt}
\cr
w(z), &\quad  $z<a$,}\nonumber\\
\eqntext{\mbox{with $\ell_{a,b}\dvtx[a,
\infty)\to\R$: $\ell_{a,b}(x)= b(x - a) + w(a)$},}
\end{eqnarray}
where $w\dvtx(-\infty,a]\to\mbb R$ is a Borel-function satisfying the
integrability condition
%
\begin{equation}
\label{aint} \forall x>a\dvtx\int_{(x-a,\infty)}\bigl|w(x-z)\bigr|\nu(\td z) <
\infty.
\end{equation}
For any such function $w$ and any $a\in\mbb R$, the operator
$_a {\Gamma}_{\infty} ^w\dvtx C^2([a,\infty))\to D((a,\infty))$ is defined
as follows: for $x>a$,
%
\begin{eqnarray}\label{eq:gamma2}
_a {\Gamma}_{\infty} ^w f(x)&=&
\frac{\s^2}{2}f''(x) + \bigl(\eta+ \ovl
\nu_1(x-a)\bigr) f'(x) - \bigl(q+\ovl\nu(x-a)\bigr)f(x)
\nonumber
\\
&&{}+ \int_{(0,x-a]} \bigl[f(x-y) - f(x) + f'(x)y
\bigr]\nu(\td y) \\
&&{}+ \int_{(x-a,\infty)}w(x-y)\nu(\td y),
\nonumber
\end{eqnarray}
where $\ovl\nu(x)=\nu((x,\infty))$ and
$\ovl\nu_1(x)=\int_{(x,\infty)}y\nu(\td y)$.
It follows by comparison with form \eqref{eq:gamma11} of the
infinitesimal generator $\Gamma$
that for any $f\in C^2_c(\mbb R)$ with $f|_{(-\infty,a]}=w$ it holds
$(\Gamma f - qf)(x) ={} _a {\Gamma}_{\infty} ^w g(x)$ for $x > a$ with
$g=f|_{[a,\infty)}$.
The form of the generator applied to $\ell_{a,b}$ is given in the
following result:

%
\begin{Lemma}\label{lem:smp}
Let $a,b\in\mbb R$ and let $w$ be any Borel function satisfying
integrability condition \eqref{aint}.
\textup{(i)} For any $x> a$, $({}_a{\Gamma}^{w}_\infty\ell_{a,b})(x)$ is
given by
%
\begin{eqnarray}\label{eq:Uy}
&&\eta\ell_{a,b}'(x) - q \ell_{a,b}(x) + \int
_{\mbb R_+\setminus\{0\}}\bigl[\WT\ell ^w_{a,b}(x-z) -
\ell_{a,b}(x) + z\ell_{a,b}'(x)\bigr]\nu(\td z)
\nonumber
\\
&&\qquad= b\eta- q\bigl(b(x-a) + w(a)\bigr) \\
&&\qquad\quad{}+ \int_{(x-a,\infty)}\bigl
\{w(x-z) - w(a) + b(z+ a -x)\bigr\}\nu(\td z).
\nonumber
\end{eqnarray}

\textup{(ii)} Suppose $({}_a{\Gamma}^{w}_\infty\ell_{a,b})(x) \leq0$ for all
$x> a$
and $\sup_{x>a}\int_{(x-a,\infty)}|w(x-z) - w(a) +
b(z+ a -x)|\nu(\td z)<\infty$.
Then $\{\mathrm{e}^{-q(t\wedge T_a^-)}\WT\ell^w_{a,b}(X_{t\wedge
T_a^-}),t\in\mbb R_+\}$ is an $\mbf F$-supermartingale.
\end{Lemma}

\begin{pf}
(i) The assertion directly follows from the form \eqref{eq:gamma2}
of the operator~${}_y{\Gamma}^{w}_\infty$.

(ii) An application of It\^o's lemma [which is justified since
$\ell_{a,b}$ is $C^2([a,\infty))$] shows that the following process is
an $\mbf F$-local martingale:
%
\begin{equation}
\mathrm{e}^{-q(t\wedge T_a\label{eq:mmart}^-)}\WT\ell ^w_{a,b}(X_{t\wedge
T_a^-})-
\int_0^{t\wedge
T_a^-}\mathrm{e}^{-qs}
{}_a{\Gamma}^{w}_\infty\ell_{a,b}(X_{s})
\,\td s.
\end{equation}
Hence the assumptions (together with the fact
$\int_0^{T^-_a}\mbf1_{\{X_s=a\}}\,\td s=0$ $\P$-a.s.)
imply the asserted supermartingale property.
\end{pf}

%
\section{Stochastic solutions of the HJB equation}\label{sec:mart}
While, as was mentioned in the \hyperref[sec:intro]{Introduction}, it is in
general not to be
expected that the HJB equation
in \eqref{eq:HJB} admits a classical solution,
it will be shown in Section~\ref{thm:repg2} that the optimal
value-function $v_*$ is the unique
\emph{stochastic solution} to the HJB equation. A real-valued function
$g$ with domain $\mbb R$ and sublinear growth,
satisfying the boundary condition~\eqref{eq:bc} and the gradient constraint
$\mathtt d_g(x)\ge1$ for all $x>0$,
will be called a stochastic solution of the HJB equation
given in \eqref{eq:HJB} if the stochastic processes
%
\begin{eqnarray}
\label{eq:ggmartI} \ovl M^{g,T_I} &:=& \bigl\{\mathrm{e}^{-q (t\wedge T_I  )}g
(X_{t\wedge
T_I} ), t\in\mbb R_+ \bigr\},
\nonumber
\\[-8pt]
\\[-8pt]
\nonumber
T_I &:=& \inf\{t\ge0\dvtx
X_t\notin I\},
\end{eqnarray}
with $\inf\varnothing=\infty$,
are $\mbf F$-martingales for any closed interval $I$ contained in
$\mathcal C_g$,
the ``no dividend region'' corresponding to the function $g$,
%
\begin{equation}
\label{eq:CD} \mc C_g:= \bigl\{x\in\mbb R_+\setminus\{0\}\dvtx{
\tt d}_g(x) > 1\bigr\},
\end{equation}
and are $\mbf F$-supermartingales for any closed interval $I$ contained
in $\mbb R_+\setminus\{0\}$.

More specifically, the notions of (local) stochastic (super-, sub-)
solutions are defined
as follows:

\begin{Def}\label{def:sss}
Let $g\dvtx\mathbb R\to\mathbb R$ be a c\`{a}dl\`{a}g function
satisfying the boundary condition \eqref{eq:bc} and the linear growth condition
%
\begin{equation}
\label{eq:gbound} \sup_{x\in\mbb R_+} \frac{|g(x)|}{x+1} < \infty.
\end{equation}
\begin{longlist}[(iii)]
\item[(i)] $g$ is a \emph{local stochastic
supersolution on the closed interval $I\subset\mbb R_+$}
of the HJB equation \eqref{eq:HJB}
if
\[
\ovl M^{g, T_I} \mbox{is a UI $\mbf F$-supermartingale and $
\mathtt d_g(x)\ge1$ for any $x\in I\setminus\{0\}$.}
\]

\item[(ii)] $g$ is called a
\emph{local stochastic
subsolution on the closed interval $I\subset\mc C_g$}
of the HJB equation \eqref{eq:HJB} if
\[
\ovl M^{g,T_I} \mbox{is a UI $\mathbf F$-submartingale.}
\]

\item[(iii)] $g$ is a \emph{stochastic supersolution} [\emph{stochastic subsolution}]
of the HJB equation
if $g$ is a local stochastic supersolution on $\mbb R_+$
[local stochastic subsolution on $I$ for all closed intervals $I\subset
\mc C_g$], respectively.

\item[(iv)] $g$ is a \emph{stochastic
solution} of the HJB equation
if $g$ is both a stochastic supersolution and a stochastic
subsolution of the HJB equation.
\end{longlist}
\end{Def}

\begin{Rem}\label{rem:stst}(i) The optimal value-function $v_*$ is
a stochastic supersolution.
This follows as a direct consequence of Lemma~\ref{lem:est}(i), (iii)
(taking $\pi$ equal
to the ``waiting strategy'' $\pi_\varnothing$ of not paying any
dividends) and Doob's Optional Stopping theorem.\vspace*{-6pt}
\begin{longlist}[(iii)]
\item[(ii)] The terms ``stochastic supersolution'' and ``stochastic
subsolution'' are justified by the fact that stochastic supersolutions
dominate stochastic subsolutions (under some regularity condition); see
Proposition~\ref{prop:stcomp}.

\item[(iii)] When $g$ is a local stochastic supersolution on a finite
partition of intervals of $\mathbb R_+$,
a global super-martingale property holds true on $\mathbb R_+$,
provided that $g$ is differentiable
at the boundaries of the intervals when $X$ has unbounded variation;
see Corollary~\ref{prop:global}.
\end{longlist}
\end{Rem}

The following global representation of the optimal value function $v_*$
in terms of the collection of stochastic supersolutions provides a
key step in the solution of the optimal control problem in \eqref{optdiv}:

\begin{Prop}\label{thm:repg}
\textup{(i)} The value function $v_*$ is the smallest stochastic
supersolution of the HJB equation \eqref{eq:HJB}
%
\begin{equation}
\label{eq:rep-g} v_*(x) = \min_{g\in\mathcal G^+}g(x),
\end{equation}
for all $x\in\mbb R_+$,
where $\mc G^+$ denotes the family of stochastic supersolutions
of the HJB equation \eqref{eq:HJB}.

\textup{(ii)} For any $a,b\in\mbb R_+$ with $a<b$, representation \eqref
{eq:rep-g} remains
valid for all $x\in(-\infty,b]$ if the set $\mc G^+$ is replaced by the
set $\mc G^+_{a,b}$
of local stochastic supersolutions $g$ on $[a,b]$ satisfying the condition
%
\begin{equation}
\cases{ g(x) = v^*(x), &\quad$ \mbox{for all $x\in[0,a)\cup\{b\}$, and in addition},$
\vspace*{2pt}
\cr
g(a) = v^*(a), & \quad$\mbox{if $X$ has unbounded variation.}$ }
\end{equation}
\end{Prop}

Proposition~\ref{thm:repg},  the proof of  which   is given in
Section~\ref{sec:pf1},
yields the following \emph{(local) verification theorem},
which is one of the main results of the paper:

%
\begin{Thm}\label{cor:repg}\textup{(i)} If there exist $a,b\in\mbb R_+$ with
$b>a\ge0$, $\pi\in\Pi$ and
$g\in\mc G^+$ satisfying $g(x) = v_{\pi,\tau^\pi_a}(x)$
for all $x\in[a, b]$, with $\tau^\pi_a=\inf\{t\ge0\dvtx U^\pi
_t<a\}$,
then it holds $v_*(x)=v_{\pi, \tau^\pi_a}(x)$
for all $x\in[a,b]$.

\textup{(ii)} In particular, if there exist $\pi\in\Pi$ and $g\in\mc G^+$
satisfying $g(x)=v_\pi(x)$ for all $x\in\mbb R_+$, then
$g = v_*$ and $\pi$ is an optimal strategy.
\end{Thm}

\begin{pf}
In view of the dynamic programming equation \eqref{eq:vstar}, it
follows that $v_*$ dominates
$v_{\pi,\tau^\pi_a}$, while the dual representation \eqref
{eq:rep-g} in
Proposition~\ref{thm:repg}
implies $v_*(x)\leq g(x)$ for all $x\in\mbb R_+$, so that when $g$ is
equal to $v_{\pi,\tau^\pi_a}$
on the interval $[a,b]$, it follows that $v^*(x)=g(x)=v_{\pi,\tau^\pi
_a}(x)$ for all $x\in[a,b]$,
which establishes part~(i). Part (ii) follows by a similar line of reasoning.
\end{pf}

This verification result will be used in the
piecewise construction of the value-function $v_*$,
in Sections \ref{aux}--\ref{sec:multib}. It can also be used to deduce
that the value function is affine for large levels of the reserves if
$\nu$ is finite.

\begin{Prop}\label{prop:large}
Let the measure $\nu$ have finite mass.
For some $y\in\R_+$, the function $v_*$ restricted to $[y,\infty)$
takes the form
%
\begin{equation}
\label{large} v_*(x) = x-y + v_*(y)\qquad \mbox{for any $x-y\in\R_+$,}
\end{equation}
and it is optimal to immediately pay out a lump-sum dividend
for all sufficiently large levels of the reserves.
\end{Prop}

\begin{pf}
The local verification theorem [Theorem~\ref{cor:repg}(i)]
in conjunction with Lemma~\ref{lem:smp} imply that condition in \eqref
{large} holds
if the supremum
$m_*:=\sup_{x>y}\int_{(x-y,\infty)}|v_*(x-z)-v_*(y) + z+y-x|\nu(\td z)$
is finite
and
%
\begin{equation}
\label{ellly} \mbox{for all $y\in\R_+$ sufficiently large} \qquad\bigl\{\forall x>y
\dvtx\bigl({}_y{\Gamma}^{v_*}_\infty
\ell_{y,1}\bigr) (x) \leq 0\bigr\}.
\end{equation}
This is verified next. The expression for ${}_y{\Gamma}^{v_*}_\infty
\ell
_{y,1}$ in \eqref{eq:Uy}
for $x > y$ can be bounded above by
\begin{eqnarray*}
&&\eta- q\bigl(x-y+v_*(y)\bigr)+\int_{(x-y,x)}
\bigl|v_*(x-z)-v_*(y)+z+y-x)\bigr|\nu(\td z)
\\
&&\qquad{}+ \int_{(x,\infty)} \bigl|w(x-z)-v_*(y)+z+y-x\bigr|\nu(\td z).
\end{eqnarray*}
Hence, in view of \eqref{eq:vlbub}, the linear bounds in Lemma~\ref
{lem:est}(iii)
and the monotonicity of $w$, the first and second integrals are
bounded above by a constant times $\lambda(1+m)$
and by $\int_{(0,\infty)}|w(-z)|\nu(\td z) + \lambda(|y-v^*(y)|) +
\lambda m$
with $\lambda=\nu(0,\infty)$ and $\lambda m=\int_{(0,\infty)}x\nu
(\td x)$.
Since the integral with $w$ as integrand is finite [as $w\in\mc P$
satisfies \eqref{eq:cw2}]
it follows that $m_*$ is finite. Moreover, as $v_*(y)\to\infty$ and
$v^*(y)-y$ is bounded as $y\to\infty$ [Lemma~\ref{lem:est}(iii)], it is
clear that \eqref{ellly} is satisfied, and the proof is complete.
\end{pf}

\subsection{Proof of the dual representation}\label{sec:pf1}

The proof of Proposition~\ref{thm:repg} is based on a
representation of $v_*$ as the point-wise minimum of a class of
``controlled'' supersolutions of the HJB equation.

\begin{Def}\label{def:H}
For any closed interval $I$, a Borel-measurable function $H\dvtx\mbb
R\to
\mbb R$ is called
a \emph{controlled supersolution for the stochastic control
problem~\textup{\eqref{optdiv}}
on the closed interval $I$} if it holds for any $\pi\in\Pi$ that
%
\begin{equation}
\label{eq:Mtildegpi} \WT M^{H,\pi}_{t}:= \mathrm{e}^{-q(\tau^\pi_I\wedge t)}H
\bigl(U^\pi _{\tau^\pi
_I\wedge t}\bigr) + \int_{[0, \tau^\pi_I\wedge t]}
\mathrm{e}^{-qs} \mu^\pi_{K}(\td s)
\end{equation}
is a UI
$\mbf F$-supermartingale, with $\tau^\pi_I=\inf\{t\ge0\dvtx U^\pi_t
\notin I\}$,
subject to boundary condition
\[
\cases{ H(x) \ge v_*(x),& \quad$\mbox{for $x<y:=\min I$ and $x=z:=\sup I$ if $z<
\infty$, and,}$\vspace*{2pt}
\cr
H(y) \ge v_*(y),&\quad  $\mbox{if $X$ has unbounded
variation}$.}
\]
The family of such functions will be denoted by $\mc H_{I}$.
\end{Def}

\begin{Prop}\label{prop:vc}
For any closed interval $I$
the value-function
$v_*$ restricted to $I$
admits the
following representation:
\begin{eqnarray*}
&& v_*(x) = \min_{H\in\mc H_I} H(x)\qquad \mbox{for all $x\in\mbb R_+$.}
\end{eqnarray*}
\end{Prop}

\begin{pf}
The proof rests on standard arguments. Fix $x\in\mbb R_+$, a closed
interval $I$ in $\mbb R_+$,
and let $H$ be
any element of $\mc H_I$, and $\pi\in\Pi$ any admissible policy. The
supermartingale property and uniform integrability (Definition~\ref
{def:H}) yield
\begin{eqnarray*}
H(x) &\ge& \lim_{t\to\infty}\E_x \biggl[
\mathrm{e}^{-q(\tau
_I^\pi\wedge
t)}H\bigl(U^\pi_{\tau_I^\pi\wedge t}\bigr) + \int
_{[0, \tau_I^\pi\wedge
t]}\mathrm{e}^{-qs}\mu_K^\pi(
\td s) \biggr]
\\
&\ge& \E_x \biggl[\mathrm{e}^{-q\tau_I^\pi}v_*\bigl(U^\pi_{\tau_I^\pi
}
\bigr) + \int_{[0, \tau_I^\pi]}\mathrm{e}^{-qs}
\mu_K^\pi(\td s) \biggr],
\end{eqnarray*}
where the convention $\exp\{-\infty\}=0$ is used.
Taking the supremum over $\pi\in\Pi$ and using the dynamic
programming equation
(Proposition~\ref{prop:mart}) show that $H(x)\ge v_*(x)$. Since $H\in
\mathcal H_I$ was
arbitrary, it holds thus
\[
\inf_{H\in\mathcal H_I}H(x)\ge v_*(x).
\]
The inequality in the display is in fact an equality since $v_*$ is a
member of $\mathcal H_I$, by virtue of Doob's optional stopping theorem and
the fact that $V^\pi$ is a UI supermartingale [Lemma~\ref{lem:est}(iii)].
\end{pf}

The proof of the representations of the value function $v_*$
in Proposition~\ref{thm:repg} rests on the fact that for any
admissible policy $\pi\in\Pi$ and stochastic supersolution there exists
a corresponding
``controlled'' supermartingale.

\begin{Lemma}[(Shifting lemma)]\label{lem:shift}
Let $I\subset\mbb R_+$ be any closed interval. If $g$ is a local
stochastic supersolution on $I$,
then $g$ is a controlled supersolution on~$I$.
\end{Lemma}

Given the shifting lemma, the proof of the dual representations in
Proposition~\ref{thm:repg}
can be completed as follows:

\begin{pf*}{Proof of Proposition~\ref{thm:repg}}
(i) The representation follows from Proposition~\ref{prop:vc} in view
of the following two
observations: (a) $\mathcal G^+$ is contained in $\mathcal H_{[0,\infty
)}$ [Remark~\ref{rem:stst}(i)] and (b)
$v_*$ is an element of the set $\mc G^+$ [by Lemma~\ref{lem:est}(iii)].

(ii) The proof is analogous to that of part (i),
using the facts $\mc G^+_{a,b}\subset\mathcal H_{[a,b]}$ [Lemma~\ref
{lem:shift}(ii)]
and $v_*\in\mc G^+_{a,b}$ [by Remark~\ref{rem:stst}(i) and Doob's
Optional Stopping theorem].
\end{pf*}

\begin{pf*}{Proof of Lemma~\ref{lem:shift}}
Fix arbitrary $\pi\in\Pi$ and
$s,t\in\mbb R_+$ with $s< t$.
Note that $\WT M^{g,\pi}$ is $\mbf F$-adapted (as $g$ is a Borel-measurable),
while $\WT M^{g,\pi}$ is UI by the linear growth condition and
Lemma~\ref
{lem:est}.
Furthermore, the following (in)equalities hold true:
\[
\E \bigl[\WT M_t^{g,\pi} |\mc F_{s\wedge\tau^\pi} \bigr]
\stackrel {\mathrm{(a)}} {=} \lim_{n\to\infty} \E \bigl[\WT
M^{g,\pi_n}_t |\mc F_{s\wedge\tau^{\pi}} \bigr] \stackrel {\mathrm{(b)}} {
\leq} \lim_{n\to\infty} \WT M^{g,\pi_n}_{s\wedge\tau^{\pi}}
\stackrel{\mathrm{(c)}} {=} \WT M_{s\wedge\tau^{\pi}}^{g,\pi} \stackrel{\mathrm{(d)}} {=} \WT
M^{g,\pi}_s,
\]
where the sequence $(\pi_n)_{n\in\mbb N}$ of strategies is defined
by $\pi_n=\{D^{\pi_n}_t, t\in\mbb R_+\}$ with $D_0^{\pi_n} =
D_0^\pi$ and
\[
D^{\pi_n}_u = \cases{ \sup\bigl\{D^\pi_{v}
\dvtx v< u, v\in\mbb T_n\bigr\}, &\quad $0<u < \tau^\pi$,
\vspace *{2pt}
\cr
D^{\pi_n}_{\tau^\pi-}, &\quad $u\ge\tau^\pi$,
}
\]
with $\mbb T_n:= ( \{t_k:=s+(t-s)\frac{k}{2^n}, k\in\mbb
Z  \}\cup\{0\}  )\cap\mbb R_+$.
Since $s$ and $t$ are arbitrary,
it thus follows that $\WT M^{g,\pi}$ is a $\mbf F$-supermartingale.

The remainder of the proof is devoted to the verification of the
(in)equalities (a)--(d) in above display.
(a) Note that the sequence $(D^{\pi_n})_n$ is monotone
($D^{\pi_n}\leq D^{\pi_{n+1}}$ for $n\in\mbb N$) and tends to $D^\pi$
as $n$
tends to infinity, and $D^{\pi_n}$ is equal to $D^{\pi_n}_{\tau^\pi-}$
on the interval
$[\tau^\pi,\infty)$, for each $n\in\mbb N$. Thus the monotone convergence
theorem (MCT) in combination with an integration-by-parts implies
$\int_{[0, \tau^{\pi}\wedge t]}\mathrm{e}^{-qs}\,\td D^{\pi_n}_s
\nearrow
\int_{[0, \tau^\pi\wedge t]}\mathrm{e}^{-qs}\,\td D^\pi_s$.
Also, in the case $K>0$, it holds
$\int_{[0, \tau^{\pi}\wedge t]}\mathrm{e}^{-qs}\,\td N^{\pi_n}_s
\nearrow
\int_{[0, \tau^\pi\wedge t]}\mathrm{e}^{-qs}\,\td N^\pi_s$. Hence, by
right-continuity of the function $g$,
it holds
%
\begin{equation}
\label{eq:point} \WT M^{g,\pi_n}_{t\wedge\tau^\pi} \longrightarrow\WT
M^{g,\pi
}_{t\wedge\tau^\pi} \qquad\mbox{as $n\to\infty$}, \mbox{$\P$-a.s.}
\end{equation}
As the collection $(\WT M^{g,\pi_n}_{t\wedge\tau^\pi})_n$ is UI,
Lebesgue's dominated convergence theorem implies that the equality (a)
holds true.
Equality (c) is a consequence of the pointwise convergence in \eqref{eq:point}
(which also holds with $t$ replaced by $s$), while
(d) follows since it holds $\WT M^{g,\pi}_s = \WT M^{g,\pi}_{s\wedge
\tau
^\pi}$
(by definition of the process $\WT M^{g,\pi}$).

Finally, inequality (b) is verified, in what constitutes the key step
of the proof.
Denote $T_i:=\tau^{\pi}\wedge t_i$
and $M = \WT M^{g,\pi_n}$,
$D = D^{\pi_n}$, and observe that the folowing decomposition holds true:
%
\begin{eqnarray}
 M_t - M_s = \sum_{i=1}^{2^n}
Y_i + \sum_{i=1}^{2^n}{Z_i}\nonumber\\
\eqntext{\mbox{with }
 Y_i =\mathrm{e}^{-q T_i}g (X_{T_i} -
D_{T_{i-1}} ) - \mathrm{e}^{-q
T_{i-1}}g (X_{T_{i-1}} -
D_{T_{i-1}} ),}
\end{eqnarray}
with $Z_i=\mathrm{e}^{-q T_i}(g(X_{T_i}-D_{T_i}) - g(X_{T_i} -
D_{T_{i-1}}) + \Delta D_{T_i} - K)\mbf1_{\{\Delta D_{i}>0\}}$
and $\Delta D_{i}= D_{T_i}
- D_{T_{i-1}}$.
The strong Markov property of $X$ and
the definition of $U$ imply that $\E[Y_i|\mc F_{T_{i-1}}]$
is equal to
%
\begin{eqnarray}
\label{eq:YEDOOB}&& \mathrm{e}^{-q T_{i-1}}\E \bigl[\mathrm{e}^{-q(T_i -
T_{i-1})}g
(U_{T_{i-1}} + X_{T_i} - X_{T_{i-1}} ) - g
(U_{T_{i-1}} ) |\mc F_{T_{i-1}} \bigr]
\nonumber
\\[-8pt]
\\[-8pt]
\nonumber
&&\qquad= \mathrm{e}^{-q T_{i-1}} \E_{U_{T_{i-1}}} \bigl[\mathrm{e}^{-q\tau_i}g
(X_{\tau
_i} ) - g(X_0) \bigr],
\end{eqnarray}
with $\tau_i = T_i\circ\theta_{T_{i-1}}$, where $\theta$ denotes the
translation-operator.
The right-hand side of \eqref{eq:YEDOOB}
is nonpositive as a consequence of
the supermartingale property \eqref{eq:ggmartI} (with $I=\mbb R_+$) and
Doob's optional stopping theorem.
Furthermore, in view of the bound $\mathtt d_g(x)\ge1$ for any $x\in
\mbb R_+\setminus\{0\}$
it follows that all the $Z_i$ are
nonpositive in the case $X_{T_i}-D_{T_i}\ge0$, while,
in the case $X_{T_i}-D_{T_i}<0$, it holds that $Z_i$ is zero,
since $T_i=\tau^\pi$, so that, by construction,
$\Delta D_i = D^{\pi_n}_{\tau^\pi} - D^{\pi_n}({\tau^+_n})=0$ with
$\tau
^+_n=\sup\{v<\tau^\pi: v\in\mbb T_n\}$.
Hence, the tower-property of conditional
expectation yields
\[
\E[M_t- M_s|\mc F_s] \leq\sum
_{i=1}^{2^n} \mbf 1_{\{T_{i-1}>s\}}\E \bigl[
\E[Y_i|\mc F_{T_{i-1}}]|\mc F_s \bigr] \leq0.
\]
This establishes inequality (b), and the proof is complete.
\end{pf*}

\section{Gerber--Shiu functions}\label{sec:RLevy}
A key-ingredient for the solution of the optimal control problem \eqref
{optdiv} is a family of martingales
given in terms of \emph{Gerber--Shiu functions}, a nonstandard
terminology; see Definitions \ref{d:GS} and \ref{Def:Fw}. While the
(homogeneous) $q$-scale function $W^{(q)}$
is defined to be equal to $0$ on the set $(-\infty,0)$, Gerber--Shiu
functions are ``inhomogeneous $q$-scale functions''
corresponding to nonzero boundary conditions $w$ on the negative half-line.

The definition of Gerber--Shiu functions is phrased in terms of $w$ and
$W^{(q)}$
of which next a number of well-known properties are recalled that will
be deployed in the sequel; refer to the review article Kyprianou et
al. \cite{KypScale}, Chapters 2, 3, for proofs and references.
The function $W^{(q)}$ [see \eqref{Wq} for its definition]
is a ``$q$-harmonic function'' for the process $X$ stopped at first
entrance into $(-\infty,0)$. Specifically, for any $a\in\mbb R$, the
stopped process
%
\begin{eqnarray}
\label{eq:Wqmart} && \bigl(\mathrm{e}^{-q(t\wedge T_a^-)} W^{(q)}
(X_{t\wedge
T_a^- }- a ), t\in \mbb R_+ \bigr)
\nonumber
\\[-8pt]
\\[-8pt]
\nonumber
&&\qquad \mbox{is an $\mbf F$-martingale,
with }T_a^-:= T_{[a,\infty)}=\inf\{t\in\mbb R_+\dvtx
X_t<a\}.
\end{eqnarray}
Furthermore,
the function $W^{(q)}$ is well-known to be continuous and
nondecreasing on $[0,\infty)$, and right-
and left-differentiable on $(0,\infty)$, with
the right-derivative and left-derivative at $x>0$
denoted by $W^{(q)\prime}(x)$ and $W^{(q)\prime}_-(x)$, respectively,
which are right- and left-continuous and satisfy
%
\begin{equation}
\label{eq:WqD} W^{(q)\prime}(x) \leq W^{(q)\prime}_-(x),\qquad x>0,
\end{equation}
by continuity and log-concavity of $W^{(q)}|_{\mbb R_+}$. In
particular, if $\nu_{0,1}$ [which was defined in \eqref{exoruinliq}]
is infinite, the function $W^{(q)}|_{(0,\infty)}$ is $C^1$, while
$W^{(q)}|_{(0,\infty)}$ is $C^2$ with $W^{(q)\prime}(0+)=\frac
{2}{\sigma
^2}$ if the Gaussian coefficient $\sigma^2$ is strictly positive.\vspace*{1pt}

A function will be referred to as a Gerber--Shiu function if it
satisfies the following conditions:

%
\begin{Def}\label{d:GS} Given $a\in\mbb R$
and a \emph{pay-off} $w\dvtx(-\infty,a]\to\mbb R$ with $w\in
\mathcal R_a$,
the function $F\dvtx\mathbb R\to\mathbb R$ is called a \emph
{Gerber--Shiu function
for payoff $w$} if $F(x-a)=w(x)$ for $x<a$, and
%
\begin{equation}
\label{eq:Fwm} \bigl(\mathrm{e}^{-q(t\wedge T_a^-)} F (X_{t\wedge
T_a^-}-a ), t\in\mbb
R_+ \bigr)\qquad \mbox{is an $\mbf F$-martingale}.
\end{equation}
\end{Def}

Of course, such a function $F$ is not unique (as the addition of
multiples of $W^{(q)}$ to a Gerber--Shiu
function yields another Gerber--Shiu function). It is shown below that
there exists a special choice $F_w$ of Gerber--Shiu function
that is \emph{continuous} on $\mbb R$
for continuous payoffs $w$ and \emph{continuously differentiable}
on $\mbb R$ if $X$ has unbounded variation and
$w$ is continuously differentiable (recall that $W^{(q)}$ is continuous
nor continuously differentiable
on $\mbb R$ in general). The function $F_w$ is defined as follows:

\begin{Def}\label{Def:Fw} Let $q\ge0$ and $w\in\mc R_0$.
The function $F_{w}\dvtx\mbb R\to\mbb R$ is given by $F_w(x)=w(x)$ for
$x< 0$, and by
%
\begin{eqnarray}\qquad
\label{eq:Fwr}  F_w(x) &=& w(0) + w'_-(0) x - \int
_0^x W^{(q)}(x-y) J_w(y)
\,\td y, \qquad x\in\mbb R_+,\mbox{ with}
\\
J_w(x) &=& \bigl({}_0\Gamma_\infty^w
\ell_{0,w'_-(0)} \bigr) (x), \label{J}
\end{eqnarray}
where $_0 \Gamma_\infty^w\ell_{0,w'_-(0)}$ is given in \eqref{eq:Uy}
[with $a=0$ and $b=w_-'(0)$].
\end{Def}

The following result confirms that the function $F_{_a w}$ is a
Gerber--Shiu function that ``inherits'' the
continuity/differentiability from the function $w$,
where, for any function $f$ and $a\in\mbb R$, $_a f$ denotes the
composition of $f$
with the translation-operator $\theta_a$,
%
\begin{equation}
\label{eq:trans} _a f:= f \circ\theta_a:= f(\cdot+a).
\end{equation}

\begin{Thm}\label{thm:sd}
Let $a\in\mbb R$ and $w\in\mc R_a$. Then $_a w\in\mc R_0$
and the function $F_{_a w}$ is a Gerber--Shiu function for payoff $w$
satisfying
%
\begin{eqnarray}
\label{dFw0} \cases{ F_{_a w}(0) = w(a),\vspace*{2pt}
\cr
F_{_a w}'(0+) = w'_-(a), & \quad $\mbox{in the
case $\sigma^2>0$ or $\nu _{0,1}=\infty$}$.}
\end{eqnarray}
Furthermore, $F_{_a w}|_{\mbb R_{+}}$ is right-differentiable,
with right-derivative at $x\in\mbb R_+$ denoted by $F'(x)$.
If $_a w$ is continuous, then $F_{_a w}$ is continuous,
and, in the case $w\in C^1(\mbb R_-)$ and $\{\sigma^2>0\mbox{ or } \nu
_{0,1}=\infty\}$,
it holds $F_{_a w}\in C^1(\mbb R)$.
\end{Thm}

An example of a Gerber--Shiu function is the
\emph{Gerber--Shiu penalty function} $\mc V_w$ corresponding to
penalty $w$
\[
\mc V_w(x) = \E_x \bigl[\mathrm{e}^{-q T_0^-}w(X_{T_0^-})
\bigr],
\]
which admits the following explicit expression
in terms of the functions $W^{(q)}$ and~$F_w$ (see Biffis and
Kyprianou \cite{Biffis}
for an equivalent representation of $\mc V_w$
in terms of $W^{(q)}$):

\begin{Prop}[(Gerber--Shiu penalty function)]
\label{prop:GS}
Let $w\in\mc R_0$. For any \mbox{$x\in\mbb R$} it holds
%
\begin{eqnarray}
\label{eq:BK} \mc V_w(x)& =& F_w(x) -
W^{(q)}(x)\kappa_w\qquad \mbox{with}
\\
\kappa_w &:=& \biggl[\frac{\sigma^2}{2}w'(0-) +
\frac{q}{\Phi
(q)}w(0) - \mc Lw_\nu\bigl(\Phi(q)\bigr) \biggr],\label{eq:kappa}
\end{eqnarray}
where $\mc Lw_\nu$ denotes the Laplace transform of the function
$w_\nu
(x)=\int_{(x,\infty)}[w(x-z)-w(0)]\nu(\td z)$, $x>0$.
\end{Prop}

For later reference two further exit identities are recorded that are
also expressed
in terms of $W^{(q)}$ and $F_w$. First, the {two-sided exit identity}
of $X$ on the interval $[a,b]$
which involves the distribution of the pair $(T_{a,b}, X_{T_{a,b}})$
where $T_{a,b}:= T_{[a,b]} = T_a^-\wedge
T_b^+$, with $T^+_b:= T_{(-\infty,b]} = \inf\{t\in\mbb R_+\dvtx X_t
>b\}$,
denotes the first exit time from the interval $[a,b]$.
Second, a \emph{absorption/reflection exit identity}
on the interval $[a,b]$ which concerns the law of the pair $(\tau_a(b),\break
Y^b_{\tau_a(b)})$
and the expected local time up to $\tau_a(b)$ at the level $b$ of $Y^b$
where $\tau_a(b) =\inf\{t\in\mbb R_+\dvtx Y^b_t < a\}$ denotes the
first-passage time into
the set $(-\infty,a)$ of the process $Y^b=\{Y^b_t, t\in\mbb R_+\}$
given by
%
\begin{equation}
\label{eq:tauab} Y^b_t = X_t - \ovl
X_t^b\qquad \mbox{with } \ovl X^b_t=
\sup_{s\leq
t} (X_t - b )\vee0.
\end{equation}
The identities are given as follows:

\begin{Prop}\label{prop:two}
Given $a\in\mbb R$
and a \emph{pay-off} $w\dvtx(-\infty,a]\to\mbb R$ with $w\in
\mathcal R_a$,
the following hold for all $b,\d,\beta\in\mbb R$ with $a<b<\infty$ and
$x\in(a,b)$:
%
\begin{eqnarray}\label{d:wxt0a}
&&\mathbb E_x \bigl[\mathrm{e}^{-q
T_{a,b}}w
(X_{T^-_{a}} )\mbf1_{\{T_a^-<T^+_b\}} \bigr] + \d \mathbb E_x
\bigl[\mathrm{e}^{-q T^+_{b}}\mbf 1_{\{T_a^->T^+_b\}} \bigr]
\nonumber
\\[-8pt]
\\[-8pt]
\nonumber
&&\qquad = F_{_aw}(x-a) + W^{(q)}(x-a)
\frac{\d- F_{_aw}(b-a)}{W^{(q)}(b-a)},
\\
\label{d:wYta}&&
\mathbb E_x \bigl[\mathrm{e}^{-q
\tau_{a}(b)}w
\bigl(Y^b_{\tau_{a}(b)} \bigr) \bigr] + \beta\E_x
\biggl[\int_{[0,\tau_{a}(b)]}\mathrm{e}^{-qs}\,\td\ovl
X^b_s \biggr]
\nonumber
\\[-8pt]
\\[-8pt]
\nonumber
&&\qquad= F_{_aw}(x-a) + W^{(q)}(x-a)\frac{\beta-F_{_aw}'
(b-a)}{W^{(q)\prime}(b-a)}.
\end{eqnarray}
\end{Prop}

The proofs of Theorem~\ref{thm:sd} and Proposition~\ref{prop:GS} rests
on the following auxiliary results
(shown in Section~\ref{sec:GSp}):

\begin{Lemma}\label{F1}Let $w\in\mc R_0$.
The function $F_w|_{\mbb R_+}$ real-valued and continuous and admits
the following alternative representation:
for $x\ge0$,
%
\begin{eqnarray}\qquad
\label{eq:repFF} F_w(x) = \frac{\sigma^2w_-'(0)}{2} W^{(q)}(x) +
w(0) Z^{(q)}(x) - \int_0^x
W^{(q)}(x-y)w_\nu(y)\,\td y
\nonumber
\\[-8pt]
\\[-8pt]
\eqntext{\mbox{with }Z^{(q)}(x)=1 + \int_0^x
W^{(q)}(y)\,\td y.}
\end{eqnarray}
In particular, it holds $F_w(0) = w(0)$ and $\int_0^x|w_\nu(y)|\,\td
y<\infty$ for any $x\ge0$, and in the case that
$X$ has bounded variation $w_\nu(0+)<\infty$.
\end{Lemma}

%
\begin{Lemma}\label{lem:} Let $w\in\mc R_0$.
\textup{(i)} $F_w(x)/W^{(q)}(x) \to\kappa_w$ as $x\to\infty$.\vspace*{-6pt}
\begin{longlist}[(iii)]
\item[(ii)] $F_w(x)$ is left- and right-differentiable at any $x>0$ with
right-derivative at $x>0$ given by
%
\begin{eqnarray}
\label{eq:Fwd0} F_{w}'(x) &=& w_-'(0) -
\int_{[0,x)} J_w(x-y) W^{(q)}(\td y)
\nonumber
\\[-8pt]
\\[-8pt]
\nonumber
&=& F_{w,-}'(x) - W^{(q)}(0)
\bigl(J_w(x+) - J_w(x-)\bigr),
\end{eqnarray}
where $F^{\prime}_{w,-}(x)$ denotes the left-derivative of $F_w$ at $x$.
In particular, $F_{w}'(0) = w_-'(0)$ if $X$ has unbounded variation,
and $F_{w}'(0) = w_-'(0) - W^{(q)}(0)J_w(0+)$ if $X$ has bounded variation.

\item[(iii)] The function $x\mapsto F_{w}'(x)$ is right-continuous on $\mbb
R_+\setminus\{0\}$,
and is $C^1$ on $\mbb R_+\setminus\{0\}$ in the case $w\in C^1(\mbb R_-)$.
\end{longlist}
\end{Lemma}

Given these two results the proofs of Proposition~\ref{prop:GS} and
Theorem~\ref{thm:sd} can
be completed as follows:

\begin{pf*}{Proof of Proposition~\ref{prop:GS}}
Writing $\mc V_w(x)=w(0)\mc V_{\mbf1}(x) +
\E_x[\mathrm{e}^{-q T_0^-}\*(w(X_{T_0^-}) - w(0))]$, where $\mbf1$
denotes the function
on $\mbb R_-$ that is constant equal to one,
and applying the compensation formula (e.g., Bertoin \cite{bert96}, Chapter~O) to the Poisson
point process $(\Delta X_t, t\in\mbb R_+)$ yields the following expressions
for any $x\in\mbb R_+$:
%
\begin{eqnarray}\label{int}
\nonumber
\mc V_w(x) - w(0) \mc V_{\mbf1}(x) &=& \int
_{[0,\infty)}\int_{(y,\infty)} \bigl(w(y-z) - w(0)\bigr)
\nu(\td z)U^q(x,\td y)
\\
&=& W^{(q)}(x) \mc Lw_\nu\bigl(\Phi(q)\bigr) - \int
_{0}^x W^{(q)}(x-y) w_\nu
(y)\,\td y,
\\
U^{q}(x,\td y) &=& \bigl[W^{(q)}(x) \mathrm{e}^{-\Phi(q) y}
- W^{(q)}(x-y)\bigr]\,\td y,\qquad y > 0,
\nonumber
\end{eqnarray}
where $U^q(x,\td y)$ denotes
the $q$-potential measure of $X$ under
$\P_x$ killed upon entering $(-\infty,0)$.
It follows from Lemmas \ref{F1} and \ref{lem:} that
the integrals in~\eqref{int} are finite.
Deploying the form of the Laplace transform of $T_0^-$,
$\mc V_{\mbf1}(x) = Z^{(q)}(x) - q\Phi(q)^{-1}W^{(q)}(x)$,
and the definition of $F_w$ leads to \eqref{eq:BK} [since
the term $\frac{\sigma^2}{2}w'(0-)W^{(q)}(x)$ cancels].
\end{pf*}

\begin{pf*}{Proof of Proposition~\ref{prop:two}}
Denote the left-hand side of \eqref{d:wYta} by $\mc U_{w,\beta}^{a,b}(x)$,
and let $e_{0,a}$ be the function with domain $(-\infty,a]$ that is
constant equal to 1.
Another application of the compensation formula
yields the following representation of $\mc U_w^{a,b}(x)$
for $x\in[a,b]$:
\begin{eqnarray*}
&&\mc U_{w,\beta}^{a,b}(x) - w(0) \mc
U^{a,b}_{e_{0,a},0}(x) - \beta\mc U^{a,b}_{0,1}(x)
\\
&&\qquad= \int_{[a, b]}\int_{(y,\infty)}
\bigl(w(y-z) - w(0)\bigr) \nu(\td z)R^q_{a,b}(x,\td y)\qquad
\mbox{with}
\\
&&R^q_{a,b}(x,\td y) = \frac{W^{(q)}(x-a)}{W^{(q)\prime}(b-a)}
W^{(q)}(b - \td y) - W^{(q)}(x-y)\,\td y,
\\
&&\mc U^{a,b}_{e_{0,a},0}(x) = \E_x\bigl[
\mathrm{e}^{-q\tau_a(b)}\bigr] = Z^{(q)}(x-a) - q\frac{W^{(q)}(x-a)}{W^{(q)\prime}(b-a)}
W^{(q)}(b-a),
\\
&&\mc U^{a,b}_{0,1}(x) = \E_x \biggl[\int
_{[0,\tau_{a}(b)]}\mathrm {e}^{-qs}\,\td \ovl
X^b_s \biggr] = \frac{W^{(q)}(x-a)}{W^{(q)\prime}(b-a)},
\end{eqnarray*}
where $R^q_{a,b}(x,\td y)$, $y\in[a,b]$, is
the $q$-resolvent measure
of $Y^b$ killed upon entering $(-\infty,a)$
(from Pistorius \cite{P}, Theorem~1)
and the final two identities in the previous display
are from Avram et al. (\cite{AKP}, Theorem~1, \cite{APP}, Proposition~1).
Combining these expressions with representation \eqref{eq:repFF} of
$F_w$ and taking note of the fact that
the term $\frac{\sigma^2}{2}{}_a w'(0-)W^{(q)}(x)$
again cancels yields that \eqref{d:wYta} holds true.
Equation \eqref{d:wxt0a} follows by a similar line of reasoning.
\end{pf*}

\begin{pf*}{Proof of Theorem~\ref{thm:sd}}
That $F_{_a w}$ is a Gerber--Shiu function follows from~\eqref{eq:BK} (with $F_w$ replaced by $F_{_a w}$),
the strong\vspace*{1pt} Markov property of $X$
and the martingale property \eqref{eq:Wqmart} of $W^{(q)}$.
The martingale property \eqref{eq:Fwm} was shown in Proposition~\ref{prop:GS}.
The asserted continuity follows from the relation \eqref{dFw0}
combined with the continuity of $_a w$ and $F_{_a w}|_{\mbb R_+}$
(Theorem~\ref{thm:sd}). The assertion that $F_{_a w}$ is $C^1(\mbb R)$
is a consequence
of the following two observations: (i) $F_{_a w}|_{\mbb R_+\setminus\{
0\}}$ is $C^1(\mbb R_+\setminus\{0\})$
[by Lemma~\ref{lem:}(ii)]; (ii) $_a w$ is $C^1(\mbb R_-)$ (by assumption)
and $w_-'(a) = {}_a w'_-(0) = F_{_a w}'(0)$
[by Lemma~\ref{lem:}(ii)].
\end{pf*}

\subsection{Proofs of Lemmas \texorpdfstring{\protect\ref{F1}}{5.6} and \texorpdfstring{\protect\ref{lem:}}{5.7}}\label{sec:GSp}
\mbox{}

\begin{pf*}{Proof of Lemma~\ref{F1}}
First it is verified that the function on the right-hand side of \eqref
{eq:repFF} is continuous on $\mbb R_+$.
This follows from the continuity on $\mbb R_+$ of $W^{(q)}(x)$, $Z^{(q)}(x)$
and of the final term in \eqref{eq:Fwr}, as functions of $x$. The
continuity of the integral
is a consequence of Lebesgue's dominated convergence theorem
and the finiteness of $\int_0^x|w_\nu(y)|\,\td y$ for any $x\ge0$,
which in turn holds as $w$ is c\`adl\`ag and left-differentiable at 0
($w\in\mc R_0$) and
$\nu$ satisfies the integrability condition $\int_0^1z^2\nu(\td
z)<\infty$.
Furthermore, in the case that $X$ has paths of bounded variation,
it holds that $\int_0^1z\nu(\td z)$ is finite, and a similar line of reasoning
yields that $w_\nu(0+)$ is finite.

As it follows by a similar argument that also $F_w$ is continuous on
$\mbb R_+$
it suffices to show that the Laplace transforms of the right-hand side of \eqref
{eq:repFF} and of \eqref{eq:Fwr}
coincide in order to verify
the representation \eqref{eq:repFF}. Note that
the Laplace transform $\mc L|\tilde w_\nu|(\theta)$ of $|\tilde w_\nu|$
is finite for any $\theta>0$ in view of the integrability condition
\eqref{eq:cw2} and since $\int_0^1|w_\nu(y)|\,\td y$ is finite.
Taking the Laplace transform of \eqref{eq:Fwr}, using the forms \eqref
{eq:psi} and
\eqref{Wq} of the Laplace exponent $\psi(\theta)$ and the
Laplace transform $\mc LW^{(q)}$ and rearranging terms yields
%
\begin{eqnarray}\label{2}
\nonumber
\mc LF_w(\theta) &=& \mc LW^{(q)}(\theta) \biggl[
\frac{\sigma^2}{2}w_-'(0) + \frac{\psi(\theta)}{\theta} w(0) - \mc
Lw_\nu(\theta) \biggr],\qquad \theta>\Phi(q),
\\
&=& \theta^{-1}\cdot w(0) + \theta^{-2}\cdot
w'_-(0) - \bigl(\psi(\theta )-q\bigr)^{-1} \mc
LJ_w(\theta),
\nonumber
\\[-8pt]
\\[-8pt]
\nonumber
\qquad\mc LJ_w(\theta) &=& \theta^{-1} \cdot\bigl[{
\psi'(0)w_-'(0) - qw(0)}\bigr] + \mc L\tilde
w_\nu(\theta) - \theta^{-2}\bigl[qw'_-(0)
\bigr],
\\
\mc L\tilde w_\nu(\theta) &=& \mc Lw_\nu(\theta) +
w_-'(0)\cdot \theta^{-2} \int_{(0,\infty)}
\bigl[\mathrm{e}^{-\theta x} - 1 + \theta x\bigr]\nu(\td x).
\nonumber
\end{eqnarray}
Termwise inverting \eqref{2} yields the expression \eqref{eq:repFF}.

By letting $x\to0$ in \eqref{eq:repFF}, in combination with the facts
$\sigma^2 W^{(q)}(0+)=0$ and $Z^{(q)}(0+)=1$ and the fact that the
integral tends to zero (again by
Lebesgue's dominated convergence theorem), it follows that $F_w(0)=w(0)$.
\end{pf*}

\begin{pf*}{Proof of Lemma~\ref{lem:}}
(i) The limit \eqref{eq:kappa} follows from \eqref{eq:Fwr} or \eqref
{eq:repFF} using $W^{(q)}(x)\sim\mathrm{e}^{\Phi(q)x}/\psi'(\Phi
(q))$ as
$x\to\infty$.

(ii) Observe first that $J_w$ is c\`adl\`ag on $\mbb R_+\setminus\{0\}$,
by noting that $w_\nu(x)$ is c\`adl\`ag at any $x> 0$
[as a consequence of the facts that $w$ is c\`adl\`ag,
left-differentiable at zero,
and satisfies the integrability condition \eqref{eq:cw2}].

The continuity of $W^{(q)}$ on $\mbb R_+$, \eqref{eq:cw2} and the
finiteness of $\int_0^1|w_\nu(y)|\,\td y$
(Lemma~\ref{F1}) imply that the integral $\int_0^x
|W^{(q)}(x-y)J_w(y)|\,\td y$ is finite for any $x>0$.
A change of the order of integration in \eqref{eq:repFF}, justified by
Fubini's theorem, implies for $x> 0$
the integral $\int_0^x J_w(x-y)W^{(q)}(y)\,\td y $ is equal to
\[
W^{(q)}(0)\int_0^x
J_w(u)\,\td u + \int_0^x \int
_0^{x-z} J_w(u)\,\td u
W^{(q)\prime}(z)\,\td z.
\]
As a consequence, it follows that the right- and left-derivatives
$F_w'(x)$ and $F'_{w,-}(x)$ are equal to
$w'_-(0) - \int_0^x J_w((x-z)\pm)W^{(q)\prime}(z)\,\td z -
W^{(q)}(0)J_w(x\pm)$,
respectively, at any $x>0$.
Thus the difference $F'_{w}(x) - F'_{w,-}(x)$ is as stated in \eqref{eq:Fwd0}.
An application of Lebesgue's dominated convergence theorem implies that
the integral in the previous line converges to zero when $x$ tends to 0.
The right-continuity of $J_w$ and the fact that $W^{(q)}(0)$ is 0
precisely if $X$ has unbounded variation,
yields the stated form of $F_w'(0)$.

(iii) The right-continuity follows from the right-continuity of $J_w$
on $\mbb R_+\setminus\{0\}$
and Lebesgue's dominated convergence theorem. In the case $w\in
C^1(\mbb R_-)$, a similar argument as at the start of part (ii) implies
that $J_w$ is continuous on $\mbb R_+$.
It follows thus from \eqref{eq:Fwd0} that $F_w'(x)$ is continuous at
any $x>0$.
\end{pf*}

\subsection{Exponential and polynomial boundary conditions}
For later reference it is noted that in the case that the payoff $w$ is
exponential, $w(x)=\mathrm{e}^{xv}$ for some $v\in\mbb R$, or is a monomial,
$w(x)=x^k$, the solutions of the two-sided and mixed
absorbing/reflected exit problems simplify and
can be expressed in terms of the functions $Z^{(q,v)}$ and $Z_k$ that
are specified as follows:

\begin{Def}\label{def:Zqv}(i) For $q,v\in\mbb R_+$, the function
$Z^{(q,v)}\dvtx\R\to\R$ is defined by $Z^{(q,v)}(x)=\mathrm
{e}^{vx}$ for
$x< 0$, and by
%
\begin{equation}
\label{eq:Zv} \qquad Z^{(q,v)}(x) = \mathrm{e}^{vx} + \bigl(q-
\psi(v)\bigr)\int_0^x \mathrm{e}^{v(x-y)}W^{(q)}(y)
\,\td y, \qquad x\in\mbb R_+.
\end{equation}
(ii) With $n_0$ the largest integer such that
$\int_{(-\infty, -1)}|x|^n\nu(\td x)<\infty$, the related family of
functions $Z_k\dvtx\R\to\R$, $k=0,\ldots, n$, is defined by
%
\begin{equation}
\label{eq:Zpower} Z_k(x) = \frac{\partial^k}{\partial v^k}\bigg\vert
_{v=0+} Z^{(q,v)}(x).
\end{equation}
\end{Def}

As suggested above, $Z^{(q,v)}$ and $Z_k$ are in fact Gerber--Shiu
functions of the exponential and monomial pay-offs $e_v, p_k\dvtx\mbb
R_-\to
\mbb R$, which for any $v\in\mbb R$ and $k=1,\ldots, n_0$ are
given by $e_v(x):=\mathrm{e}^{vx}$ and $p_k(x):=x^k$.

%
\begin{Cor}\label{prop:polpen}
For any $q>0$, $v\in\mbb R$ and $k=1,\ldots, n_0$, $Z^{(q,v)}$ and
$Z_k$ are Gerber--Shiu functions with payoffs
$e_{v,a}:= {}_{a} e_v$ and $p_{k,a} ={} _{a} p_k$, the translations of
$e_v$ and $p_k$, respectively.
\end{Cor}

\begin{pf} The assertion concerning $Z^{(q,v)}$ directly follows from
Theo-\break rem~\ref{thm:sd}
since the function $Z^{(q,v)}$ is equal to
the Gerber--Shiu function $F_w$ corresponding to $w=e_v$.
The two functions coincide since both are continuous on
$\mbb R_+$ and it holds
%
\begin{equation}
\label{eq:LTZ} \mc L F_{e_v}(\theta) = \mc LZ^{(q,v)}(\theta) =
\bigl(\psi(\theta) - q\bigr)^{-1}\frac{\psi(\theta) - \psi(v)}{\theta- v}.
\end{equation}
The proof of the assertion concerning $Z_k$ is similar and omitted.
\end{pf}


%
\begin{Rem}\label{rmZ}
(i)
For $v\ge0$, the function $x\mapsto Z^{(q,v)}(x)$ is strictly
increasing on $\mbb R_+$. In particular, for $x>0$ and
$v>\Phi(q)$, $Z^{(q,v)\prime}(x)$ is equal to
%
\begin{equation}
\label{eq:Zd} Z^{(q,v)\prime}(x) = \bigl(\psi(v)-q\bigr)\int
_{x}^{\infty}\mathrm {e}^{v(x-y)}W^{(q)\prime
}(y)
\,\td y,
\end{equation}
which can be derived from \eqref{Wq} and \eqref{eq:Zv}
by integration by parts.

(ii) The map $v\mapsto v^{-1} Z^{(q,v)\prime}(x)$ is completely
monotone\footnote{A function $f\dvtx(a,\infty)\to\mbb R_+\setminus
\{0\}$,
$a\in\mbb R$,
is \emph{completely monotone} if $(-1)^{k-1}f^{(k)}(x)\ge0$ for
all $k\in\mbb N$ and $x>a$, where $f^{(k)}$ denotes the $k$th
derivative with respect to $x$. }
on $(\Phi(q),\infty)$, for any $x>0$.
That this is the case follows from the observation that $v\mapsto
v^{-1} Z^{(q,v)}(x)$ is the Laplace transform of some measure on $\mbb
R_+$ which is shown next. From the definition of
$Z^{(q,v)}$ it follows that the derivative $Z^{(q,v)\prime}(x)$ at
$x>0$ satisfies
\[
Z^{(q,v)\prime}(x) = v Z^{(q,v)}(x) + \bigl(q-\psi(v)\bigr)
W^{(q)}(x).
\]
Inserting the forms of the Laplace transforms of $W^{(q)}|_{\mbb R_+}$
and $Z^{(q,v)}|_{\mbb R_+}$ [given in~\eqref{Wq}
and \eqref{eq:LTZ}], it follows
%
\begin{eqnarray}
\label{eq:ltZtilde}  \mc L Z^{(q,v)\prime}(\theta) &=& \frac{q}{\psi(\theta) - q}
\nonumber
\\[-8pt]
\\[-8pt]
\nonumber
&&{} +
\frac{\theta v}{\psi(\theta) - q} \biggl[ \frac{\sigma^2}{2} + \int_{0}^\infty
\frac{\mathrm{e}^{-\theta y} - \mathrm
{e}^{-vy}}{v-\theta} \ovl \nu(y)\,\td y \biggr].
\end{eqnarray}
Inversion of the Laplace transform in \eqref{eq:ltZtilde}
and the observation
\[
\int_0^\infty\frac{\mathrm{e}^{-\theta y} - \mathrm
{e}^{-vy}}{v-\theta} \ovl\nu (y)\,\td
y = \int_0^\infty\int_{0}^\infty
\mathrm{e}^{-\theta s- v
t}\ovl\nu(s+t)\,\td t\,\td s,
\]
yield the following
expression for $v^{-1} Z^{(q,v)\prime}(x)$ at any $x>0$:
\[
\frac{q}{v} W^{(q)}(x) + \frac{\sigma^2}{2}W^{(q)\prime}(x) +
\int_{0}^\infty\int_{[0,x]}
\mathrm{e}^{-v t}\ovl\nu(x-y+t) W^{(q)}(\td y) \,\td t.
\]
By inspection it follows that, for any $x>0$, the function $v\mapsto
v^{-1} Z^{(q,v)\prime}(x)$
is the Laplace transform of a measure on $[0,\infty)$, which
implies the stated complete monotonicity.
\end{Rem}

\section{Single dividend-band strategies}\label{aux}
The analysis of various strategies starts with the case
of single dividend-band strategies.
In the absence of transaction costs
such a barrier strategy at level $b=(b_-,b_+)$, denoted by $\pi_b$,
specifies to pay out the minimal amount of dividends
to keep the reserves $U^b:= U^{\pi_b}$ below the level $b_+=b_-$,
while, in the case $K>0$, $\pi_b$ prescribes to
pay out a lump-sum $b_+-b_->0$ each time that the reserves
$U^b$ reach the level $b_+$. More formally, in the cases
$K=0$ and $K>0$ the forms of the strategy $\pi_b = \{D^b_t, t\in\mbb
R_+\}$ are
given by \eqref{Deebee} [with $b=b_+=b_-$]
and by
\[
D_t^{b} = \bigl(U_0^b-b_-\bigr)
+ (b_+-b_-) N^b_t,\qquad N_t^b=\#\bigl
\{s\in(0,t]\dvtx U_{s-}^b = b_+\bigr\},\qquad t\in\mbb R_+,
\]
respectively.
As a consequence, it follows that the value $v_{b}(x):=v_{\pi_b}(x)$
associated to the
single dividend band strategy $\pi_b$ at a nonzero level $b$ when
$X_0$ is equal to $x$ is given by
\[
v_{b}(x) = \E_x \biggl[\int_0^{\tau_b}
\mathrm{e}^{-qt}\mu _K^{b}(\td t) +
\mathrm{e}^{-q\tau_b}w\bigl(U^{b}_{\tau_b}\bigr) \biggr],
\]
with $\mu_K^{b}:= \mu_K^{\pi_{b}}$, $U^b:=U^{\pi_b}$ and
$\tau^{b}=\tau^{\pi_b}=\inf\{t\in\mbb R_+\dvtx U^{b}_t<0\}$.
The function $v_b$ can be
expressed in terms of the homogeneous and inhomogeneous scale
functions $W^{(q)}$ and $F_w$ as follows:

\begin{Prop}\label{prop:lp} For $b_+>b_-\ge0$ and $x\in[0,b_+]$ and
with $F=F_w$ it holds
%
\begin{eqnarray}
\label{eq:wapap}  v_{b}(x) &=& \cases{w(x), &\quad $ x < 0,$\vspace*{2pt}
\cr
W^{(q)}(x) G(b_-,b_+) + F(x), &\quad $x\in[0,b_+],$\vspace*{2pt}
\cr
x - b_+ +
v_b(b_+), & \quad $x > b_+,$}
\\
\label{eq:Gw}\qquad  G(b_-,b_+) &:=& \cases{\displaystyle\frac{b_+-b_- - K -
(F(b_+) - F(b_-))}{W^{(q)}(b_+)-W^{(q)}(b_-)}, &\quad  $K>0, b_+>b_-$,
\vspace*{2pt}
\cr
\displaystyle\frac{1 - F'(b_+)}{W^{(q)\prime}(b_+)}, &\quad  $K=0, b_+=b_-$.}
\end{eqnarray}
\end{Prop}

\begin{Rem}
Note that in the case $K>0$ and $X_0=x>b_+$ the strategy $\pi_b$
prescribes an immediate lump-sum dividend payment of size
$x-b_-$, which is in agreement with the value $v_b(x)$ for $x>b_+$,
\[
v_b(b_+) = v_{b}(b_-) + b_+-b_- - K \Rightarrow
v_b(x) = x-b_- - K + v_b(b_-),\qquad x>b_+.
\]
\end{Rem}

\begin{pf*}{Proof of Proposition~\ref{prop:lp}}
Consider the case $K>0$. Since
no dividend payment
takes place before $X$ reaches the level $b_+$ it follows that
$\{X_t, t\leq T_{0,b_+}\}$ and $\{U^{b_+}_t, t\leq\tau^{\pi_b}\}$
have the same law. The strong Markov property of $X$ and the absence of
positive jumps then yield that for $x\in[0, b_+]$ $v_b(x)$ is equal to
\begin{eqnarray*}
&&\E_x \bigl[\mathrm{e}^{-q T^+_{b_+}}\bigl(v_b(b_-) +
\Delta b - K\bigr)\mbf1_{\{T^+_{b_+} < T_0^-\}} \bigr] + \E_x \bigl[
\mathrm{e}^{-q T_0^-}w(U_{T_0^-}) \mbf1_{\{T^+_{b_+} >
T_0^-\}
} \bigr]
\\
&&\qquad= \frac{W^{(q)}(x)}{W^{(q)}(b_+)}\bigl[v_b(b_-) + \Delta b - K\bigr] +
\biggl[F(x) - F(b_+) \frac{W^{(q)}(x)}{W^{(q)}(b_+)} \biggr],
\end{eqnarray*}
with $F = F_{w}$, where the second line follows from Proposition~\ref
{prop:two} (applied with $w\equiv0$ and
with $\d=0$). Evaluating the expression in the display
at $x=b_-$, solving the resulting linear equation for
$v(b_-)$ and inserting the result yields the stated form.
The case $K=0$ follows by a similar line of
reasoning, using \eqref{d:wYta} in Proposition~\ref{prop:two}.
\end{pf*}

Next the candidate optimal levels are described.
The form of $G$ suggests to define the level $b^*=(b^*_-, b^*_+)$
as a maximizer of $G(x,y)$ over all $x,y\ge0$ in the case $K>0$, and
similarly, to define $b^*_+$ as a maximizer of $G(x,x)$ over all
$x\ge0$ in the case $K=0$.

\begin{Rem}\label{rem:gopt} Observe that in the case $K>0$ and $G$
is $C^1$,
the partial
right derivatives of $G(x,y)$ are given by
%
\begin{eqnarray}
\label{eq:GH}
\frac{\partial G}{\partial x}(x,y) &=& \frac{W^{(q)\prime}(x)}{W^{(q)}[x,y]}
\bigl[G(x,y) - G^\#(x)\bigr],
\nonumber
\\[-8pt]
\\[-8pt]
\nonumber
\qquad\frac{\partial G}{\partial y}(x,y) &=& -
 \frac{W^{(q)\prime}(y)}{W^{(q)}[x,y]}\bigl[G(x,y) - G^\#(y)\bigr],\qquad
G^\#(x):=\frac{1- F'(x)}{W^{(q)\prime}(x)},
\end{eqnarray}
and with $W^{(q)}[x,y]:= W^{(q)}(y)-W^{(q)}(x)$.
Therefore, in this case, an
interior maximum $(x^*,y^*)$ will satisfy $G(x^*,y^*)= G^\#(x^*) =
G^\#(y^*)$, and a candidate optimum may be found by fixing $d=y-x$,
and optimizing the left endpoint $x(d)$ for fixed $d$
[graphically, this would amount to determining the highest value
of the function $G^\#$ where the ``width'' $y(d)-x(d)$ of the
function $G^\#$ is $d$].
\end{Rem}

In the case $K>0$, fix
therefore $d>0$, and
let
%
\begin{equation}
\label{aodd} b^*(d) = \sup\bigl\{b\ge0\dvtx G(b, b+d)\geq
G(x,x+d)\ \forall x\ge
0\bigr\}
\end{equation}
denote the last global maximum of $G(x,x+d)$.

Define next $d^*$ to be the last global maximum of $G(b^*(y),b^*(y) +y)$
\[
d^* = \sup\bigl\{d\ge0\dvtx G\bigl(b^*(d),b^*(d) + d\bigr)\geq G
\bigl(b^*(y),b^*(y) +y\bigr)\ \forall y\ge0\bigr\},
\]
where $\inf\varnothing=+\infty$.

The candidate optimal levels are then defined as follows:
%
\begin{equation}
\label{eq:astar} b^*=\bigl(b^*_{-},b^*_{+}\bigr)\qquad
\mbox{with } b^*_{-} = b^*\bigl(d^*\bigr), b^*_{+} = b^*
\bigl(d^*\bigr) + d^*.
\end{equation}
In the absence of transaction cost ($K=0$), set
%
\begin{equation}
\label{eq:bK0} b_{+}^* = b_{-}^* = \sup\bigl\{b\ge0\dvtx
G^\#(b) \ge G^\#(x)\ \forall x\ge0\bigr\}.
\end{equation}

\begin{Thm}\label{thm:cm}
It holds $b_+^*<\infty$ and
%
\begin{equation}
\label{eq:v*b} v_{*}(x) = W^{(q)}(x)G^\# \bigl(b^*_+ \bigr)
+ F(x),\qquad x\in\bigl[0,b^*_+\bigr],
\end{equation}
where $F=F_w$.
In particular, it is optimal to adopt the
strategy $\pi_{b^*}$ while the reserves are not larger than $b^*_+$.
\end{Thm}

The proof rests on the following auxiliary result
that concerns explicit expressions linking the operator ${}_a {\Gamma
}_{\infty} ^w$ with
the function $G$ and the scale functions $F_w$ and~$W^{(q)}$. This relation
is also deployed in the formulation of necessary and sufficient
optimality conditions for optimality of band policies in Sections \ref
{ssec:singo}--\ref{sec:multib}.

\begin{Lemma}\label{prop:key}
Let $c>0$, and for any $b_+\ge b_-\ge0$ (with $b_+\neq b_-$ in the
case $K>0$)
define
\[
{} J_{v_b}\dvtx\mbb R_+\setminus\{0\}\to\mbb R, \qquad
J_{v_b}(y) = \bigl({}_{b_+} \Gamma^{v_b}_\infty
v_{b}\bigr) (y),\qquad y> 0.
\]
\begin{longlist}[(ii)]
\item[(i)] The following identity holds true:
%
\begin{eqnarray}
\label{eq:key}
&&W^{(q)\prime}\bigl(b^*_++c\bigr)\bigl[G
\bigl(b^*_-,b^*_++c\bigr) - G\bigl(b^*_-,b^*_+\bigr)\bigr]\nonumber \\
&&\qquad=\int
_{[0,c)} {} _{b^*_+}J_{v_{b^*}}(c-y)
W^{(q)}(\td y)
\\
&&\qquad= v_{b^*,-}'\bigl(b^*_+\bigr) - F'_{{}_{b^*_+} v_b}(c).\nonumber
\end{eqnarray}
In particular, it holds
%
\begin{equation}
\label{eq:idFdW} \int_{[0,c)} {} _{b^*_+}J_{v_{b^*}}(c-y)
W^{(q)}(\td y) < 0\qquad \forall c>0,
\end{equation}
and the functions $y\mapsto G(b^-,y)$
and $y\mapsto G^\#(y)$ are decreasing
for all $y$ sufficiently large.

\item[(ii)] Denoting $G_{b_-}(x):= G(b_-,x)$,
the Laplace transform of the function $g\dvtx\mbb R_+\setminus\{0\}
\to
\mbb
R$ given by
$g(x)={} _{b_+}J_{v_b}(x)$ is equal to
%
\begin{equation}
\label{eq:Xi2} \mc L g(\theta) = + \frac{\mathrm{e}^{\theta b_+}}{\theta} \int_{(b_+,\infty)}
\mathrm {e}^{-\theta z} Z^{(q,\theta)\prime}(z) G_{b_-}(\td z),\qquad
\theta>\Phi(q).
\end{equation}
In particular, $g$ is nonpositive precisely if
$\theta\mapsto- \mc L g(\theta+\Phi(q))$ is completely monotone.
\end{longlist}
\end{Lemma}

%
\begin{Rem}
The integral in \eqref{eq:Xi2} is to be interpreted as a
Lebesgue--Stieltjes integral.
This follows as a consequence of the form of $G_{b_-}$ and the fact that
the functions $W^{(q)}$ and $1/W^{(q)\prime}$ are of bounded variation
(which follows in turn
as $W^{(q)}$ is increasing and $W^{(q)\prime}$ is logconcave).
\end{Rem}

The proof of Lemma~\ref{prop:key} is given in  Appendix \ref{pf:pr:key}.

\begin{pf*}{Proof of Theorem~\ref{thm:cm}}
\emph{$b^*_+$ is finite, and the supremum is attained:}
Note that, for any $x>0$, it holds $G^\#(x)\ge G^\#(x-)$,
by virtue of the form \eqref{eq:GH} of $G^\#(x)$, and
the inequalities $W^{(q)\prime}(x)\ge W^{(q)\prime}_-(x)$
[from \eqref{eq:WqD}] and $F^{\prime}(x)\ge F^{\prime}_-(x)$
[from~\eqref{eq:Fwd0}], where $W^{(q)\prime}_-(x), F^{\prime}_-(x)$
denote the left-derivatives at $x$.
In view of the facts that the map $x\mapsto G^\#(x)$ defined in \eqref
{eq:Gw} is
right-continuous and monotone decreasing for all $x$ sufficiently
large (Lemma~\ref{prop:key}),
it then follows that
there exists an $x^*\in\mbb R_+$ such that $\sup_{x\ge0} G^\#(x)=G^\#(x^*)$.
In the case that $K$ is strictly positive, $G$
attains its maximum at some $(x^*,y^*)\in Q:=\{(z_1,z_2)\in\mbb
R^2\dvtx0\leq z_1<z_2\}$, since
(a) $G(x,y)$ is continuous at any $(x,y)$ in $Q$, (b) monotone
decreasing for $y$ sufficiently large and fixed $x$
[Proposition~\ref{prop:key}(iii)], (c) tends to minus infinity if
$y\searrow x$ and
(d) tends to the constant $\kappa_w$ in \eqref{eq:kappa} if $|x|+|y|$
tends to infinity such that
$x<y$.

\emph{Verification of optimality:}
Assume for the moment that the function $h\dvtx\mbb R_+\to\mbb R$
defined by
the right-hand side of \eqref{eq:v*b} is a supersolution in the sense of
Definition~\ref{def:sss}. Under this assumption $h$ dominates the
value-function $v_*$ (by Proposition~\ref{thm:repg}). In
fact, since $h(x)$ is equal to the value $v_{b^*}(x)$ of the
strategy $\pi_{b^*}$ for any level $x$ of initial reserve smaller
or equal to $b^*_+$, the local verification theorem, Theorem~\ref
{cor:repg}(i), implies that $h(x)$
is equal to the optimal value $v_*(x)$ for all $x\in[0,b^*_+]$.

Next it is shown that $h$ is a supersolution by verifying the following
two facts:
(a) $\mathrm{e}^{-q(t\wedge
T_0^-)}h(X_{t\wedge T_0^-})$ is a martingale, and (b) $h$
satisfies the inequality
\[
h(x) - h(y) \ge x-y - K \qquad\mbox{for any $0\leq y< x$}.
\]
Fact (a) follows from the martingale properties of
$F_w$ and $W^{(q)}$ (see Proposition~\ref{prop:GS}),
while (b) follows on account of the
definitions of $b^*$ and $G^\#$. Indeed, if $K=0$ and $x>0$,
$h'(x) = W^{(q)\prime}(x) G^\#(b^*) - F_w'(x)$ is bounded below by
%
\begin{equation}
\label{disp1} W^{(q)\prime}(x) G_*(x) - F_w'(x)=1,
\end{equation}
while, if $K>0$ and $x>y>0$,
$h(x)-h(y)=(W^{(q)}(x)-W^{(q)}(y))G(b^*_-,b^*_+) -
F_w(x) + F_w(y)$ is bounded below by
%
\begin{eqnarray}
\label{disp2} h(x)-h(y) &\ge& \bigl(W^{(q)}(x)-W^{(q)}(y)
\bigr)G(y,x) - F_w(x) + F_w(y)
\nonumber
\\[-8pt]
\\[-8pt]
\nonumber
& =& x - y -K.
\end{eqnarray}
Displays \eqref{disp1} and \eqref{disp2}
imply $h(x) - h(y)\ge x-y-K$ for any $K\ge0$
and $x,y\ge0$ with $x\ge y$. This completes the proof.
\end{pf*}

%
\section{Two-band strategies and a mixed optimal stopping/control
problem}
\label{aux2}
The policy $\pi_{b^*}$ considered in the previous section may be
optimal for any level of the reserves,
and not just for small levels as shown in Theorem~\ref
{thm:cm}---necessary and sufficient conditions
for this to be the case are given in
Section~\ref{ssec:singo}.
In this section the complementary case is considered that it is optimal
to have a \emph{second} dividend band.
The problem of finding the optimal levels of the second dividend band
differs from the single-band optimization problem in the following two respects:
\begin{longlist}[(ii)]
\item[(i)] at any time $t$ prior to the time of ruin it is possible to
make a lump-sum payment to
bring the reserves down to the level
$b^*_-$ defined in \eqref{eq:astar}, yielding a pay-off of $U_t-b^*_- +
v_{b^*}(b^*_-) - K$, and
\item[(ii)] it will not be optimal to place a dividend band at levels
close to $b^*_+$.
\end{longlist}
The observation in (i) in combination with the dynamic programming
principle (Proposition~\ref{prop:mart})
and Theorem~\ref{thm:cm} yield the representation
%
\begin{equation}
\label{eq:ooptst}  v_*(x) = \sup_{\pi\in\Pi,\tau\in\mc T} \E_x
\biggl[ \int_{[0, \tau\wedge\tau)}\mathrm{e}^{-qt}
\mu_K^\pi(\td t) + \mathrm{e}^{-q (\tau^\pi_{b^*}\wedge\tau
)}v_{b^*}
\bigl(U^\pi _{\tau^\pi
_{b^*}\wedge\tau} \bigr) \biggr],
\end{equation}
where $\tau^\pi_{b^*}=\inf\{t\ge0\dvtx U^\pi_t < b^*_+\}$.
This section is devoted to a stochastic control problem that is closely
related to \eqref{eq:ooptst},
$V_*^f(x) = \sup_{\pi\in\Pi,\tau\in\mc T} V^f_{\tau,\pi}(x)$, where
%
\begin{equation}
\label{eq:optstop} V^f_{\tau,\pi}(x) = \E_x
\biggl[ \int_{[0,\tau^\pi\wedge\tau)}\mathrm{e}^{-qt}
\mu_K^\pi(\td t) + \mathrm{e}^{-q (\tau^\pi\wedge\tau  )}f
\bigl(U^\pi _{\tau^\pi
\wedge\tau} \bigr) \biggr],
\end{equation}
where, as before $\tau^\pi=\inf\{t\ge0\dvtx U^\pi_t < 0\}$, and
$f\dvtx\mbb R\to\mbb R$ is assumed to satisfy the following conditions:
%
\begin{eqnarray}
\label{eq:as1} && \mbox{ $f|_{\mbb R_+}$ is given by $f(x)=x+c$ for $x\in
\mbb R_+$, for some $c\in\mbb R$,}
\\
&& f_-'(0)\ge1, \label{eq:as1b}
\\
\label{eq:as2} && \mbox{$J_{\bar w}(u):={}_0 {
\Gamma}_\infty^{\bar w} f(u) >0$ for some $u>0$, with $\bar w =
f|_{\mbb R_-}$},
\\
\label{eq:as3} && \mbox{for all $c\in\mbb R_+\setminus\{0\}$, $ \int
_{[0,c)}J_{\bar w}(c-y)W^{(q)}(\td y) < 0 $}.
\end{eqnarray}
It will be shown that, under \eqref{eq:as2}, it is not optimal to stop
immediately ($V_*^f\not\equiv f$),
while, under \eqref{eq:as3}, the dividend barrier strategy with level
at 0 is not optimal
($V_*^f\not\equiv V^f_{\tau^\pi,\pi_0}$). In particular, in the setting
of the stochastic control problem
in~\eqref{eq:ooptst} conditions
in \eqref{eq:as1}--\eqref{eq:as3} are satisfied:

\begin{Lemma}\label{lem:f} If it holds $v_{\pi_{b^*}}(x) < v_*(x)$ for
some $x>b^*_+$,
then the function $f\dvtx\mbb R\to\mbb R$ defined by $f(x) = v_{b^*}(b^*_++x)$
satisfies the stated conditions in \eqref{eq:as1}--\eqref{eq:as3}.
\end{Lemma}

\begin{pf}
First, note that the conditions in \eqref{eq:as1}--\eqref{eq:as1b} hold
since $v_{b^*}|_{[b^*_+,\infty)}$ is affine with unit slope
and $v_{b^*,-}'(b^*)$ is larger or equal to one (with equality when
$W^{(q)}$ and $F_w$
are differentiable at $b^*$). Also, condition \eqref{eq:as3} holds by
\eqref{eq:idFdW} in Lemma~\ref{prop:key}.
Furthermore, it is shown in Theorem~\ref{thm:cm2} in Section~\ref
{ssec:singo} that
if condition \eqref{eq:as2} was not satisfied, then
$v_{b^*}=v_*$, which would be in contradiction
with the assumed existence of an $x$ larger than $b_+^*$ satisfying
$v_{b^*}(x) < v_*(x)$.
\end{pf}

Next a candidate optimal policy is specified for the mixed optimal
stopping/optimal control problem in \eqref{eq:optstop}.
Strategies for this optimization problem consist of pairs $(\tau,\pi)$
of an $\mbf F$-stopping time $\tau$ and
a policy $\pi$ from the set $\Pi$. The discussion at the beginning of
the section [especially item (ii)] in conjunction with Lemma~\ref
{lem:f} suggests
to consider candidate optimal strategies of the form $(\tau^{\pi
_b}_a,\pi^b)$, $a<b_+$:
such policies specify to pay out dividends according to
a single dividend-band strategy $\pi_b$ at levels $(b_-, b_+)$ until
the first moment $\tau_a^{\pi_b}=\inf\{t\ge0\dvtx U^{\pi_b}_t<a\}$
that $U^{\pi_b}$ falls below the level $a > 0$ at which moment one
should stop.
Another strategy that is worth considering in the case $K>0$
is to refrain from paying dividends until the first moment that the
reserves process exits a finite interval $[a,b_+]$ and to stop then;
such strategies are denoted by
$(\pi^\varnothing, T_{a,b_+})$ for $a<b_+$.
The value functions associated to the strategies
$(\tau^{\pi_b}_a,\pi^b)$ and $(\pi^\varnothing, T_{a,b_+})$
are given by
\[
V^{f}_{a,b_-, b_+}(x) = \E_x \biggl[\int
_{[0, \tau_a^{\pi
_{b}})}\mathrm{e}^{-qt}\mu _K^{b}(
\td t) + \mathrm{e}^{-q\tau_a^{\pi_{b}}}f \bigl(U^{b}_{\tau_a^{\pi
_{b}}} \bigr)
\biggr],
\]
and $V^{f,\varnothing}_{a, b_+}(x) = \E_x [\mathrm
{e}^{-qT_{a,b_+}}f
(X_{T_{a,b_+}}  )  ]$,
with $\mu_K^{b} = \mu_K^{\pi_{b}}$. In the following result,
which can be derived by a line of reasoning that is similar to the one
used in the proof of Proposition~\ref{prop:lp}, the functions
$V^{f}_{a,b_-,b_+}$ and
$V^{f,\varnothing}_{a,b_+}$
are explicitly expressed in terms of
scale functions and the families of functions $(y,z)\mapsto
G_{f}^{(a)}(y,z)$, $G_{f,\varnothing}^{(a)}(y,z)$, $a\ge0$,
that are defined as follows:
%
\begin{eqnarray}
\label{eq:Gfabb}\qquad G_{f}^{(a)}(b_-,b_+) &= &\cases{
\displaystyle\frac{
b_+ - b_- - K - F^{(a)}[b_--a,b_+-a]}{W^{(q)}[b_--a, b_+-a]}, & \quad
$K>0,$ \vspace*{2pt}
\cr
\displaystyle G_{f,\#}^{(a)}(b_+):= \frac{1 - F^{(a)\prime}(b_+-a)}{
W^{(q)\prime}(b_+-a)}, &\quad $K=0,$}
\\
 G_{f,\varnothing}^{(a)}(b_+) &=& \frac{f(b_+) -
F^{(a)}(b_+-a)}{W^{(q)}(b_+-a)}, \label{eq:Gfem}
\end{eqnarray}
where $F^{(a)} = F_{_a f}$ is the Gerber--Shiu function for payoff
$_a f= f(a+\cdot)$, $F^{(a)}[x,y] = F^{(a)}(y)-F^{(a)}(x)$ and, as before,
$W^{(q)}[x,y] = W^{(q)}(y) - W^{(q)}(x)$.

\begin{Prop}\label{prop:lp2}
For any $b_-,b_+, a\in\mbb R_+$ satisfying $b_+\ge b_-\ge a$
the following representations hold true:
\begin{eqnarray*}
V^{f}_{a,b_-, b_+}(x) &=& \cases{F^{(a)}(x-a) = f(x), &\quad
$x \in[0, a),$\vspace*{2pt}
\cr
W^{(q)}(x-a)G^{(a)}_f(b_-,b_+)
+ F^{(a)}(x-a), &\quad $x\in[a,b_+],$\vspace *{2pt}
\cr
x-b_+ +
V^f_{a,b_-, b_+}(b_+), & \quad $x \in(b_+,\infty)$;}
\\
V^{f,\varnothing}_{a, b_+}(x) &=& \cases{ F^{(a)}(x-a) = f(x), &\quad
$x \notin[a, b_+]$,\vspace*{2pt}
\cr
W^{(q)}(x-a)G^{(a)}_{f,\varnothing}(b_+)
+ F^{(a)}(x-a), & \quad $x\in[a,b_+]$.}
\end{eqnarray*}
\end{Prop}

Next the candidate optimal levels are described.
Focusing first on the case that dividends are paid and fixing the level
$a$ for the moment,
and similarly as in the case of the single dividend-band
strategies, let $\b_f^*(a)= (\beta^*_{f,-}(a), \beta
^*_{f,+}(a)  )$ denote
the (largest) maximizer of the function $G_f^{(a)}$.
In the case $K>0$ we set
\begin{eqnarray*}
\beta^*_{f,-}(a) &=& \b_f^*\bigl(a,\d_f^*(a)
\bigr),\qquad \beta^*_{f,+}(a) = \b_f^*\bigl(a,
\d_f^*(a)\bigr)+\d_f^*(a),
\\
\b_f^*(a,d) &=& \sup \bigl\{b\ge a\dvtx G^{(a)}_f(b,b+d)
\ge G^{(a)}_f(x,x+d)\ \forall x\ge0 \bigr\},
\\
\d_f^*(a) &=& \sup \bigl\{d\ge0\dvtx G_f^{(a),*}(d)
\leq G_f^{(a),*}(y)\ \forall y\ge0 \bigr\},
\end{eqnarray*}
with $G_f^{(a),*}(d):= G_f^{(a)} (\b_f^*(a,d),\b
_f^*(a,d)+d  )$,
while, in the case $K=0$, we define
\[
\b^*_{f,+}(a) = \b_{f,-}^*(a) = \b^*_{f,\#}(a): =
\sup \bigl\{b\ge a\dvtx G_{f,\#}^{(a)}(b) \ge
G_{f,\#}^{(a)}(x)\ \forall x\ge0 \bigr\}.
\]

The candidate optimal specification $\a^*_f$ of the stopping level $a$
and the candidate optimal
level $\beta_f^*$ are given by
%
\begin{eqnarray}
\label{eq:astar2}  \a_f^* &=& \inf \bigl\{a\ge0\dvtx
G_f^{(a,*)} \bigl(\d_f^*(a) \bigr) > 0 \bigr\}\qquad
\mbox {in the case $K>0$},
\\
\label{eq:astar22}  \a_f^* &=& \inf \bigl\{a\ge0\dvtx
G_{f,\#}^{(a)} \bigl(\b_{f,\#}^*(a) \bigr) > 0 \bigr\}\qquad
\mbox{in the case $K=0$},
\\
\label{eq:bstar2}  \b^*_{f} &=& \bigl(\b^*_{f,-},
\b^*_{f,+} \bigr), \qquad\b^*_{f,-} = \b ^*_{f,-} \bigl(
\a^*_f \bigr), \qquad\b^*_{f,+} = \b^*_{f,+} \bigl(
\a^*_f \bigr).
\end{eqnarray}

Next consider the strategy to continue without paying dividends and
stop upon exiting a finite interval.
It will turn out that in the case $K=0$ such a strategy is never
optimal; see Remark~\ref{rem:aa}.

In the case $K>0$ define
%
\begin{eqnarray}
\label{eq:bstare}  \b_{f,\varnothing}^*(a) &=& \sup \bigl\{b\ge a\dvtx
G^{(a)}_{f,\varnothing
}(b)\ge G_{f,\varnothing}^{(a)}(x)
\ \forall x\ge0 \bigr\},
\\
 \a^*_{f,\varnothing} &=& \inf \bigl\{a\ge0\dvtx G_{f,\varnothing
}^{(a)}
\bigl(\b ^*_{f,\varnothing}(a) \bigr) > 0 \bigr\},\qquad \b^*_{f,\varnothing} =
\b^*_{f,\varnothing}\bigl(\a^*_{f,\varnothing
}\bigr).\label
{eq:bstare2}
\end{eqnarray}
The levels $\b_{f,+}^*$, $\beta_{f,\varnothing}^*$ and $\a
^*_{f,\varnothing
}$ given above are
finite and strictly positive.

\begin{Lemma}\label{lem:a}
Suppose that $f$ satisfies the conditions in \eqref{eq:as1}--\eqref{eq:as3}
and denote $\bar w = f|_{\mbb R_-}$.
\begin{longlist}[(ii)]
\item[(i)] $K=0$: $0<\a^*_f\leq\b^*_{f,+}<\infty$ and
$ G^{(\a^*_f)}_{f,\#}(\b^*_{f})=0$, and
$_0\Gamma^{\bar w}_\infty f(u) \leq0$ for all $u\in (0,\alpha
^*_f  )$.

Furthermore, if $X$ has unbounded variation, it holds $\a^*_f<\b^*_{f,+}$.

\item[(ii)] $K>0$: $0<\a^*_{f,\varnothing}\leq\b^*_{f,\varnothing}<\infty
$ and
$G_{f,\varnothing}^{(\a^*_{f,\varnothing})}(\beta^*_{f,\varnothing
})=0$, and
it holds\break $_0\Gamma^{\bar w}_\infty f(u) \leq0$ for all $u\in
(0,\alpha^*_{f,\varnothing}  )$.

Furthermore, if it holds in addition $\a^*_f<\infty$, then $0<\a
^*_f<\b
^*_{f,+}<\infty$ and
$G_f^{(\a^*_f)} (\beta^*_{f}  )=0$.
\end{longlist}
\end{Lemma}

%
\begin{Rem}[(Smooth and continuous fit)] The choice of $\a^*_f$
coincides with what would be obtained
by applying
the \emph{principles of
continuous} and \emph{smooth fit} from the theory of optimal stopping
(see Peskir and Shiryaev \cite{PesShir}, Chapter IV.9), which suggest
that in the
mixed optimal stopping/stochastic control problem \eqref{eq:optstop}
it can be expected
that $V^f$ be \emph{continuous}/\emph{continuously
differentiable} at a level $\a^*_f$ if $\a^*_f$ is irregular/regular for
$(-\infty,\a^*_f)$ for $X$, respectively, where $\pi_*$ denotes the
optimal strategy.
Since it is well-known that $\a^*_f$ is regular for
$(-\infty,\a^*_f)$ for $X$ if and only if $X$ has unbounded variation,
this heuristic yields
\[
\cases{ \mbox{$\a^*_f$ satisfies $V^{f\prime}_{\a^*,\b^*}
\bigl(\a^*_f +\bigr)=f'\bigl(\a ^*_f-
\bigr)$}, &\quad $\mbox{if $X$ has unbounded variation,}$\vspace*{2pt}
\cr
\mbox{$
\a^*_f$ satisfies $V^f_{\a^*,\b^*}\bigl(
\a^*_f\bigr)=f\bigl(\a^*_f\bigr)$}, & \quad $\mbox {if $X$ has
bounded variation.}$}
\]
The first equation in the display is equivalent to the expression in
\eqref{eq:astar2} in view of the form of $V^f_{a,b}$ and the facts
(i) $F_{_a
f}'(0)=f_-'(a)$ for any $a>0$ and (ii)
$W_+^{(q)\prime}(0)\in(0,\infty]$. The second equation in the display
can also be
equivalently expressed as~\eqref{eq:astar2}, in view of (i$'$) the
form of $V^f_{a,b_-, b_+}$ in Proposition~\ref{prop:lp2} and (ii$'$) the
fact that $W^{(q)}(0)$
is strictly positive precisely if $X$ has
bounded variation. A similar remark holds true for the level $\a
^*_{f,\varnothing}$.
\end{Rem}

\begin{Rem}\label{rem:aa}
(i) In the case $K=0$ it is straightforward to verify that
any strategy of the form $(\pi_\varnothing, T_{a,b_+})$, for $a,b\in
\mbb
R_+$ with
$0< a<b_+$, is not optimal
[indeed, the minimal slope of the value function $u$ of such a strategy
is smaller than one,
since $u$ satisfies $u(b_+)-u(0)=b_+$, given that $u(b_+)=f(b_+)$,
$u(0)=f(0)$ and
$f$ is affine with unit slope].\vspace*{-6pt}
\begin{longlist}[(ii)]
\item[(ii)] In the case $K>0$ and $\alpha^*_{f,\varnothing}<\alpha^*_f$,
the definition of $\alpha^*_f$, Proposition~\ref{prop:lp2} and
Lemma~\ref{lem:a}(ii) imply
\[
V(x):= V^{f,\varnothing}_{\a^*_{f,\varnothing},\beta
^*_{f,\varnothing}}(x) \ge V^f_{\a^*_{f,\varnothing},
\beta^*_{f}(\a^*_{f,\varnothing})}(x),\qquad
x\in\bigl[0,\beta^*_{f,\varnothing}\bigr].
\]
Note that the nonpositivity of $G_f^{(\a^*_{f,\varnothing})}(\beta
^*(\a
^*_{f,\varnothing}))$
implies $\mathtt d_{V}(x)\ge1$ for all $x>0$.

\item[(iii)] In the case $K>0$ and $\alpha^*_{f,\varnothing}\ge\alpha^*_f$ a
similar argument using
the definition of $\alpha^*_{f,\varnothing}$ in conjunction with
Proposition~\ref{prop:lp2} and Lemma~\ref{lem:a}(ii) implies
\[
V^{f}_{\a^*_{f},\beta^*_{f}}(x) \ge V^{f,\varnothing}_{\a
^*_{f},\beta
^*_{f,\varnothing}(\a^*_{f})}(x),\qquad x
\in\bigl[0,\beta^*_{f}\bigr].
\]
\end{longlist}
\end{Rem}

\begin{pf*}{Proof of Lemma~\ref{lem:a}}
(i) Consider the function $\ovl G\dvtx\mbb R_+\to\mbb R$ defined by
$\ovl
G(a) = \sup_{b\ge0}G_{f,\#}^{(a)}(b)$. The fact that
$\a_f^*$ is positive and finite is a consequence of the intermediate
value theorem and
the following three assertions concerning $\ovl G$:
\begin{longlist}[(a)]
\item[(a)] $\ovl G(0)<0$;
\item[(b)] there exists an $a_0\in\mbb R_+\setminus\{0\}$ such
that $\ovl G(a_0)>0$;
\item[(c)] {the function $a\mapsto\ovl
G(a)$ is continuous at $a\in[0,a_0]$.}
\end{longlist}
Next these three assertions are verified.
Assertion (a) follows from the definition of $G^{(0)}$ in \eqref{eq:Gfabb},
the form of $F^{(a)\prime}$ [in \eqref{eq:Fwd0}]
and conditions \eqref{eq:as1b} and \eqref{eq:as3}.

To verify assertion (b) it suffices to find $a_0$ and $b$ with $a_0<b$
satisfying $G_{f,\#}^{(a_0)}(b)>0$, or equivalently $F^{(a_0)\prime}(b-a_0)<1$
(in view of the form of $G_{f,\#}^{(a_0)}$). By the form of
$F^{(a_0)\prime}$ and the fact $f'(a_0)\ge1$
it suffices to show $\int_{[0,b-a_0)}J_{\tilde w}(b-a_0-y)W^{(q)}(\td
y)>0$ with $\tilde w = {}_{a_0} f$
for some $a_0<b$, which is equivalent to the condition
$\int_{[0,b-a_0)}J_{\ovl w}(b-y)W^{(q)}(\td y)>0$ for some $a_0<b$,
as it holds $J_{\ovl w}(b-y) = J_{\tilde w}(b-a_0-y)$.

To see that the latter condition is satisfied, note that
right-continuity of the map $J_{\ovl w}$ and \eqref{eq:as2} imply that
there exists an interval $I=[u_-,u_+]$, with
$0<u_-<u_+$, such that $J_{\ovl w}(y)>0$ for all $y\in I$;
taking $a_0:=u_-$ and $b:=u_+$ it thus follows that the integral
$\int_{[0,b-a_0)}J_{\ovl w}(b-y)W^{(q)}(\td y)$ is strictly positive,
and the proof of assertion (b) is complete.

To verify that assertion (c) holds fix $a\ge0$, and note
$V^f_{a,\beta^*(a)}(x) = W^{(q)}(x)\times \ovl G(a) + F^{(a)}(x-a)$
for $x\in[a,\beta_+^*(a)]$. By reasoning analogous to the proof of
Theorem~\ref{thm:cm}
the following identity can be shown to hold:
\[
V^f_{a,\beta^*(a)}(x) = \sup_{(\pi,\tau)\in\Pi(\beta^*_+)} \E
_x \biggl[\int_{[0, \tau^\pi_{a}\wedge\tau]}\mathrm{e}^{-qt}\,
\td D^\pi_t + \mathrm{e}^{-q(\tau^\pi_{a}\wedge\tau)} f
\bigl(U^\pi_{\tau^\pi
_{a}\wedge
\tau}\bigr) \biggr],
\]
where $\Pi(\beta^*_+)$ is the set of the strategies $(\pi,\tau)$
that is such that the stochastic process $\{U^\pi_{t\wedge\tau},
t\in
\mbb R_+\}$ stays below the level
$\beta^*_+$.
Let $a_1, a_2\in\mbb R_+$ be such that $a_2<a_1<\min\{\beta
^*(a_1),\beta
^*(a_2)\}$
and fix $x_0\in(a_1, \min\{\beta^*(a_1),\beta^*(a_2)\})$. To show
the continuity
of $\ovl G(a)$ we show next that $V^f_{a_1,\beta^*(a_1)}(x_0) -
V^f_{a_2,\beta^*(a_2)}(x_0)$ tends to $0$ when $a_2-a_1\to0$.

By an application of the triangle inequality it follows that the difference
$\llvert V^f_{a_1,\beta^*(a_1)}(x_0) - V^f_{a_2,\beta
^*(a_2)}(x_0)\rrvert $
is bounded above by
%
\begin{eqnarray}
\label{eq:supa}&& \sup_{\pi\in\Pi} \E_{x_0} \biggl[\int
_{[\tau^\pi_{a_1}, \tau
^\pi_{a_2}]} \mathrm{e}^{-qt}\,\td D^\pi_t
+ \bigl\llvert \mathrm{e}^{-q\tau^\pi
_{a_2}}f\bigl(U^\pi_{\tau
^\pi_{a_2}}
\bigr) - \mathrm{e}^{-q\tau^\pi_{a_1}}f\bigl(U^\pi_{\tau^\pi_{a_1}}\bigr)
\bigr\rrvert \biggr].
\end{eqnarray}
Since $\P_{x_0}(U_{\tau^\pi_{a_1}}\in[a_2,a_1)) = \P_{x_0}(\tau
^\pi
_{a_1}< \tau^\pi_{a_2})$ converges to zero if
$a_1-a_2\searrow0$, it follows that also the random variable under
the expectation
tends to zero $\P_{x_0}$-a.s. if $a_1-a_2\searrow0$. Since this random
variable is dominated by an integrable random variable,
uniformly
for all $(\pi,\tau)\in\Pi(\beta^*_+)$, Lebesgue's dominated convergence
theorem implies that the
right-hand side of \eqref{eq:supa} tends to zero when
$a_1-a_2\searrow0$.
To see that the random variable is dominated, recall that $f$ is
affine, and
note that $\mathrm{e}^{-q\tau^\pi_{a_1}}D^\pi_{\tau^\pi
_{a_1}}\vee
\mathrm{e}^{-q\tau^\pi_{a_2}}D^\pi_{\tau^\pi_{a_2}} \vee
\int_{[\tau^\pi_{a_1},\tau^\pi_{a_2}]}\mathrm{e}^{-qt}\,\td D_t^\pi$
is bounded above by
\[
\int_{[0,\infty)}\mathrm{e}^{-qt}\,\td
D^\pi_t \leq\int_0^\infty
q\mathrm{e}^{-qt} D^\pi _t\,\td t \leq\int
_0^\infty q\mathrm{e}^{-qt} \ovl
X_t\,\td t
\]
with $\ovl X_t = \ovl X^0_t = \sup_{s\in[0,t]}X_s\vee0$,
which is equal to $\E_{x_0}[\ovl X_{\mbf e_q}] = x_0 + \Phi(q)^{-1}$,
where ${\mbf e}_q$ is an independent
exponential random time, and
\[
\E_{x_0} \bigl[\bigl\llvert \mathrm{e}^{-q\tau^\pi_{a}}X_{\tau^\pi
_{a}}
\bigr\rrvert \bigr] \leq \E_{x_0} \bigl[\mathrm{e}^{-q\tau^\pi_{a}}(\ovl
X_{\tau^\pi_{a}} -\unl X_{\tau^\pi_{a}}) \bigr] \leq2x_0 +
\E_{x_0}[\ovl X_{\mbf
e_q} - \unl X_{\mbf e_q}] <\infty,
\]
with $\unl X_t=\inf_{0\leq s\leq t}X_s\wedge0$, where the finiteness
follows from the bound
$\E_{x_0}[\unl X_{\mbf e_q}]\ge\E_{0}[\unl X_{\mbf e_q}] = \E
_{0}[X_{\mbf e_q}] - \E_0[\ovl X_{\mbf e_q}]$
(which follows from the Wiener--Hopf factorization) and the fact $\E
_{0}[X_{\mbf e_q}] = \psi'(0)/q$.

The finiteness of $\b_{f,+}^*(\a_f^*)$ follows by a line of reasoning
that is analogous to the one that was used in the proof of Theorem~\ref
{thm:cm}, while
the relation $\beta_{f,+}^*(\a_f^*)\ge\a^*_f$
follows by definition of $\beta_{f,+}^*(\a^*_f)$. Finally, in the case
$K=0$ and $\{\sigma^2>0$ or $\nu_{0,1}=\infty\}$ the equality
$\alpha
^*=\b^*_+(\a^*)$
would imply that $V^f_{\a^*,\b^*}\equiv f$; however, since there exists
a $u$
such that $_0 \Gamma_{\infty}^f f(u)>0$ by \eqref{eq:as2}, an argument
as above shows that,
for some $\a,\beta$, $V^f_{\a,\b}(x)>f(x)$ for $x\in(\a,\beta)$, which
yields a contradiction.
A similar argument shows $_0\Gamma^{\bar w}_\infty f(u) \leq0$ for all
$u\in(0,\alpha^*_f)$.

The proof of part (ii) is analogous to that of part (i), and is omitted.
\end{pf*}

The solution of the stochastic control problem in \eqref{eq:optstop}
for small levels of the reserves is given as follows:

\begin{Thm}\label{thmos}
Suppose that $f$ satisfies conditions \eqref{eq:as1}--\eqref{eq:as3}.
\begin{longlist}[(ii)]
\item[(i)] When either $K=0$ or $\{K>0$ and $\alpha^*_{f,\varnothing} \ge
\alpha
^*_{f}\}$,
it holds
$V^f_*(x)
= V^f_{\a^*_f,\beta^*_f}(x)$ for any $x\in [0,\beta
^*_{f,+}  ]$.
While the reserves are smaller than $\b^*_{f,+}$
it is optimal
to adopt the policy $ (\tau^{\pi_{\b^*}}_{\a^*},\pi_{\b
^*}  )$.

\item[(ii)] In the case $\{K>0$ and $\alpha^*_{f,\varnothing} < \alpha
^*_{f}\}$
it holds
$V^f_*(x) = V^{f,\varnothing}_{\a^*_{f,\varnothing},\beta
^*_{f,\varnothing}}(x)$
for any $x\in [0,\beta^*_{f,\varnothing}  ]$. While the
reserves are
smaller than
$\b^*_{f,\varnothing}$ it is optimal to adopt the policy
$ (T_{\a^*_{f,\varnothing},\beta^*_{f,\varnothing}},\pi
^\varnothing
  )$.

In particular, it holds
%
\begin{eqnarray}
V^f_*(x) = \cases{ f(x), &\quad $x\in[0,a^*),$\vspace*{2pt}
\cr
F^{(a^*)}\bigl(x-a^*\bigr), &\quad $x\in\bigl[a^*, b^*\bigr],$}
\end{eqnarray}
where $F^{(a^*)}=F_{{}_{a^*} f}$ and $ (a^*,b^*  )= (\a
^*_f,\beta
^*_{f,+}  )$
in the cases $K=0$ or $\{K>0$ and $\alpha^*_{f,\varnothing} \ge
\alpha
^*_{f}\}$, and
$ (a^*,b^*  ) =  (\alpha^*_{f,\varnothing}, \beta
^*_{f,\varnothing}  )$
in the case $\{K>0$ and $\alpha^*_{f,\varnothing} < \alpha^*_{f}\}$.
\end{longlist}
\end{Thm}

The proof of Theorem~\ref{thmos} rests an auxiliary result concerning
the combination
of locally defined martingales into a globally defined one, which is
developed in the next section.

\section{Pasting lemma}\label{sec:paste}
The verification that a given stochastic solution satisfies a global
martingale property
relies on ``martingale pasting,'' which is the property
(shown below) that, for a given function $g$,
the combination of two supermartingales of type \eqref{eq:ggmartI}
on two adjacent closed intervals
$I_1$ and $I_2$ gives rise to a supermartingale defined on the union
$I_1\cup I_2$,
provided that, in the case that $X$ has unbounded variation,
$g$ is differentiable at the intersection $I_1\cap I_2$ of $I_1$ and $I_2$.

\begin{Lemma}\label{lem:paste}
Let $(I_i)_{i=1}^n$ be
a finite collection of closed intervals with disjoint interiors
satisfying $\bigcup_{i=1}^n I_i = \mbb R_+$, and let
$g\dvtx\mbb R\to\mbb R$ be a c\`{a}dl\`{a}g function satisfying boundary
condition \eqref{eq:bc}
and growth condition \eqref{eq:gbound}.
Assume in addition that $g$ is differentiable at any $x>0$ with $x\in
\bigcup_{i=1}^n \partial I_i$\footnote{For any set $A$, $\partial A =
\ovl A\setminus A^o$ is the boundary of $A$, where $\ovl A$, $A^o$
denote the closure and interior of~$A$.}
if $X$ has unbounded variation. If
%
\begin{equation}
\label{eq:Ytau}  \mbox{$S^{T_{I_i}}= \bigl\{\mathrm{e}^{-q(t\wedge T_{I_i})}g
(X_{t\wedge
T_{I_i}} ), t\in\mbb R_+ \bigr\}$\qquad are $\mbf F$-supermartingales, }\hspace*{-20pt}
\end{equation}
for $i=1,\ldots, n$, then
%
\begin{equation}
\label{eq:mmmm} S = \bigl\{\mathrm{e}^{-q(t\wedge T_{\mbb R_+})}g (X_{t\wedge
T_{\mbb
R_+}} ), t
\in\mbb R_+ \bigr\}\qquad \mbox{is a UI $\mbf F$-supermartingale}.\hspace*{-20pt}
\end{equation}
\end{Lemma}

The pasting lemma implies in particular that a global super-martingale
property holds for sufficiently regular stochastic supersolutions:

\begin{Cor}\label{prop:global}
Assume that $g$ is a local stochastic supersolution on $I_i$,
$i=1,\ldots, n$,
for some finite collection of closed intervals $(I_i)_{i=1}^n$ with
$\bigcup_{i=1}^n I_i=\mbb R_+$
and $I^o_i\cap I^o_j=\varnothing$ for $i\neq j$. If $X$ has unbounded
variation,
suppose in addition that $g$ is differentiable at any $x>0$ with $x\in
\bigcup_{i=1}^n \partial I_i$.
Then \eqref{eq:mmmm} holds true.
\end{Cor}

\begin{pf*}{Proof of Lemma~\ref{lem:paste}}
In view of the observations that $S$ is $\mbf F$-adapted
and UI [by Lemma~\ref{lem:est}(ii), as $g$ satisfies the linear growth
condition],
it suffices to show that $\E[S_t|\mc F_s]\leq S_s$ for any $s,t\in
\mbb
R_+$ with $s<t$.
For the ease of presentation, only the verification in the case of a
collection of closed intervals
the form $\{[0,a], [a,\infty)\}$ for some $a>0$
is considered, as the general case follows by a similar line of reasoning.

Fix thus $s,t\in\mbb R_+$ arbitrary with $s<t$ and suppose first that
$X$ has bounded variation.
Then $a$ is irregular for $(-\infty,a)$ for $X$, so that the following
collection of stopping times $(T_i)_{i\in\mbb N\cup\{0\}}$ forms a
discrete set:
%
\begin{equation}
\label{eq:Ti}  T_0:=0, \qquad T_{2i}:= T_{[0,a]}
\circ\theta_{T_{2i-1}},\qquad T_{2i-1} = T_{[a,\infty)}\circ
\theta_{T_{2i-2}},\qquad i\in\mbb N,\hspace*{-30pt}
\end{equation}
where $\theta$ denotes the translation operator.
The strong Markov property of $X$ and the tower property of conditional
expectation
imply that,
on the event $\{s\leq T_{i-1}, T_{i-1}<\infty\}$, $i\in\mbb N$,
$\E [S_{t\wedge T_i}-S_{t\wedge T_{i-1}}|\mc F_s  ]$ is
equal to
%
\begin{eqnarray}\label{eq:smi}
&&\E \bigl[\E[S_{t\wedge T_i}-S_{t\wedge T_{i-1}}|\mc F_{T_{i-1}}]|\mc
F_s \bigr]
\nonumber
\\
&&\qquad= \E \bigl[\mbf1_{\{t>T_{i-1}\}}\mathrm{e}^{-qT_{i-1}}\\
&&\hspace*{42pt}{}\times \bigl\{\E
_{X_{t\wedge T_{i-1}}} \bigl[\mathrm{e}^{-q R_v}g(X_{R_v}) |\mc
F_s \bigr] |_{v=T_{i-1}\wedge t} - g(X_{t\wedge T_{i-1}}) \bigr\} \bigr],\nonumber
\end{eqnarray}
with $R_v=(T_i\wedge t) \circ\theta_v$, where the expectation on the
right-hand side is
nonpositive in view of Doob's optional stopping theorem [which holds
in view of the uniform integrability of $S$
and the assumed supermartingale property \eqref{eq:Ytau}].
Since $T_n\to\infty$ $\P$-a.s. as $n\to\infty$ (recalling $\inf
\varnothing
=\infty$
and $X_t\to\infty$ as $t\to\infty$)
and $S$ is UI, it follows $\E [S_t - S_s|\mc F_s  ] = \lim_{n\to
\infty} \E [S^{T_n}_t - S^{T_n}_s|\mc F_s  ]$
is equal to the limit as $n\to\infty$ of
\[
\sum_{j=1}^n \mbf1_{\{T_{j-1}<s\leq T_{j}\}}
\Biggl\{\E \bigl[(S_{T_j\wedge
t}-S_{T_j\wedge s}) |\mc F_s
\bigr] + \sum_{i=j+1}^n \E
\bigl[(S_{t\wedge
T_i} - S_{t\wedge
T_{i-1}}) |\mc F_s \bigr] \Biggr
\},
\]
which is nonpositive.

\begin{figure}[b]

\includegraphics{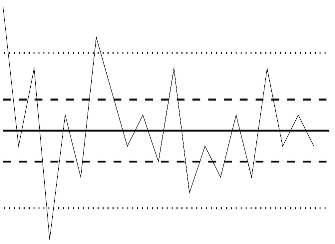}

\caption{The martingale increments commence when $X$ enters the inner
band (dashed) and stop when $X$ leaves the outer band (dotted).}
\label{fig:bandm}
\end{figure}

Suppose next that $X$ has unbounded variation. For any given
$\varepsilon>0$,
denote by
$(T'_i)_{i\in\mbb N\cup\{0\}}$ the sequence of subsequent entrance
times into the sets $[a-\varepsilon,a+\varepsilon]$ and $\mathbb
R\setminus[a-2\varepsilon,a+2\varepsilon]$,
\begin{eqnarray*}
T'_0&:=& 0, \qquad T'_{2i-1}:=
T_{\mbb R\setminus[a-\varepsilon,a+\varepsilon]} \circ \theta_{T'_{2i-2}},
\\
T'_{2i}&:=& T_{[a-2\varepsilon,a+2\varepsilon]}\circ\theta_{T'_{2i-1}},\qquad i\in\mbb N,
\end{eqnarray*}
(see Figure~\ref{fig:bandm}). For any $t\in\mbb R_+$,
decompose $S_t$ as $S_t-S_0=S_t^{(1,\varepsilon)} +
S_t^{(2,\varepsilon)}$ with
\[
S^{(1,\varepsilon)}_t = \sum_{i\ge1}
[S_{t\wedge T'_{2i}} - S_{t\wedge
T'_{2i-1}} ],\qquad S^{(2,\varepsilon)}_t = \sum
_{i\ge1} [S_{t\wedge T'_{2i-1}} - S_{t\wedge T'_{2i-2}} ].
\]
The conditional expectation $\E [S^{(1,\varepsilon
)}_t-S^{(1,\varepsilon)}_s|\mc F_s  ]$,
which concerns increments of $S$ during the periods that $X$ spends
in the band $[a-2\varepsilon,a+2\varepsilon]$, vanishes as
$\varepsilon
\searrow0$, as shown in the
following result:


\begin{Lemma}\label{lem:m1e}
We have $\lim_{n\to\infty}\E [S^{(1,\varepsilon
_n)}_t-S^{(1,\varepsilon_n)}_s|\mc F_s
]\leq0$ a.s.
for some sequence $(\varepsilon_n)_n$ with $\varepsilon_n\searrow0$.
\end{Lemma}

The proof of Lemma~\ref{lem:m1e} is given below.
Since $S^{(2,\varepsilon)}$ is a UI super-martingale for any
$\varepsilon>0$
(which follows by the line of the reasoning given in the first part of
the proof),
we thus have that $\E[S_t|\mc F_s]$ is equal to
\[
\operatorname{lim}_{n\to\infty}\E \bigl[S^{(1,\varepsilon
_n)}_t|\mc
F_s \bigr] + \operatorname{lim}_{n\to\infty}\E
\bigl[S^{(2,\varepsilon
_n)}_t|\mc F_s \bigr] \leq
\operatorname{lim}_{n\to\infty} \bigl(S^{(1,\varepsilon_n)}_s +
S^{(2,\varepsilon)}_s\bigr),
\]
which is equal to $S_s$. As $s$ and $t$ were arbitrary, the proof is complete.
\end{pf*}

Lemma~\ref{lem:m1e} can be established deploying the properties of
Gerber--Shiu functions:

\begin{pf*}{Proof of Lemma~\ref{lem:m1e}}
Let $\varepsilon>0$ be given and, for any $t\ge0$ write
$S_t^{(1,\varepsilon)}=\Sigma
_t^{(1,\varepsilon)} + \Sigma_t^{(2,\varepsilon)}
+ \Sigma_t^{(3,\varepsilon)}$ with
$\Sigma^{(1,\varepsilon)}_t = \sum_{i\ge1} g(X_{t\wedge
T'_{2i}})[\mathrm{e}^{-q(t\wedge T'_{2i})} -\break \mathrm{e}^{-q(t\wedge
T'_{2i-1})}]$,
\[
\Sigma^{(2,\varepsilon)}_t = \sum_{i\ge1}
\mathrm{e}^{-q(t\wedge T'_{2i-1})} \bigl[\E\bigl[g(X_{t\wedge T'_{2i}})|\mc
F_{t\wedge T'_{2i-1}}\bigr] - g(X_{t\wedge
T'_{2i-1}})\bigr]
\]
and $\Sigma^{(3,\varepsilon)}_t = \sum_{i\ge1}\mathrm
{e}^{-q(t\wedge T'_{2i-1})}
[g(X_{t\wedge T'_{2i}}) - \E[g(X_{t\wedge T'_{2i}})|\mc F_{t\wedge
T'_{2i-1}}]]$.
We next estimate these three sums.

In view of growth condition \eqref{eq:gbound}, it follows that there exist
positive real numbers $a$ and $b$ satisfying $\{\forall x\in\mbb R_+,
|g(x)|\leq ax + b\}$,
so that the following estimate holds:
\[
\bigl\llvert \Sigma^{(1,\varepsilon)}_t\bigr\rrvert \leq(a \ovl
X_{t\wedge
\tau
_\pi} + b) \int_0^{t\wedge\tau_\pi}
\mathrm{e}^{-qs} \mbf1_{\{X_s\in
(a-2\varepsilon,a+2\varepsilon)\}}\,\td s,\qquad  t\ge0.
\]
The absolute continuity of the potential measure of $X$ and the
integrability of $\ovl X_t$ for any $t\ge0$
implies that, as $\varepsilon\searrow0$,
the left-hand side tends to zero $\P$-a.s. and in $L^1(\P)$
(by Lebesgue's dominated convergence theorem).

The next step is the observation that the following estimate
holds (as a consequence of the differentiability of $g$ at $a$):

%
\begin{Lemma}\label{lemL}
Let $\eta>0$ and $q\ge0$.
There exists a $\WT C>0$ such that for all $\varepsilon>0$
sufficiently small,
$L(x) = \E_x[\mathrm{e}^{-q T_{a-2\varepsilon,a+2\varepsilon
}}g(X_{T_{a-2\varepsilon,a+2\varepsilon}})] - g(x)$ satisfies
%
\begin{equation}
\label{eq:Lest} \sup_{x\in[a-2\varepsilon,a+2\varepsilon]} L(x)
 \leq\varepsilon \cdot C(
\varepsilon),\qquad C(\varepsilon):= \WT C\bigl[\eta + W^{(q)}(4\varepsilon)
\bigr].
\end{equation}
\end{Lemma}

The proof of Lemma~\ref{lemL} is given below.

The triangle inequality and the strong Markov property
imply that $|\Sigma^{(2,\varepsilon)}_t|$ is bounded by the sum
$\sum_{i\ge1}\mathrm{e}^{-q(t\wedge T'_{2i-1})}|(\WT L_1 + \WT
L_2)(t -
t\wedge T'_{2i-1}, X_{t\wedge T'_{2i-1}})|$
where $\WT L_1(t,x) = \E_x[(g(X_t)-g(x))\mbf1_{\{T>t\}}]$ and
$\WT L_2(t,x) = \E_x[(g(X_{T})-\break g(x))\mbf1_{\{T\leq t\}}]$ with $T=
T_{a-2\varepsilon,a+2\varepsilon}$
may be decomposed as $\WT L_2(t,x) = A_1 -A_2$ with $A_1 = L(x)$,
and
\[
A_2 = \E_x\bigl[\bigl(g(X_T)-g(x)\bigr)
\mbf1_{\{t< T\}}\bigr] = \E_x\bigl[L(X_t)
\mbf1_{\{t <
T\}}\bigr].
\]
To estimate $|\Sigma^{(2,\varepsilon)}_t|$ we split it into two sums.
It is
straightforward to check that
the sum involving the terms $\WT L_1$ is bounded by $\E_x[|g(X_t) -
g(X_\rho)|\mbf1_{\{t<\rho'\}}]$
where $\rho=\sup\{u\leq t\dvtx X_u\in(a-\varepsilon,a+\varepsilon
)\}$ and
$\rho'=\inf\{t>\rho
\dvtx X_t\notin[a-2\varepsilon,a+2\varepsilon]\}$,
which in turn is bounded by $C'\varepsilon$ for some constant $C'$ (as
$g$ is
differentiable in $a$).

Furthermore, it follows from Lemma~\ref{lemL} that
$\WT L_2(t,x)$ is bounded by $2\varepsilon C(\varepsilon)$.
Observe next that the number of terms in the sum $\Sigma
^{(2,\varepsilon
)}$ is
bounded by $1 + D^-_t(\varepsilon) + U^+_t(\varepsilon)$,
where $D^-_t(\varepsilon)$ and $U^+_t(\varepsilon)$ denote the
numbers of
down-crossings
of the band
$(a-2\varepsilon,a-\varepsilon)$ and
upcrossings of $(a+\varepsilon,a+2\varepsilon)$ by $X$ before time $t$.
Thus the expectation of $|\Sigma^{(2,\varepsilon)}_t|$ can be bounded
as follows:
%
\begin{equation}
\label{keyest} \E_x \bigl[\bigl\llvert \Sigma^{(2,\varepsilon)}_t
\bigr\rrvert \bigr] \leq 2\varepsilon \E _x\bigl[1 +
D^-_t(\varepsilon) + U^+(\varepsilon)\bigr] C(\varepsilon) +
C'\varepsilon.
\end{equation}
Since $X$ is a L\'{e}vy process with positive drift, $X$ is a submartingale,
so that the upcrossing lemma implies that
the expected number of upcrossings of the band
$(c,d)=(a+\varepsilon,a+2\varepsilon)$
by time $t$ does not grow faster than $\varepsilon^{-1}$,
\[
\varepsilon\cdot\E_x\bigl[U^+_t(\varepsilon)\bigr]\leq
\E_x\bigl[(X_t-d)^+\bigr]-\E_x
\bigl[(X_0-c)^+\bigr].
\]
Thus, it follows that
$\varepsilon\cdot\E_x[U^+_t(\varepsilon)]$ remains bounded as
$\varepsilon\to0$.
As the number of downcrossings $D^-_t(\varepsilon)$ of the band
$(a-2\varepsilon,a-\varepsilon)$
is bounded by $2 + U^+_t(\varepsilon)$; also
$\varepsilon\cdot\E_x[D^-_t(\varepsilon)]$ remains bounded as
$\varepsilon\to0$.
Since $C(\varepsilon)$ tends to $\eta$
as $\varepsilon\to0$ [as $W^{(q)}(0)=0$ when $X$ has unbounded variation],
it thus follows from \eqref{keyest}
that $\E_x[|\Sigma^{(2,\varepsilon)}_t|]$ tends to $2\eta$ as
$\varepsilon$ tends to zero.
As $\eta$ is arbitrary, we conclude $\lim_{\varepsilon\searrow0}\E
_x[|\Sigma
^{(2,\varepsilon)}_t|]=0$.

Next we turn to the sum $\Sigma^{(3,\varepsilon)}$. We have the decomposition
$\E[\Sigma^{3,\varepsilon}_t - \Sigma^{3,\varepsilon}_s|\mc
F_s]=\sum_{j\ge1}\mbf1_{\{
T_{2j-2}\leq s< T_{2j}\}} B_j$
with $B_j =\mathrm{e}^{-q(t\wedge T_{2j-1}}(E[g(X_{T_{t\wedge
T_{2j}}})|\mc
F_s] - E[g(X_{T_{t\wedge T_{2j}}})|\mc F_{t\wedge T_{2j-1}}]$.
Reasoning as above we find that the sum convergences to 0 in $L^1(\P)$
when $\varepsilon\to0$.
Finally, an application of the Borel--Cantelli lemma (recalling
$S^{(1,\varepsilon)}= \sum_{i=1}^3\Sigma^{(i,\varepsilon)}$)
yields the existence of a sequence
$(\varepsilon_n)$, $\varepsilon_n\to0$, such that $\E
[S^{(1,\varepsilon
_n)}_t-S^{(1,\varepsilon_n)}_s|\mc
F_s]\to0$ a.s. as $n\to\infty$.
\end{pf*}

\begin{pf*}{Proof of Lemma~\ref{lemL}}
By rearranging terms observe that $L(x)$ can be written as
$L(x) = g(a) R_0(x) + g'(a)R_1(x) + R(x) - \tilde w(x)$
with $\tilde w(x):= g(x) - g(a) - g'(a)(x-a)$,
$R(x):= \E_x [\mathrm{e}^{-q T_{a-2\varepsilon,a+2\varepsilon
}}\tilde
w
(X_{T_{a-2\varepsilon,a+2\varepsilon
}}  )  ]$,
$R_0(x):=  \E_x[\mathrm{e}^{-q T_{a-2\varepsilon,a+2\varepsilon}}]-
1$ and
\[
R_1(x):= \E_x\bigl[\mathrm{e}^{-q T_{a-2\varepsilon,a+2\varepsilon
}}(X_{T_{a-2\varepsilon, a+2\varepsilon}}
- a)\bigr] - (x-a).
\]
Next the terms $R_0(x)$, $R_1(x)$ and $R(x)$ are estimated.
Given $\eta>0$, let $\delta>0$ satisfy $|\tilde w(y)/(y-a)|<\eta$,
whenever $|y-a|<\delta$ (such a $\delta$ exists as $g$ is assumed to be
differentiable at $a$).
Then, for any $\varepsilon$ sufficiently small and any $x\in
[a-2\varepsilon, a+2\varepsilon]$, the
bounds $|\tilde w(x)|\leq2\eta\varepsilon$ and $|R(x)| \leq|R_2(x)| +
\eta
|R_3(x)|$ hold, with
%
\begin{eqnarray}\qquad
\label{w2w3}
R_i(x) &=& \E_x\bigl[
\mathrm{e}^{-q T_{a-2\varepsilon,a+2\varepsilon
}}w_i (X_{T_{a-2\varepsilon,a+2\varepsilon}} )\bigr],\qquad i=2,3,
\nonumber
\\[-8pt]
\\[-8pt]
\nonumber
w_2(x) &= &\tilde w(x)\mbf1_{(-\infty,a-\d]}(x),\qquad w_3(x)
= (x-a)\mbf 1_{(a-\d,0]}(x), \qquad x\leq a.
\end{eqnarray}

From expression \eqref{d:wxt0a}, with the replacements $a\to
a-2\varepsilon$,
$b\to a+2\varepsilon$ and $w \to\tilde w_i\in\mc R_0$ for
$i=0,\ldots, 3$
given by $\tilde w_i = {}_{a-2\varepsilon} w_i$ with $w_i\dvtx
(-\infty,a-2\varepsilon]\to\mbb
R$ specified in \eqref{w2w3} and
by $w_0(x):= 1$ and $w_1(x):=x-a+2\varepsilon$, and the fact that
$W^{(q)}$ is
increasing, it is straightforward to verify that, for any $x\in
[a-2\varepsilon,
a+2\varepsilon]$,
%
\begin{equation}
\label{Ri} \bigl|R_i(x)\bigr| \leq2 \max_{z\in[0,4\varepsilon]}\bigl|F_{\tilde w_i}(z)
- \tilde w_i(0) - \tilde w_{i,-}'(0)z\bigr|,\qquad
i=0,1,2.
\end{equation}

Since the functions $J_{\tilde w_i}$, $i=0,1,2$, given in \eqref{J}
with $w\to\tilde w_i$, are
bounded, by $J_\infty$ say, and $W^{(q)}$ is increasing,
it follows from the form \eqref{eq:Fwr} of $F_w$ that
$|F_{\tilde w_i}(z) - \tilde w_i(0) - \tilde w_{i,-}'(0)z|$, $i=0,1,2$,
$z\in[0,4\varepsilon]$, is bounded by
%
\begin{equation}
\label{Fwi} J_\infty\int_0^z
W^{(q)}(z-y)\,\td y \leq J_\infty\cdot4\varepsilon\cdot
W^{(q)}(4\varepsilon).
\end{equation}
Combining \eqref{Ri} and \eqref{Fwi} yields that the functions
$R_i(x)$, $i=0,1,2$, are each bounded by
$J_\infty\cdot8\varepsilon W^{(q)}(4\varepsilon)$ for any $x\in
[a-2\varepsilon,a+2\varepsilon]$.
Similarly, it follows from the facts $F_{\tilde w_3}(0) = \tilde
w_3(0)=0$ and $F_{\tilde w_3}'(0+) = \tilde w_{3,-}'(0)=1$
(Theorem~\ref
{thm:sd})
that, for all $\varepsilon$ sufficiently small, $|R_3(x)|\leq C_1
\varepsilon$, for all
$x$ in the interval $[a-2\varepsilon,a+2\varepsilon]$
for some constant $C_1>0$. Combining the estimates for $\tilde w(x)$
and $R_0(x),\ldots, R_3(x)$
with the form of $L(x)$ completes the proof.
\end{pf*}

%
\section{Optimality conditions for single dividend-band
strategies}\label{ssec:singo}
A necessary and sufficient condition for the optimality of the single band
policy $\pi_{\unl b^*}$ at levels $\unl b^*:= b^*_1 = (b^*_-, b^*_+)$
defined in \eqref{eq:astar}--\eqref{eq:bK0}
can be expressed
in terms of the function $G^*\dvtx(b^*_-,\infty)\to\mbb R$ given by
%
\begin{equation}\qquad
\label{eq:Gstar} G^*(y) = G\bigl(b^*_-, y\bigr) = \cases{\displaystyle \frac{y-b^*_- -
K - (F(y) - F(b^*_-))}{W^{(q)}(y)-W^{(q)}(b^*_-)}, &\quad
$\mbox{if $K>0$},$\vspace*{2pt}
\cr
\displaystyle G^\#(x) = \frac{1 - F'(x)}{W^{(q)\prime}(x)}, &\quad $\mbox{if
$K=0$}$.}
\end{equation}
This condition can be expressed in terms of the function $Z^{(q,v)}$
that was defined in Definition~\ref{def:Zqv}.

%
\begin{Thm}\label{thm:cm2}
\textup{(i)} The single-band policy $\pi_{\unl b^*}$ at level $\unl b^*=b_1^*$
is optimal for the stochastic control problem \eqref{optdiv} if and
only if
%
\begin{equation}
\label{HJBs} {} _{b^*_+} \bigl(\Gamma_\infty^{\ovl w}
v_{b^*} - q v_{b^*}\bigr) (x) \leq0 \qquad\mbox{for all $x>b^*_+$
and with $\ovl w=v_{b^*}$,}
\end{equation}
where the operator $_{b^*_+} \Gamma^{\ovl w}_{\infty}$ is defined in
\eqref{eq:gamma2}, or equivalently, if and only if
$\Xi^*\dvtx(\Phi(q),\infty)\to\mbb R$ is completely monotone, where
%
\begin{equation}
\label{eq:Xi} \Xi^*(\theta) = - \frac{\mathrm{e}^{\theta b^*_+}}{\theta} \int_{(b^*_+,\infty)}
\mathrm{e}^{-\theta z} Z^{(q,\theta)\prime}(z) G^*(\td z), \qquad\theta>\Phi(q).
\end{equation}
\textup{(ii)} In particular, if $G^*$ is nonincreasing on
$(b^*_+,\infty)$, then the strategy $\pi_{b^*}$ is optimal.
\end{Thm}

Theorem~\ref{thm:cm2}(ii) yields a useful simple sufficient
optimality condition:

\begin{Cor} \label{Cor:one}\textup{(i)} The unimodality of the function
$G^*$ implies the optimality of single dividend-band policies.

\textup{(ii)} In particular, in the case $K=0$ and if $G^\#$ is monotone
decreasing, then the ``lump-sum'' strategy $\pi_0$ is optimal.
\end{Cor}

\begin{Rem}
In the absence of transaction costs, the function $\Xi^*$
in \eqref{eq:Xi} can be equivalently expressed as
\begin{eqnarray*}
 \Xi^*(\theta)& =& G^\#\bigl(b^*_+\bigr) L_0(\theta) +
\frac{(\psi(\theta)-q)}{\theta^{2}}\E\bigl[F'\bigl(b^*_+ + \mbf e_\theta
\bigr) - F'\bigl(b^*_+\bigr)\bigr],
\\
 L_0(\theta)&:=& \frac{\psi(\theta) - q}{\theta^2} \E\bigl[W^{(q)\prime}
\bigl(b_+^*+\mbf e_\theta\bigr) - W^{(q)\prime}\bigl(b_+^*\bigr)
\bigr],
\end{eqnarray*}
where $\mbf e_\theta$ denotes an independent exponential random
variable with mean $\theta^{-1}$.
In particular, if the penalty is zero and there are no transaction
cost $(w=K=0)$, the necessary and sufficient optimality condition
simplifies to the complete monotonicity of $L_0(\theta)$ on the
interval $(\Phi(q),\infty)$. This observation appears new even in
this particular case.
\end{Rem}

\begin{Rem}[({Lump-sum strategy})]
In the absence of transaction cost \mbox{($K=0$)}, the ``lump-sum''
strategy $\pi_0$ is to ``pay out all the reserves to the
beneficiaries and subsequently pay all the premiums as dividends,
until the moment of ruin.'' Note that $\pi_0$ is a
single dividend-band strategy at level $0$. In the case that $X$
is given by the Cram\'{e}r--Lundberg model, the first jump (claim)
arrives after an independent exponential time $\mbf e_\lambda$ with
finite mean
$\lambda^{-1}$, so that the value $v_0$ is equal to
\begin{eqnarray*}
v_0(x) &=& \E_x \biggl[x + p \int_0^{\mbf e_\lambda}
\mathrm{e}^{ -q t} \,\td t + \mathrm{e}^{ -q
\mbf e_\lambda} w(\Delta
X_{\mbf e_\lambda}) \biggr]
\\
&=&\E_x \biggl[x+ \frac{p}{q}\bigl(1-\mathrm{e}^{-q{\mbf e_\lambda}}
\bigr) + \mathrm{e}^{ -q {\mbf e_\lambda}}\bigl(w(\Delta X_{\mbf e_\lambda})-w(0)\bigr)+
w(0)\mathrm{e}^{
-q {\mbf e_\lambda}} \biggr],
\end{eqnarray*}
which is equal to $x + \frac{p + w_\nu(0)+\lambda w(0)}{\lambda+q}$,
where $\Delta X_{\mbf e_\lambda} =
X({\mbf e_\lambda})-X({\mbf e_\lambda}-)$, and $w_{\nu}\dvtx\mbb
R_+\setminus\{0\}\to\mathbb R$ is defined
in Proposition~\ref{prop:GS}. If $X_0$ is zero and $X$ has infinite
activity or nonzero Gaussian component, ruin
occurs immediately if strategy $\pi_0$ is followed ($\tau^{\pi_0}=0$,
$\P_0$-a.s.)
and $v_0(x)=x +w(0)$.

Hence, the value of the lump-sum strategy is equal to
$v_0(x) = (x + \gamma_w)\times\mbf1_{[0,\infty)}(x) + w(x) \mbf
1_{(-\infty,0)}(x)$ with $\gamma_w =
v_0(0)$ given by
\[
\cases{ \displaystyle\frac{1}{q + \ovl\nu} \bigl[p + w_\nu(0) +
\ovl\nu w(0) \bigr],
& \quad $\mbox{if $\ovl\nu:=\nu(\mbb R_+)<\infty$ and $\sigma =0$}$,\vspace *{2pt}
\cr
w(0), &\quad  $\mbox{if $\ovl\nu=\infty$ or $\sigma>0$}$.}
\]
If $G^\#$ is monotone decreasing, it attains its maximum over $\mbb R_+$
at zero,
and the function $\Xi$ is completely monotone, so that
$\pi_0$ is optimal [Theorem~\ref{thm:cm2}(ii)].
\end{Rem}

\begin{Rem}
In the following result (proved in  Appendix \ref{ssec:12band}) explicit
sufficient conditions
are given in terms of the penalty $w$ and the L\'{e}vy density $\nu$
for optimality of a single barrier strategy at a positive level:

\begin{Cor}\label{thm:expl}
In the case $\{K=0$ and $b_1^*>0\}$, if $\nu$ admits a
convex density $\nu'$ and the penalty $w$ is \emph{severe} [i.e.,
$w(0)\leq\gamma_w$ and $w(x+y)-w(y)\leq x$ for all $x,y\in\mbb R_-$],
then the strategy $\pi_{b_1^*}$ is optimal.
\end{Cor}

Note that a penalty $w$ is severe if (i) the penalty at 0 is at least
the value of the lump-sum strategy at 0
and (ii) the slope of the penalty is at least one.
\end{Rem}

\begin{pf*}{Proof of Theorem~\ref{thm:cm2}, part (i)}
The equivalence of the conditions \eqref{HJBs} and \eqref{eq:Xi}
directly follows due to Lemma~\ref{prop:key}(iii).

\textit{Proof of sufficiency of} \eqref{HJBs}:
It suffices to show that $v_{b^*}$
is a stochastic supersolution, as then the local verification theorem
(Theorem~\ref{cor:repg}) implies
that $v_{b^*}$ is equal to the value-function $v_*$. The supersolution
property of
$v_{b^*}$ follows by combining the pasting lemma (Lemma~\ref
{lem:paste}) with the following facts:
\begin{longlist}[(a)]
\item[(a)]
$\exp\{-q(t\wedge T^-_{b^*_+})\}v_{b^*} (X(t\wedge
T^-_{b^*_+})  )$
is an $\mbf F$-supermartingale
[by \eqref{HJBs} and Lemma~\ref{lem:smp}(ii)],
\item[(b)] $\exp\{-q(t\wedge T_{0,b^*_+})\}v_{b^*} (X(t\wedge
T_{0,b^*_+})  )$
is an $\mbf F$-martingale [by the form of $v_{b^*}$
in \eqref{eq:v*b} and the martingale properties of $W^{(q)}$ and $F_w$
in Proposition~\ref{prop:mart}] and
\item[(c)] if $X$ has unbounded variation, $v_{b^*}$ is
differentiable at
$b^*_+$ [in view of the form of $v_{b^*}$
in \eqref{eq:v*b}].
\end{longlist}

\textit{Proof of necessity of} \eqref{HJBs}: Suppose that
the condition in \eqref{HJBs} is not satisfied. Since
$x\mapsto({} _{b^*_+} \Gamma^{\ovl w}_\infty v_{b^*} - q v_{b^*})(x)$
is right-continuous
at any $x$ with $x> b^*_+$, it follows that there exists an open interval
$(\a,\b)$ contained in $(b^*_+,\infty)$ with $({} _{b^*_+} \Gamma
^{\ovl
w}_\infty v_{b^*} - qv_{b^*})(x) > 0$ for $x\in(a,b)$. Define a
strategy $\tilde\pi$ as
follows: whenever $U_t$ does not take a value in the interval
$(\a,\b)$, operate according to $\pi_{b^*}$, and while the reserve
process $U_t$ takes a value in the interval $(\a,\b)$, do not pay
any dividends. Then $S_t:=\mathrm{e}^{-q(t\wedge
T_{\a,\b})}(v_{\tilde\pi}(X_{t\wedge
T_{\a,\b}})-v_{b^*}(X_{t\wedge T_{\a,\b}}))$ is an {$\mbf
F$}-supermartingale,
and the following holds true [cf. \eqref{eq:mmart}] for any
$x\in(\a,\beta)$:
\[
v_{\tilde\pi}(x) - v_{b^*}(x) \ge \E_x[S_t-S_0]
= \E_x \biggl[\int_0^{t\wedge T_{\a,\b}}
\mathrm{e}^{-qs}{} \bigl( {}_{b^*_+} \Gamma^{\ovl w}_{\infty}
v_{b^*} - qv_{b^*}\bigr) (X_s)\,\td s \biggr] >
0.
\]
Hence it follows that $\pi_{b_*}$ is not an optimal policy, and the
proof is complete.
\end{pf*}

\begin{pf*}{Proof of Theorem~\ref{thm:cm2}, part (ii)}
The statement follows by combining part (i) with the next result.
\end{pf*}

\begin{Lemma}\label{lem:xil}
If $x\mapsto G^*(x)$ is nonincreasing on $(b^*_+,\infty)$, then
$\Xi(\theta)$ is completely monotone on $(\Phi(q),\infty)$.
\end{Lemma}

\begin{pf}
If the function $G^*$ is nonincreasing, then the function $\Xi$ is completely
monotone in view of the form
of $\Xi$ given in \eqref{eq:Xi}, the complete monotonicity
of $\theta^{-1}\mathrm{e}^{\theta(b-x)}Z^{(q,\theta)\prime}(x)$ [cf.
Remark~\ref
{rmZ}(ii)]
and the following facts:
\begin{longlist}[(iii)]
\item[(i)] A function $f\dvtx(c,\infty)\to\mbb R_+$, $c>0$, is completely monotone
if and only if $f$ is the Laplace transform of a measure supported on
$[0,\infty)$.

\item[(ii)] If $f(\theta)$ is the
Laplace transform of the measure $\mu$ supported on $[0,\infty)$, then
for any
$c>0$, $\mathrm{e}^{-\theta c}f(\theta)$ is the Laplace
transform of the translated measure $y\mapsto\mbf1_{\{y\ge c\}}\mu
(\td(y-c))$.

\item[(iii)] The Laplace transform of the measure
$n(\td y) = \int_{[b,\infty)}\mu_x(\td y)m(\td x)$ supported on
$[0,\infty)$
is given by $\mc L n(\theta) = \int_{[b,\infty)}\mc L\mu_x(\theta
)m(\td x)$
where $(\mu_x, x>b)$, $b\in\mbb R$, is a collection measures supported
on $[0,\infty)$.\quad\qed
\end{longlist}
\noqed\end{pf}
%

\section{\texorpdfstring{Optimality conditions for solutions to the mixed optimal stopping/\break control~problem}{Optimality conditions for solutions to the mixed optimal
stopping/control~problem}}\label{sec:doublop}
The Hamilton--Jacobi--Bellman equation associated to the
stochastic control problem in \eqref{eq:optstop}
differs from \eqref{eq:HJB} by the inclusion of
the additional requirement that the value-function
should be larger than the function $f$ [reflecting the fact that \eqref
{eq:optstop} is a
mixed optimal stopping/control problem];
hence, the HJB equation corresponding to \eqref{eq:optstop} is given by
%
\begin{eqnarray}
\label{HJBos} &&\max\bigl\{\mathcal L g(x) - q g(x), f(x) - g(x), 1-\mathtt
d_g(x)\bigr\} = 0,\qquad  x>0,
\\
&& \cases{ g(x) = f(x), & \quad $\mbox{for all $x<0$}$, \vspace*{2pt}
\cr
g(0) = f(0),
&\quad $\mbox{in the case $\bigl\{\sigma^2>0$ or $ \nu_{0,1} =
\infty\bigr\}$}$,} \label{HJBosbc}
\end{eqnarray}
where $\mathtt d_g(x)$ is defined in \eqref{qgamma}.
Stochastic supersolutions
$g$ of the HJB equation in~\eqref{HJBos} and \eqref{HJBosbc}
are defined as
in Definition~\ref{def:sss}, with the additional
requirement $g\ge f$. By a line of reasoning similar to that
used in the proof of Theorem~\ref{cor:repg},
it follows that a local verification result
for the stochastic control problem \eqref{eq:optstop} holds true:

\begin{Cor}\label{cor:optstop}
Let $g$ be a stochastic supersolution of the HJB equation
in \eqref{HJBos} and \eqref{HJBosbc}. If there exist
$c,a,b_-,b_+$ satisfying
$0\leq c\leq a\leq b_-\leq b_+$ and $g(x) = V^f_{a,b_-,b_+}(x)$ \{$g(x)
= V^{f,\varnothing}_{a,b_+}(x)$\}
for any $x\in[c,b_+]$,
then it holds $V_*^f(x) = V^f_{a,b_-,b_+}(x)$ for all $x\in[c,b_+]$
\{$V_*^f(x) = V^{f,\varnothing}_{a,b_+}(x)$ for all $x\in[c,b_+]$\},
respectively.
\end{Cor}

Given this verification result the proof of Theorem~\ref{thmos} can be
completed.
A key step in the proof is the following property of the function $f$:

\begin{Lemma}\label{lem:a2}
Suppose that $f$ satisfies the conditions in \eqref{eq:as1}--\eqref{eq:as3},
and denote $\bar w = f|_{\mbb R_-}$. It holds $_0\Gamma^{\bar
w}_\infty
f(u) \leq0$ for all $u\in (0,\alpha(K)  )$ with $\alpha
(0):=
\alpha
^*_f$ and $\alpha(K):=\alpha^*_{f,\varnothing}$ for $K>0$.
\end{Lemma}

\begin{pf*}{Proof of Theorem~\ref{thmos}}
(i) Since $V^f_{\a^*_f,\b^*_f}$ is the value-function of the strategy
$ (\tau^{\pi_{\b^*}}_{\a^*},\pi_{\b^*}  )$,
Corollary~\ref{cor:optstop} implies
that, to prove the assertion, it suffices to show that
$V^f_{\a^*_f,\b^*_f}$ is a supersolution of the
HJB equation in \eqref{HJBos} and \eqref{HJBosbc}. Next the various
conditions are verified.

Analogously to the proof of Theorem~\ref{thm:cm}, it follows from
the definition of $\b^*_f$ and the form of the function
$V=V^f_{\a^*_f,\beta^*_f}$ given in
Proposition~\ref{prop:lp2} that the following inequality holds:
%
\begin{equation}
\label{eq:Vf} V(x) - V(y) \ge x-y-K
\end{equation}
for all $x,y\ge0$ satisfying $x\ge y\ge\a^*_f$.
In view of the fact $V'(x)=f'(x)=1$ for $x\in(0,\a^*_f)$, it follows
that the inequality in \eqref{eq:Vf}
is in fact valid for all $x$ and $y$ satisfying $x\ge y\ge0$.

To see that the $V$ dominates the function $f$,
%
\begin{equation}
\label{eq:Vff} V(x) \ge f(x),\qquad x\ge0,
\end{equation}
note first that it holds $V(0)=f(0)$ (a direct consequence of the form
of $V$ in Proposition~\ref{prop:lp2} and $\a^*_f>0$ by Lemma~\ref{lem:a}).
In the case $K=0$, \eqref{eq:Vff} is hence a special case of
\eqref{eq:Vf} (with $y=0$).
In the case $\{K>0$ and $\alpha^*_{f,\varnothing} \ge\alpha^*_{f}\}$,
the definitions of $\alpha^*_{f,\varnothing}$, $\beta
^*_{f,\varnothing}$
and $G^{(a)}_{f,\varnothing}$,
the positivity of $W^{(q)}(x)$ imply
\begin{eqnarray*}
&& G^{(a)}_{f,\varnothing}(b)\leq0\qquad\mbox{for all $a\in\bigl[0,\alpha
^*_{f,\varnothing}\bigr]$ and $b\in\bigl[0,\beta^*_{f,\varnothing}\bigr]$}
\\
&&\qquad \Longleftrightarrow\quad F^{(a)}(x-a)\ge f(x) \qquad\mbox{for all $x\in
\bigl[0, \beta^*_{f,\varnothing}(a) \bigr]$ and $a\in \bigl[0, \alpha
^*_{f,\varnothing} \bigr]$},
\end{eqnarray*}
which yields the inequality in \eqref{eq:Vff}, in view of the facts
$V(x)=F^{(a)}(x-a)$ for all $x\leq b:= \beta^*_{f,+}$
[by Proposition~\ref{prop:lp2} and Lemma~\ref{lem:a}(i) and the fact
$\b^*_{f,+} \leq\b^*_{f,\varnothing}$ which holds by Lemma~\ref
{lem:a}(ii)], and
$V|_{[b,\infty)}$ is affine (Proposition~\ref{prop:lp2}).

In view of the observations
%
\begin{eqnarray}
\label{eq:supermartingalef} &&\mathrm{e}^{-q (t\wedge T_{0,\a^*_f} )}f (X_{t\wedge
T_{0,\a
^*_f}} ) \qquad\mbox{is an $
\mbf F$-supermartingale, and}
\\
&&\mathrm{e}^{-q  (t\wedge T^-_{\a^*_f} )}F^{(\a^*_f)} \bigl(X_{t\wedge
T^-_{\a^*_f}}-
\a^*_f \bigr)\qquad\mbox{is an $\mbf F$-martingale,} \label{eq:martingaleF}
\end{eqnarray}
and the differentiability of
$F^{(\a^*_f)}(x)$ at $x=0$ if $X$ has unbounded variation
[$F^{(\a^*_f)\prime}(0)=f'_-(\a^*_f)$, by Lemma~\ref{lem:}], it follows
from the pasting lemma (Lemma~\ref{lem:paste})
%
\begin{equation}
\label{eq:smopts} \mathrm{e}^{-q  (t\wedge T^-_{0} )}F^{(\a^*_f)} \bigl(X_{t\wedge
T^-_{0}}-
\a^*_f \bigr)\qquad\mbox{is an $\mbf F$-supermartingale.}
\end{equation}
Here, the supermartingale property in \eqref{eq:supermartingalef}
follows from
Lemma~\ref{lem:a}(i), by a line of reasoning that is similar to the
one used
in the proof of Lemma~\ref{lem:smp}, while the martingale property in
\eqref{eq:martingaleF} follows from Proposition~\ref{prop:GS}.

The supermartingale property in \eqref{eq:smopts} and
the inequalities in \eqref{eq:Vf} and \eqref{eq:Vff}
imply that $F^{(\a^*_f)}(x-\a^*_f)$ is a stochastic supersolution
for the stochastic control problem in \eqref{eq:optstop}, which
completes the proof of (i).

(ii) The line of reasoning is analogous to the one in part (i) (see
Remark~\ref{rem:aa})
and is therefore omitted.
\end{pf*}

\subsection{Optimality conditions for two-band policies}
When a single band strategy is not globally optimal for the stochastic control
problem in \eqref{optdiv}, it is not optimal to pay out a lump-sum
dividend at all levels
above $b^*_+$ but is instead optimal to postpone paying dividends when
the reserves process is in
a certain subset of $(b^+_*,\infty)$. This section is concerned with the
necessary and sufficient conditions for optimality of a policy with
only one additional band.
Consider the candidate optimal two-band strategy $\pi_{\unl a^*,\unl b^*}$
at the levels $\unl a^*=(0,a^*_2)$ and
$\unl b^*=(b_1^*, b_2^*)$ where the levels $b_1^*=(b^*_-, b^*_+)$
associated to the first band have been defined
in \eqref{eq:astar}--\eqref{eq:bK0},
and where the levels associated to the second band are given by
\[
\bigl\{a_2^*,b_2^* \bigr\} = b_{1,+}^* +
\cases{ \bigl\{\a^*_{w^*}, \bigl(\beta^*_{w^*,-},
\beta^*_{w^*,+} \bigr) \bigr\}, &\quad $\mbox{if $K=0$}$\vspace*{2pt}
\cr
&\quad $
\mbox{or $\bigl\{K>0$ and $\a^*_{w^*,\varnothing}\ge\a^*_{w^*}\bigr\}
$}$;\vspace *{2pt}
\cr
\bigl\{\a^*_{v_{b^*_1},\varnothing}, \bigl(b^*_-,\beta
^*_{w^*,\varnothing
} \bigr) \bigr\}, & \quad $\mbox{if $\bigl\{K>0$ and $
\a^*_{w^*,\varnothing} < \a^*_{w^*}\bigr\}$}$, }
\]
where $w^*:={} _{b_{1,+}^*} v_{b^*_1}$
and the levels $\a^*_{w^*}$, $\a^*_{w^*,\varnothing}$, $\beta^*_{w^*,-},
\beta^*_{w^*,+}$
and $\beta^*_{w^*,\varnothing}$ are defined in \eqref
{eq:astar2}--\eqref
{eq:bstare}.

Necessary and sufficient conditions for the two-band policy $\pi_{\unl
a^*,\unl b^*}$ to be (globally) optimal
are expressed in terms of the functions $\Xi^*$ defined in \eqref
{eq:Xi} and the function
\[
\Xi^{**} = \cases{ \Xi_{a^*_2,b^*_2}\bigl(w^*\bigr),&\quad $\mbox{if
$K=0$ or $\bigl\{K>0$ and $\a ^*_{w^*,\varnothing}\ge\a^*_{w^*}\bigr
\}$,}$\vspace*{2pt}
\cr
\Xi^{\varnothing}_{a^*_2,b^*_2}\bigl(w^*\bigr), &\quad $
\mbox{if $\bigl\{K>0$ and $\a^*_{w^*,\varnothing} < \a^*_{w^*}\bigr\}$.}$
}
\]
Here for any $a$, $b_-$ and $b_+$
with $a\leq b_-\leq b_+$ and $f\in\mc R_0$ the functions
$\Xi_{a,b_-, b_+}(f)$ and $\Xi^{\varnothing}_{a,b_+}(f)$ are given by
\begin{eqnarray*}
 \Xi_{a,b_-, b_+}(f)\dvtx\theta&\mapsto& - \frac{\mathrm{e}^{\theta b_+}}{\theta} \int
_{(b_+,\infty)} \mathrm{e}^{-\theta z} Z^{(q,\theta)\prime}(z)
G^{(a)}_{f,b_-}(\td z),
\\
\Xi^{\varnothing}_{a, b}(f)\dvtx\theta&\mapsto& - \frac{\mathrm{e}^{\theta b}}{\theta}
\int_{(b,\infty)} \mathrm{e}^{-\theta z} Z^{(q,\theta)\prime}(z)
G^{(a)}_{f,\varnothing}(\td z),
\end{eqnarray*}
where, for any $z\ge b_-$, $G^{(a)}_{f,b_-}(z):= G^{(a)}_f(b_-,z)$,
and the functions $G^{(a)}_{f,\varnothing}$ and $G^{(a)}_f$ have
been defined in \eqref{eq:Gfem} and \eqref{eq:Gfabb}.

Before stating the optimality condition for this two-band policy, we
first state
a condition for (global) optimality of the policies $ (\tau^{\pi
_{\b
_f^*}}_{\a^*_f},\pi_{\b^*_f} )$
and\break $ (T_{\a^*_{f,\varnothing},\beta^*_{f,\varnothing}},\pi
^\varnothing  )$
in the auxiliary stochastic control problem in \eqref{eq:optstop}.

\begin{Thm}\label{thm:cm3}
Suppose that $f$ satisfies the conditions in \eqref{eq:as1}--\eqref{eq:as3}.
\begin{longlist}[(ii)]
\item[(i)] Suppose that it holds either $K=0$ or $\{K>0$ and $\a
^*_{f,\varnothing
}\ge\a^*_{f}\}$.
Then the strategy $(\tau^{\pi_{\b_f^*}}_{\a^*_f},\pi_{\b^*_f})$
is optimal for the stochastic
optimal control problem in \eqref{eq:optstop} if and only if the
function $\Xi_{\a_f^*,\b^*_{f,-}, \b^*_{f,+}}(f)$ is completely monotone.

\item[(ii)] Suppose that it holds $\{K>0$ and
$\a^*_{f,\varnothing} < \a^*_{f}\}$.
Then the strategy $ (T_{\a^*_{f,\varnothing},\beta
^*_{f,\varnothing
}},\pi
^\varnothing  )$
is optimal for the stochastic
optimal control problem in \eqref{eq:optstop} if and only if the
function $\Xi_{\a_{f,\varnothing}^*,\b^*_{f,\varnothing}}(f)$ is completely
monotone.
\end{longlist}
\end{Thm}

The proof of Theorem~\ref{thm:cm3} is omitted as it is analogous to the
proof of Theorem~\ref{thm:cm2}(i).

\begin{Rem}\label{rem:gencf}
As in the proof of Lemma~\ref{prop:key}, it can be shown that
the complete monotonicity of the function
$\Xi_{\a_f^*,\b^*_{f,-}, \b^*_{f,+}}(f)$
is equivalent to the condition
%
\begin{equation}
\label{eq:gencf} _0 \Gamma^w_\infty
V^f_*(x) - qV^f_*(x) \leq0 \qquad\mbox{for all $x>
\b^*_{f,+}$}.
\end{equation}
Similarly, it follows that the complete monotonicity of
$\Xi_{\a_{f,\varnothing}^*,\b^*_{f,\varnothing}}(f)$ is equivalent to
\eqref
{eq:gencf}
with $\b^*_{f,+}$ replaced by $\b^*_{f,\varnothing}$.
\end{Rem}

The relationship between the stochastic control problems in
\eqref{optdiv}
and \eqref{eq:optstop} (cf. the discussion at the beginning of
Section~\ref{aux2})
immediately yields necessary and sufficient optimality conditions
for the two-band strategy $\pi_{\unl a^*,\unl b^*}$:

%
\begin{Cor}\label{twoband}\textup{(i)} The two-band strategy $\pi_{\unl
a^*,\unl b^*}$
at finite levels $\unl a = (0,a_2^*)$ and $\unl b = (b_1^*, b_2^*)$
is optimal for \eqref{optdiv}
if and only if $\Xi^*$ is \emph{not} completely monotone and
$\Xi^{**}$ is completely monotone.

\textup{(ii)} If $\Xi^*$ is not completely monotone then the levels $a_2^*$ and
$b_{2,+}^*$
are finite, and it is optimal to adopt the two-band
strategy $\pi_{\unl a^*,\unl b^*}$ while the reserves are below $b^*_{2,+}$,
and it holds (with $F^{(a^*_{2,+})}_* = F_{_{a^*_{2,+}} v_*}$)
%
\begin{eqnarray}
v_*(x) = \cases{\displaystyle W^{(q)}(x) \frac{1 - F_w'(b^*_{1,+})}
{W^{(q)\prime}(b^*_{1,+})} + F_w(x), &\quad
$x\in\bigl[0,b^*_{1,+}\bigr]$,\vspace*{2pt}
\cr
x - b_{1,+}^*
+ v_*\bigl(b^*_{1,+}\bigr), & \quad $x\in\bigl(b^*_{1,+},a^*_{2,+}
\bigr)$,\vspace *{2pt}
\cr
F^{(a^*_{2,+})}_*\bigl(x-a^*_{2,+}\bigr), &\quad
$x\in\bigl[a^*_{2,+},b^*_{2,+}\bigr]$. }
\end{eqnarray}
\end{Cor}

\section{Multi dividend-band policies: The recursion for the
dividend-band levels}\label{sec:multib}

A flexible class of dividend strategies are the so-called multi
dividend-band strategies, which generalize the single and two-band
strategies, and are specified as follows:

\begin{Def}\label{def:multib}
The \emph{multi dividend-band strategy} $\pi_{\unl a,\unl b}$,
associated to sequences $\unl a=(a_n)_{n}$, $\unl b^- =
(b^-_n)_n$, $\unl b^+ = (b^+_n)_n$ with $a_n, b^-_n, b_n^+\in
[0,\infty]$ satisfying the intertwining conditions
\[
a_1=0 \leq b^+_1 < a_2 \leq
b^+_2< \cdots <a_{n} \leq b^+_n < \cdots,\qquad
b_n^-\leq b_n^+,
\]
is described as follows:
\begin{longlist}[(iii)]
\item[(i)] when $U^{\underline{a},\underline{b}}:=
U^{\pi_{\underline{a},\underline{b}}} = y \in(b_n^+, a_{n+1})$,
make a lump-sum payment $y-b_n^-$;

\item[(ii)] when $U^{\underline{a},\underline{b}}=b_n^+$, make a
lump-sum payment $b_n^+-b_n^-$, if $K>0$, and pay the minimal
amount to keep $U^{\underline{a},\underline{b}}$ below
$b_n^-=b_n^+$ if $K=0$;

\item[(iii)] while $U^{\underline{a},\underline{b}} \in[a_n,b^+_{n})$,
do not pay any dividends.
\end{longlist}
The strategy $\pi^{\unl a,\unl b}$ is called an $N$-dividend-bands
strategy if
$b_N^+<\infty= a_{N+1}$.
\end{Def}

\begin{figure}[b]

\includegraphics{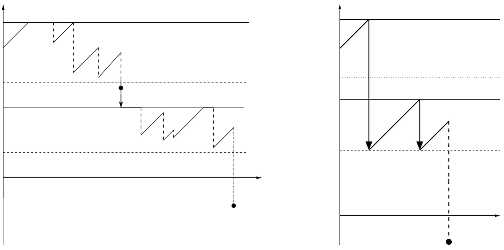}

\caption{Illustrated in the figure on the left is a
path of
the risk process $U^\pi$ in the absence of transaction
cost ($K=0$) for a three-band
strategy with the lowest level $b_1^+$ equal to zero.
The figure on the right pictures a path of the risk process
$U^\pi$ in
the case $K>0$, and $\pi$ is a two-band strategy with $b^-_2=b^-_1$.
The vertical dashed stretches represent the claims, while
lump-sum dividend payments are indicated by arrows.
At the moment $\tau$ of ruin a penalty payment $w(U_\tau)$
is required that is a function of the shortfall $U_\tau$.}\label
{fig:band}
\end{figure}

%

A multi dividend-band strategy $\pi_{\unl a,\unl b}$ consists of paying
out ``the minimal amount
to keep $U^{\underline{a},\underline{b}}_t$ below the boundary
$b(t)$,'' where
\[
b(t):=b^+_{\rho(t)} \qquad\mbox{with $\rho(t) = \min\bigl\{i \in\mbb N\dvtx
U^{\unl a, \unl b}_t< a_i\bigr\}$}.
\]
In this case, while the boundary $b(t)$ is constant,
$U_t^{\underline{a},\underline{b}}$ is equal to
the process $X$ reflected at the level $b(t)$ and the
corresponding cumulative dividend payments $D_t^{\underline
{a},\underline{b}}$
are equal to a local time of $U_t^{\unl a, \unl b}$ at $b(t)$. In the
case of a positive fixed
transaction cost $K$ the ``reflection boundaries'' $b_n^+$
widen to strips $[b_n^-, b_n^+]$, and the ``local time'' type
payments are replaced by lump-sum payments $b_n^+-b_n^-$ where
$b_n^-$ may lie below $a_{n-1}$; see Figure~\ref{fig:band}.

\subsection{Construction of the candidate solution of the stochastic
control problem}\label{sec:multialg}
The dynamic programming equation satisfied by the optimal value
function is recursive in nature,
due to the presence of only negative jumps in both the uncontrolled
reserves process $X$ and the controlled reserves process $U^\pi$ for
any admissible policy $\pi$. In conjunction with the form of the
optimal strategy
of the mixed optimal stopping/stochastic control problem \eqref
{eq:ooptst}, this suggests that the candidate optimal
policy for the stochastic control problem takes in general the form of
a multi-dividend-band strategy $\pi_{\unl a^*,\unl b^*}$ at certain
levels $\unl a^*$, $\unl b^*$. By
repeatedly solving mixed-optimal stopping/stochastic control problems
of the form \eqref{eq:optstop}
with suitably updated reward functions $f$, these levels $\unl a^*$,
$\unl b^*$
can be identified, as summarized in
the following recursive procedure:\vspace*{6pt}

%
\hspace*{-12pt}{\fontsize{10}{12}\selectfont
\tabcolsep=0pt
\begin{tabular*}{\textwidth}{@{\extracolsep{\fill}}lp{340pt}@{}}
\hline
& \textbf{Recursion to construct the candidate optimal band levels}\\
\hline
[0.] & \textbf{Set} $i\leftarrow1$, $\unl a^* \leftarrow
\{0\}$, $\unl b^* \leftarrow\{b^*\}$, $f\leftarrow{} _{b^*_+}
v_b^*$ and $\Xi\leftarrow\Xi^*(f)$, where $\Xi^*(f)$ is given by
\eqref{eq:Xi}.\\
{[1.]} & \textbf{If} $\Xi$ is completely monotone, \textbf{set}
$\unl a^*\leftarrow\unl a^*\cup\{\infty\}$. \textbf{Return} $\{
\unl a,\unl b\}$.\\
{[2.]} & \textbf{Else if} $K=0$ or \textbf{if} $\{K>0$ and $\a
^*_{f,\varnothing
}\ge\a^*_f\}$
\textbf{define} 
$
 (a^*_{i+1}, b^*_{i+1}  ) \leftarrow (b^*_{i,+} + \a
^*_f,
b^*_{i,+} +\b^*_f  )$,\\
&
where the levels $\a^*_f$ and $\beta^*_f$ are defined in \eqref
{eq:astar2} and \eqref{eq:bstar2}.\\
& \textbf{Else if} $\{K>0$ and $\a^*_{f,\varnothing}<\a^*_f\}$
\textbf
{define} 
$
 (a^*_{i+1}, b^*_{i+1}  ) \leftarrow (b^*_{i,+} + \a
^*_{f,\varnothing},  \{b^{**}_{i,-},
b^*_{i,+} + \b^*_{f,\varnothing}  \}  )$\\
& with 
$b^{**}_{i,-}=\inf \{b^*_{i,-}\dvtx V_{\unl a^*,\unl
b^*}(b^*_{i,+} +
\b
^*_{f,\varnothing}) -
V_{\unl a^*,\unl b^*}(b^*_{i,-}) = \b^*_{f,\varnothing} + b^*_{i,+} -
b^*_{i,-} - K  \}$,\\
& where the levels $\a^*_{f,\varnothing}$ and $\beta
^*_{f,\varnothing}$
are defined in \eqref{eq:bstare}.\\
{[3.]} & \textbf{Set} $\unl a^*\leftarrow\unl
a\cup\{a^*_{i+1}\}$, $\unl b^*\leftarrow\unl b\cup\{b^*_{i+1}\}$,
$f\leftarrow{} _{b^*_{i+1,+}} V_{\unl a^*, \unl b^*}$, $\Xi
\leftarrow\Xi_{\unl a^*, \unl b^*}(f)$, $i \leftarrow i+1$.\\
{[4.]}& \textbf{Go to} step 1.\\
\hline
\end{tabular*}}
%

\begin{Rem}\label{rem:res} There may exist
a limit point $\gamma_*=\lim_{i\to\infty} b_{i,+}^*=\break
\lim_{i\to\infty}a^*_i$ of the
band levels. In this case the procedure will converge
to the value-function
$V_{\unl{\tilde a}^*,\unl{\tilde b}^*}$ corresponding to the
levels $\unl{\tilde a}^*=(a^*_i)$, $\unl{\tilde b}^*=(b^*_{i})$,
and needs to be re-started as follows:
\begin{itemize}
\item[{[0.$^\prime$]}] \textbf{Set} $i\leftarrow1$, $\unl a^*
\leftarrow
\unl{\tilde a}^*$, $\unl b^* \leftarrow\unl{\tilde b}^*$,
$f\leftarrow{} _{\gamma^*}
V_{\unl{\tilde a}^*,\unl{\tilde b}^*}$,
$\Xi\leftarrow\Xi_{\unl{\tilde a}^*,\unl{\tilde b}^*}(f)$.
\end{itemize}
\end{Rem}

In the following result (proved at the end of the section) it is
confirmed that the
constructed candidate policy $\pi_{\unl a^*,\unl b^*}$ is indeed optimal:

\begin{Thm}\label{thm:explopt}
The multi-dividend-band strategy $\pi_{\unl
a^*,\unl b^*}$ is an optimal strategy for the control problem
in \eqref{optdiv} and the optimal value function is given by
$v^*= v_{\pi_{\unl
a^*,\unl b^*}} = V_{\unl a^*, \unl b^*}$, with
%
\begin{equation}
\label{eq:VASTAR}\qquad V_{\unl a^*, \unl b^*}(x):= \cases{ W^{(q)}(x)
C^*_i + F_w(x), & \quad $x \in\bigl[a^*_{i},
b^*_{i,+}\bigr], i\ge 1$,\vspace*{2pt}
\cr
x - b^*_{i,+} +
V_{\unl a_*, \unl b_*}\bigl(b^*_{i,+}\bigr), &\quad  $x\in \bigl(b^*_{i,+},a^*_{i+1}
\bigr), i\ge1$, }
\end{equation}
for some constants $C^*_i$,
where the functions $f_i\dvtx\mbb R_-\to\mbb R$ are given by $f_i(x) =
V_{\unl a^*, \unl b^*}(a^*_{i-1} + x), i >1$,
with $f_1=w$.
\end{Thm}

\begin{Rem}
In Shreve et al. (\cite{Shreve}, page 74), an explicit example is given of
an optimal control problem
in a diffusion setting in which a multi-dividend-band strategy is
optimal with
countably many bands. Azcue and Muler \cite{AM2} provide an example of
an optimal strategy with infinitely many bands below a finite level,
for the classical De Finetti dividend problem with bounded dividend
rates in the setting of a compound Poisson process. It is an open problem
to construct an explicit example
in which a multi-dividend-band strategy with countably many bands is
optimal in the dividend-penalty problem.
\end{Rem}

\subsection{Proof of Theorem \texorpdfstring{\protect\ref{thm:explopt}}{11.3}}
Denote by $\unl v_* = (v_{i,j})_{(i,j)}$, $\unl
a^*=(a^*_{i,j})_{(i,j)}$ and $\unl b^* = (b^*_{i,j})_{(i,j)}$
the sequence of value-functions and band levels
generated by the algorithm in Section~\ref{sec:multialg},
where the index $(i,j)$ refers to the $i$th iteration
of the algorithm in the $j$th run of the algorithm (i.e., it has been
restarted $j-1$ times; cf. Remark~\ref{rem:res}). In particular, it
follows that
$v_{i,j}$ is given by
%
\begin{eqnarray}
\label{eq:vijx} v_{i,j}(x) = \cases{ V_{\unl a^*, \unl b^*}(x), &\quad $x\in
\bigl[0,b^*_{i,j,+}\bigr]$,\vspace*{2pt}
\cr
x - b^*_{i,j,+} +
v_{i,j}\bigl(b^*_{i,j,+}\bigr), &\quad  $x > b^*_{i,j,+}$.}
\end{eqnarray}
In the following result (which implies Theorem~\ref{thm:explopt})
it is established that $\pi_{\unl a^*,\unl b^*}$
is an optimal strategy for \eqref{optdiv}:

\begin{Prop}
\textup{{(i)}} For a given pair $(i,j)$ of iteration and run, $v_{i,j}$ is
equal to the
value-function $v_{\unl a^*_{i,j},\unl b^*_{i,j}}$
of the multi-dividend-band strategy $\pi_{\unl a^*_{i,j},\unl b^*_{i,j}}$
at levels $\unl a^*_{i,j}=(0, a_{1,1}^*,\ldots, a_{i-1,j}^*, \infty)$
and $\unl b^*_{i,j}= (b_{1,1}^*,\ldots, b_{i,j}^*)$.
\begin{longlist}[(iii)]
\item[(ii)] For each pair $(\ell,k)$ that is smaller than $(j,i)$ in the
lexico-graphical order,
$v_{(k,\ell)}(x) = v_*(x)$ for all $x\leq b^*_{k,\ell,+}$.

\item[(iii)] The optimal value function $v_*$ is equal to
the value function $V_{\unl a^*,\unl b^*}$ of the strategy
$\pi_{\unl a^*,\unl b^*}$.
\end{longlist}
\end{Prop}

\begin{pf}
(i) The strong Markov property of the process
$U = U^{\pi_{\unl a^*_{i,j},\unl b^*_{i,j}}}$
applied at the stopping time $\tau= \tau^\pi_{a^*_{i-1,j}}$ implies
the relation
%
\begin{equation}
\label{eq:vijxR} v_{k,\ell}(x) = \E_x \biggl[\int
_{[0, \tau]}\mathrm{e}^{-qt}\mu ^\pi_K(
\td t) + v_{k-1,\ell} (U_{\tau} ) \biggr],
\end{equation}
for $k\leq j$, $\ell\leq i$, with $\pi= \pi_{\unl a^*_{i,j},\unl b^*_{i,j}}$.
As $v_{k,\ell}(x)$ is increasing in $k$, it follows that $v_{\infty,\ell
}(x):= \lim_{k\to\infty}v_{k,\ell}(x)$
exists, for any $\ell\leq j-1$. By applying again the strong Markov
property it follows that $v_{1,\ell+ 1}$ satisfies,
for any $l\leq j-1$, $\pi= \pi_{\unl a^*_{i,j},\unl b^*_{i,j}}$,
%
\begin{equation}
\label{eq:vijxR2} v_{1,\ell+1}(x) = \E_x \biggl[\int
_{[0, \tau]}\mathrm{e}^{-qt}\mu ^\pi_K(
\td t) + v_{\infty,\ell} (U_{\tau} ) \biggr].
\end{equation}
The form of $v_{i,j}$ then follows by induction, starting
from the expression for a single
dividend band strategy and using the form of the value-function
of the auxiliary stochastic control problem in \eqref{eq:optstop}
[subsequently applied with pay-off functions
$f(x) = v_{\pi_{\unl a^*_{k,\ell},\unl b^*_{k,\ell}}}(b^*_{k,\ell,+} + x)$,
and performing induction\vspace*{1pt} in $k$ for fixed $\ell$ and using the relation
\eqref{eq:vijxR2}].

(ii) By induction it follows that, for any $k$, $v_*(x) = v_{(k,1)}(x)$
for all $x\leq b^*_{k,1,+}$.
Indeed, note that Corollary~\ref{twoband} implies
$v_{(2,1)}(x) = v_*(x)$ for all $x\leq b^*_{2,1,+}$.
Furthermore, that the induction step holds is verified as follows:
Assuming that $v_{(k-1,1)}(x) = v_*(x)$ for all $x\leq b^*_{k-1,1,+}$
for some pair $k$, Theorem~\ref{thmos} with $f= {}_{b^*_{k-1,1,+}} v_*$
in conjunction with the relation in \eqref{eq:vijxR}
implies that $v_{(k,1)}(x)=v_*(x)$ for $x\leq b^*_{k,1,+}$.

The assertion in (ii) thus follows by induction in $\ell>1$,
following a line of reasoning that is analogous to the one applied in
the previous paragraph
but with the function $w$ replaced by $v_{\infty,\ell-1}$.

(iii) Since $v_{i,j}(x) = V_{\unl a^*, \unl b^*}(x)$
for all $x\leq a^*_{i-1,j}$ [from \eqref{eq:vijx}],
it follows by virtue of part (ii) that
$v_*(x) = V_{\unl a^*, \unl b^*}(x)$
for all $x\leq a^*_{i-1,j}$. Since the sequence $(a_{i,j})_{i,j}$
is strictly increasing and ultimately tends
to infinity (cf. step 2 of the algorithm
and Lemma~\ref{lem:a}), it follows that
$v_*(x)$ is equal to $V_{\unl a^*,\unl b^*}(x)$, for any fixed
$x\in\mbb R_+$.
\end{pf}

\section{Existence and uniqueness of stochastic solutions}\label{sec:ex}

In this section the optimal value function $v_*$, which was identified
in the previous section,
is shown to be a stochastic solution of the HJB equation \eqref{eq:HJB}.
From the form \eqref{eq:VASTAR} and properties of $W^{(q)}$ and of
Gerber--Shiu functions,
it follows that $v_*(x)$ is left- and right-differentiable at any $x>0$.
Furthermore, it was shown in Lemma~\ref{lem:est} that $v_*(x)$ is
continuous at any $x\in\mbb R_+$.
In particular, the function $g=v_*$ is continuous and
left-differentiable at the ``right-boundary''
$\partial^+\mc C_g:=\{b_1, b_2, \ldots\}$
of the set
$\mc C_{g}$ (which was defined in \eqref{eq:CD} and
where the interior $\mc C^o_g$ of $\mc C_g$ is denoted by $\mc C^o_g =
\bigcup_n(a_n,b_n)$
for some $a_n,b_n\in[0,\infty]$ with $a_n<b_n$) and thus satisfies the
following property:
%
\begin{equation}\qquad
\label{eq:unism} \mbox{If $K=0$, $g(x)$ is continuous and left-differentiable at
any $x\in\partial^+\mc C_g$.}
\end{equation}

The HJB equation \eqref{eq:HJB} admits a unique stochastic solution
satisfying the regularity condition \eqref{eq:unism}:

\begin{Thm}\label{thm:repg2}
The value function $v_*$ is the unique stochastic solution
of the HJB equation \eqref{eq:HJB} satisfying \eqref{eq:unism}.
\end{Thm}

\begin{pf*}{Proof (Existence)}
As $v_*$ is a stochastic supersolution [by Remark \ref{rem:stst}(i)]
and $v_*$ satisfies \eqref{eq:unism} (as discussed in above paragraph),
it suffices to show that $v_*$ is also a stochastic subsolution.

Note that, in view of the form \eqref{eq:VASTAR}, the interior $\mc
C^o_{v_*}$ of the set $\mc C_{v_*}$
is identified as $\mc C^o_{v_*} = \bigcup_n (a^*_n, b^*_{n,+})$.
Therefore, in view of \eqref{eq:VASTAR} and the martingale properties
of $W^{(q)}$ and of the Gerber--Shiu functions (Proposition~\ref{prop:mart}),
Doob's optional stopping theorem implies that $v_*$ is a local
stochastic subsolution
of the HJB equation \eqref{eq:HJB} on any closed interval $I\subset
\mc
C_{v_*}$, which shows that
$v_*$ is a stochastic subsolution.
\end{pf*}

\subsection{Proof of uniqueness}

Given a stochastic supersolution $g$ of the HJB equation,
an admissible candidate optimal strategy $\pi(g)$
can be described as follows:

\begin{Def}\label{def:optpi}
To a stochastic solution $g$ of
HJB equation \eqref{eq:HJB} are associated:
\begin{enumerate}[(iii)]
\item[(i)] the policy $\pi(g)=\{D^{\pi(g)}_t, t\in\mbb R_+\}\in\Pi$,
given in
terms of
the sets $\mc C_g$ and $\mc D_g:=\mbb R_+\setminus\mc C_g$,

\item[(ii)] the controlled process $U = U^{\pi(g)}$ and

\item[(iii)] the level $y^*(v):= \sup\{u\in[0,v]\dvtx g(v)- g(v-u) +K = u\} $
(with $\sup\varnothing= 0$), that are specified as follows:
\begin{enumerate}[(a)]
\item[(a)] In the case $K=0$, let $D=D^{\pi(g)}$ be the increasing
right-continuous
\mbox{$\mbf F$-adapted} process that satisfies
\[
\cases{ U_t = X_t - D_t \in\ovl{
\mc C}_g,&\quad $\mbox{for any $t\in[0,\tau ^{\pi
(g)})$},$
\vspace*{2pt}
\cr
\displaystyle\int_{[0,\tau^{\pi(g)})}\mbf1_{ \{s\dvtx X_s - D_{s^-}
\notin\ovl{\mc D}_g  \}}(t) \,\td D_t =
0, &}
\]
where $\mbf1_A$ denotes the indicator function of the set $A$ and
$\ovl{\mc C}_g$ and $\ovl{\mc D}_g$ denote the closures of $\mc C_g$
and $\mc D_g$;
\item[(b)] in the case $K>0$, pay out $\Delta D_t = y^*(X_t - D_{t^-})$
at time $t$ if $X_t - D_{t^-}\in\mc D_g$
and $y^*(X_t - D_{t^-}) > 0$;
\item[(c)] otherwise, pay no dividends.
\end{enumerate}
\end{enumerate}
\end{Def}

\begin{Rem}\label{rem:skor}
The Skorokhod embedding lemma implies that the strategy
$\pi(g) = \{D^{\pi(g)}_t, t\in\mbb R_+\}$ described in
Definition~\ref
{def:optpi}(iii)(a)
is equal to
\[
D^{\pi(g)}_t = \sup_{s\in[0,t\wedge\tau^{\pi(g)}]}
\bigl(X_s - b(s)\bigr)\vee0, 
\qquad b(s)=b_{\iota(s)}
\]
with $\iota(s) = \inf\{n\in\mbb N\dvtx X_s-D^{\pi(g)}_{s^-} \leq
a_n\}$,
given the representation $\ovl{\mc D}_g = \bigcup_{n\ge1}[b_n,a_n]$.
In particular, it follows that the policy defined in
Definition~\ref{def:optpi}
is a multi-dividend band strategy.
\end{Rem}

\begin{lemma}\label{lem:shift2}
Let $g$ be a stochastic solution of
the HJB in \eqref{eq:HJB} satisfying~\eqref{eq:unism}. Then the process
$\WT M^{g,\pi_*,\tau_{\mbb R_+}^{\pi_*}}$
with $\pi_* = \pi(g)$, defined in Lemma~\ref{lem:shift} and
Definition~\ref{def:optpi},
is a UI $\mbf F$-submartingale.
\end{lemma}

The proof of Lemma~\ref{lem:shift2} is based on the following
auxiliary result:

\begin{lemma}\label{lemma:reflect}
Let $a>0$ be given and
suppose that the function $g\dvtx\mbb R\to\mbb R$ is such that
$g|_{\mbb
R_-}\in\mc P$,
$g|_{\mbb R_+}$ is c\`{a}dl\`{a}g, and $g$ is continuous and
left-differentiable at $a>0$.
If $M=\{M_t, t\in\mbb R_+\}$ with
$
M_t = \mathrm{e}^{-q(t\wedge T_{0,a})}g(X_{t\wedge T_{0,a}})
$
is an $\mbf F$-martingale, then $Z=\{Z_t, t\in\mbb R_+\}$ with
\[
Z_t = \mathrm{e}^{-q(t\wedge\tau_0)}g \bigl(Y^a_{t\wedge\tau
_0}
\bigr) - g\bigl(Y^a_0\bigr) - g_-'(a)\int
_{[0, t\wedge\tau_0]}\mathrm{e}^{-qs}\,\td\ovl X^a_s
\]
is an $\mbf F$-martingale, where $g'_-(a)$ denotes the left-derivative
of $g$ at $a$.
\end{lemma}

The proof of this result rests on an application of It\^{o}'s
lemma and
a density argument. Details are omitted since these follow straightforwardly
from \cite{NY}, Proposition~1.

\begin{pf*}{Proof of Lemma~\ref{lem:shift2}}
The proof is a modification of the proof of Lemma~\ref{lem:shift}.
As, by Lemma~\ref{lem:shift}, $\WT M^{g,\pi(g)}$ is a UI
supermartingale, it suffices
to verify that $\WT M^{g,\pi(g)}$ is in fact a martingale.
Note that the set of distinct epochs
$\tilde{\mbb T}$ at which lump-sum dividend payments occur is countable,
\[
\tilde{\mbb T}=\{\tilde T_i\dvtx\Delta D_{\tilde T_i} > 0 \}
\qquad\mbox{with $\tilde T_i=\inf\bigl\{t> \tilde T_{i-1}\dvtx
X_t- D_{t-}^{\pi(g)} \in\mc D_g\bigr
\}$,}
\]
for $i\in\mbb N$ with $\tilde T_0=0$ and $\inf\varnothing=\infty$. The
form of the strategy $\pi(g)$ implies that the sequence
$(U_{\tilde T_i})_i$ is decreasing with $U_{\tilde T_i} - U_{\tilde T_{i-1}}>0$
on the set $\{\tilde T_i<\infty\}$. In particular, it follows that,
also in this case, $\tilde{\mbb T}$ is countable.

Writing $D = D^{\pi(g)}$ and $M = \WT M^{g,\pi(g)}$, fixing arbitrary
$t,s\in\mbb R_+$ with $s<t$ and denoting
$T_i = \tilde T_i\wedge t$, we have $M_t = \sum_{i\ge1} Y_i + \sum_{i\ge0} Z_i$ with $Y_i$ given by
%
\begin{equation}
\label{eq:Yi}\qquad \mathrm{e}^{-q T_i}g(X_{T_i} - D_{T_{i-}})
- \mathrm{e}^{-q T_{i-1}} g(X_{T_{i-1}} - D_{T_{i-1}}) -\int
_{(T_{i-1}, T_i)}\mathrm{e}^{-qs} \,\td D_s,\hspace*{-20pt}
\end{equation}
and $Z_i=\mathrm{e}^{-q T_i}(g(X_{T_i} - D_{T_i}) - g(X_{T_i} -
D_{T_i-}) +
\Delta D_{i} - K)\mbf1_{\{\Delta D_i>0\}}$
with $\Delta D_i = D_{T_i}-D_{T_{i-1}}$. By definition of the strategy
$\pi(g)$ it is straightforward to verify that $Z_i=0$ for all $i$.

In the case $K>0$ the integral term in \eqref{eq:Yi}
vanishes, and we have $D_{T_{i-1}}=D_{T_i-}$ for $i\ge0$. By reasoning
as in Lemma~\ref{lem:shift} it follows
that the equality in \eqref{eq:YEDOOB} holds.
By combining \eqref{eq:YEDOOB} with the fact that $g$ is a stochastic
solution,
Doob's optional stopping theorem and the definition of $T_i$, we have
\[
\E[Y_i|\mc F_{T_{i-1}}] = \mathrm{e}^{-q T_{i-1}}
\E_{U_{T_{i-1}}} \bigl[\mathrm{e}^{-q\tau_i}g (X_{\tau
_i} ) -
g(X_0) \bigr] = 0,
\]
with $\tau_i = T_i\circ\theta_{T_{i-1}}$. The tower property hence
yields $\E[M_t- M_s|\mc F_s] =0$.
Since $s,t$ were arbitrary, it thus follows that $M$ is a martingale.

If $K=0$,
the definition of $\pi(g)$ implies that
the process $\{U_{T_{i-1} +t}, t< T_{i}-T_{i-1}\}$ conditional on $\mc
F_{{T_{i-1}}}$
has the same law as the process $\{Y^b_t, t < \tau_b(a)\}$
with $X_0=b=U_{T_{i-1}}$ and
$\tau_b(a) = \inf\{t\ge0\dvtx Y_t^b < a\} $,
conditional on $U_{T_{i-1}}$,
where $Y^b$ is independent of $U_{T_{i-1}}$.
The strong Markov property of $Y^a$ implies that $\E[Y_i|\mc
F_{T_{i-1}}]$ is equal to
\[
\mathrm{e}^{-q T_{i-1}}\E_{U_{T_{i-1}}} \biggl[\mathrm{e}^{-q\tau
_b(a)}g
\bigl(Y^b_{\tau
_b(a)}\bigr) - g(Y_0) - \int
_{(0,\tau_b(a))} \mathrm{e}^{-qs} \,\td\ovl
X^b_s \biggr].
\]
This expectation is positive in view of Lemma~\ref{lemma:reflect} and
the fact that
$g_-'(a) \ge1$ [as $\mathtt d_g(a)\ge1$ and $g$ is
left-differentiable at $a$].
Again, an application of the tower property yields $\E[M_t- M_s|\mc
F_s] \ge0$,
and it follows that, in this case, $M$ is a sub-martingale.
\end{pf*}

The stated uniqueness follows as a consequence of the following
comparison principle:

\begin{Prop}
\label{prop:stcomp}
Let $h$ be any stochastic subsolution satisfying \eqref{eq:unism},
and let $g$ be any stochastic supersolution of the HJB equation \eqref{eq:HJB}.
Then $g\ge h$.
\end{Prop}

\begin{pf*}{Proof of Theorem~\ref{thm:repg2} (uniqueness)}
Let $h$ be any stochastic solution of the HJB equation. Since, by the
dual representation in Proposition~\ref{thm:repg}, $v_*$ is the minimal
stochastic supersolution of the HJB and $h$ is a stochastic
supersolution, it follows $v_*\leq h$. Furthermore, the stochastic
comparison principle in
Proposition~\ref{prop:stcomp} implies $v_*\ge h$ (as $h$ and $v_*$ are
stochastic sub- and supersolutions of the HJB).
Thus it holds $v_*=h$, and uniqueness is established.
\end{pf*}

\begin{pf*}{Proof of Proposition~\ref{prop:stcomp}}
Let $g$ and $h$ be a stochastic supersolution and stochastic
subsolution, and denote by $\pi(h)$ the policy corresponding to $h$
given in Definition~\ref{def:optpi}.
Since the processes $\WT M^{v_*,\pi(h)}$ and $\WT M^{h,\pi(h)}$
[defined in
\eqref{eq:Mtildegpi}],
are a supermartingale and a submartingale (by Lemmas \ref{lem:shift}
and \ref{lem:shift2}),
Doob's optional stopping theorem implies for $ x\in\mbb R_+$
%
\begin{equation}
\label{keyeq} v_*(x) - h(x) \geq\lim_{t\to\infty} \E_x
\bigl[\WT M^{v_*,\pi(h)}_{t\wedge\tau^{\pi(h)}} - \WT M^{h,\pi
(h)}_{t\wedge\tau^{\pi(h)}}
\bigr].
\end{equation}
The right-hand side of \eqref{keyeq} is equal to 0, since $\WT M^{v_*,\pi(h)}$ and
$\WT M^{h,\pi(h)}$ are UI,
and satisfy the boundary condition
\[
\WT M^{v_*,\pi(h)}_{\tau^{\pi(h)}} = \WT M^{h,\pi(h)}_{\tau^{\pi(h)}} =
\mathrm{e}^{-q\tau^{\pi(h)}}w\bigl(U^{\pi(h)}_{\tau^{\pi(h)}}\bigr),
\]
and $\P_x(\tau^{\pi(h)}<\infty)=1$ for all $x\in\mbb R_+$. This
completes the proof.
\end{pf*}

\section{Examples}\label{sec:exa}
\subsection{General computations for processes with rational Laplace exponent}

The determination of the optimal policy starts with
the identification of the last global maximum of the barrier influence
function $G$. For example,\vadjust{\goodbreak} in
the presence of an exponential penalty $w(x)=c\mathrm{e}^{vx}$ or a linear
penalty $w(x)=c x + c_0$, we must compute the extrema of the functions
%
\begin{equation}\qquad
\label{eq:GvG1} G^{(v)}(x):= \frac{1 - c Z^{(q,v)\prime}(x)}{W^{(q)\prime}(x)},\qquad
 G_1(x):=
\frac{1 - c Z'_1(x)- c_0 q W^{(q)}(x)}{W^{(q)\prime}(x)},
\end{equation}
respectively.

Therefore, the first step will be computing the homogeneous and
generating scale functions $W^{(q)}(x)$, $Z^{(q,v)}(x)$, for processes
with rational Laplace exponent.
Assume the typical case
\[
W^{(q)}(x) = \sum A_i \mathrm{e}^{\zeta_i(q) x},
\]
with $A_i\in\mbb R$ and the roots $\zeta_i(q)$ of the Cram\'
er--Lundberg equation
$\psi(\zeta)=q $ being distinct.

This implies $Z^{(q,v)}(x) = \mathrm{e}^{vx}(1 + (q-\psi(v))\int_0^x
\mathrm{e}^{-v
y}W^{(q)}(y)\,\td y)$ is equal to
\[
\mathrm{e}^{vx} + \bigl(q-\psi(v)\bigr)\sum
_i A_i \frac{\mathrm{e}^{\zeta_i(q) x}-\mathrm{e}^{vx}}{\zeta_i(q)-v} = \bigl(\psi(v)-q
\bigr)\sum_i \frac{A_i}{v-\zeta_i(q)}
\mathrm{e}^{\zeta_i(q) x},
\]
using that $ \sum\frac{A_i}{v-\zeta_i(q)
}=\frac{1}{\psi(v)-q}$. In particular, $Z^{(q)}(x) = q \sum_i A_i
\frac
{\mathrm{e}^{\zeta_i(q) x}}{\zeta_i(q)}$ and
\begin{eqnarray*}
Z_1(x) &=& \ovl Z^{(q)}(x) - \psi'(0) \ovl
W^{(q)}(x)= q \sum_i A_i
\frac{\mathrm{e}^{\zeta_i(q) x}}{\zeta_i^2(q)}- \psi'(0) \sum_i
A_i \frac{\mathrm{e}^{\zeta_i(q) x}}{\zeta_i(q)},\label{Z1}
\\
Z^{(q,v)}(x)& =& Z^{(q)}(x)+ \sum_i
A_i {\mathrm{e}^{\zeta_i(q)
x}}\frac
{v}{v-\zeta_i(q)} \biggl(
\frac{\psi(v)}{v}-\frac{q}{\zeta
_i(q)} \biggr).
\end{eqnarray*}

The simplest examples may be completely analyzed by studying the sign
of the functions that are given by $D^\#(x)=-G^{\#\prime}(x)
W^{(q)\prime}(x)^2$, and $D^*(x)=-G^{*\prime}(x) W^{(q)\prime}(x)^2$,
which determine the critical point $b^*$ (in particular whether it is
$0$), and the eventual unimodality
after $b^*$, which implies
optimality of the single barrier policy. To alleviate notation, the $\#,*$ will be omitted in this section, since the function considered can
always be inferred from the absence/presence of transaction costs.

For exponential and affine penalties, the corresponding functions are
given by
$D^{(v)}(x)=-G^{(v)\prime}(x) W^{(q)\prime}(x)^2$ and $D_1(x)
=-G_1^{\prime}(x) W^{(q)\prime}(x)^2$. By straightforward calculations
we find
\begin{eqnarray*}
D^{(v)}(x)&=&W^{(q)\prime\prime}(x) \bigl(1-c Z^{(q,v)\prime}(x)
\bigr) + c Z^{(q,v)\prime\prime}(x) W^{(q)\prime}(x)
\\
&=& \sum_j A_j
\zeta_j(q)^2\mathrm{e}^{\zeta_j(q)x} + c\bigl(\psi
(v)-q\bigr)\sum_j \sum
_{k>j} d^{(v)}_{j,k} A_jA_k
\mathrm{e}^{(\zeta_j(q) + \zeta_k(q))x},
\\
D_1(x) &=& \sum_j A_j
\zeta_j(q)^2\mathrm{e}^{\zeta_j(q)x} - c q \sum
_j \sum_{k>j}d_{1;j,k}
A_jA_k \mathrm{e}^{(\zeta_j(q) + \zeta_k(q))x}
\\
&&{}+ \bigl(c \psi '(0) - c_0 q\bigr) \sum
_j \sum_{k>j} \bigl(
\zeta_j(q)-\zeta_k(q)\bigr)^2A_jA_k
\mathrm{e}^{(\zeta_j(q) + \zeta_k(q))x},
\end{eqnarray*}
with $d^{(v)}_{j,k}\frac{\zeta_j(q)\zeta_k(q)(\zeta_j(q)-\zeta
_k(q))^2}{(v-\zeta_j(q))(v-\zeta_k(q))}$
and $d_{1;j,k} = \frac{(\zeta_j(q)+\zeta_k(q))}{\zeta_j(q)\zeta_k(q)}
(\zeta_j(q)-\zeta_k(q))^2$.

[Note that the coefficients of $c$ and $c \psi'(0)- c_0 q$
are the intervening Wronskians, and that the function
$D^{(v)}(x)-W^{(q)\prime\prime}(x)$ is a generating function for the
corresponding functions obtained with polynomial penalties.]
\subsection{Cram\'{e}r--Lundberg model with exponential jumps}

Consider next the Cram\'{e}r--Lundberg model
(\ref{cramermod}) with exponential jump sizes with mean $1/\mu$, jump
rate $\lambda$, and Laplace exponent $\psi(s)=p s-\lambda
s/(\mu+s)$. The homogeneous scale function is
\[
W^{(q)}(x) = A_+ \mathrm{e}^{\zeta^+(q)x} - A_- \mathrm{e}^{\zeta^-(q)x},
\]
where $A_\pm= p^{-1} \frac{\mu+
\zeta^\pm(q)}{\zeta^+(q)-\zeta^-(q)}$, and $\zeta^+(q)=\F(q)$,
$\zeta^-(q)$ are the largest and smallest roots of the polynomial
$(\psi(s)-q)(s + \mu)= p s^2 + s(p \mu-\lambda-q) - q \mu$:
\[
\zeta^\pm(q)=\frac{{q} + \lambda-\mu p \pm
\sqrt{({q} + \lambda-\mu p)^2 + 4 p {q}\mu}}{2 p}.
\]

Hence, it follows
\begin{eqnarray*}
Z^{(q)}(x) &=& {q} \biggl(\frac{A_+}{\zeta^+(q)} \mathrm{e}^{\zeta^+(q)x}
- \frac{A_-}{\zeta^-(q)} \mathrm{e}^{\zeta^-(q)x} \biggr)
\\
&=& \frac{(q-\zeta^-(q)) \mathrm{e}^{\zeta^+(q)x}+(\zeta^+(q)-q)
\mathrm{e}^{\zeta^-(q)x}} {\zeta^+(q)-\zeta^-(q)},
\\
Z^{(q,v)}(x) &=& Z^{(q)}(x)+ \lambda\frac{v}{v+ \mu}
\frac{ \mathrm{e}^{\zeta^+(q)x}-
\mathrm{e}^{\zeta^-(q)x}} {\zeta^+(q)-\zeta^-(q)},
\\
D^{(v)}(x) &=& \a_+\mathrm{e}^{\zeta^+(q) x} - \a_-
\mathrm{e}^{\zeta^-(q) x} + c\a_{v}\mathrm{e}^{(\zeta^+(q)+\zeta^-(q)) x},
\end{eqnarray*}
with $\a_+ = A_+ (\zeta_+(q))^2 > 0$, $\a_- = A_-(\zeta_-(q))^2> 0$,
$C = (\mu+\break
\zeta_+(q))(\mu+\zeta_-(q))=\frac{\lambda\mu}{p}
>0$, and
\[
\a_{v} = \frac{p}{v+\mu} \frac{C}{p^2} \frac{q\mu}{p} =
\frac{\lambda q \mu
^2}{p^3} \frac{1}{v+\mu}>0.
\]

Then, differentiating $v\mapsto Z^{(q,v)}(x)$, $v\mapsto\a_v$
or by \eqref{Z1} and using that $(\zeta^+(q) +\zeta^-(q)
)/(\zeta^+(q) \zeta^-(q) )=\psi'(0)/q -1/\mu$ yields
\begin{eqnarray*}
Z_1(x)&=& \lambda\mu^{-1} \frac{ \mathrm{e}^{\zeta^+(q)x}-
\mathrm{e}^{\zeta^-(q)x}} {\zeta^+(q)-\zeta^-(q)} =
C_{+}\mathrm{e}^{\zeta^+(q)x} + C_{-}
\mathrm{e}^{\zeta^-(q)x},
\\
D_1(x) &=& \a_+\mathrm{e}^{\zeta^+(q) x} - \a_-
\mathrm{e}^{\zeta
^-(q) x} + \a_{1} \mathrm{e}^{(\zeta^+(q)+\zeta^-(q)) x},
\end{eqnarray*}
where $C_{\pm} = \pm\lambda\mu^{-1}
(\zeta^+(q)-\zeta^-(q))^{-1}$ and
\begin{eqnarray*}
\a_{1} &=& A_+ A_- \bigl(\zeta^+ -\zeta^-\bigr)^2
\biggl(c q \frac
{\zeta^+ +\zeta^-}{\zeta^+ \zeta^-}- c\psi'(0) + c_0 q
\biggr)
\\
&=& \frac
{C}{p^{2}} \biggl( c_0 q -c \frac{q}{\mu}\biggr)=
\frac{\lambda q}{p^3} (c_0 \mu - {c}).
\end{eqnarray*}

Recall next that in the absence of penalty and costs
[$w(x)=K=0$], the function $W^{(q)\prime}(x)=G(x)^{-1}$ is
unimodal (see Avram et al. \cite{APP}) with global minimum at $b^*$
given by
\[
\frac{1}{\zeta^+(q) - \zeta^-(q)}
\cases{\displaystyle \log \frac{\zeta^-(q)^2(\m+\zeta^-(q))}
{\zeta^+(q)^2(\m+\zeta^+(q))}, \vspace*{2pt}\cr \qquad
\hspace*{10pt}\mbox{if
$W^{(q)\prime\prime}(0) < 0 \Eq(q+\lambda)^2< p\lambda\m$,}
\vspace*{2pt}
\cr
0,\qquad  \mbox{if $W^{(q)\prime\prime}(0)\geq0 \Eq (q+
\lambda)^2\leq p\lambda\m$}.}
\]
[Since
$W^{(q)\prime\prime}(0) \sim
{\zeta^+(q)^2(\m+\zeta^+(q))}-\zeta^-(q)^2(\m+\zeta^-(q))/(\zeta
^+(q)-\zeta^-(q))
= (q+\lambda)^2- p\lambda\m$, the optimal strategy is
always the barrier strategy at level $b^*$.]\looseness=-1

It is verified next that the functions $G^{(v)}$ and $G_1$ continue to
be unimodal when $w$ is exponential or affine and $K=0$, as a
consequence of Lemma~\ref{lem:uni} below, and hence single barrier
policies continue to be optimal, in view of
Lemma~\ref{Cor:one} (in the case of affine penalties
this has already been established in \cite{LR,APPjaca}).

\begin{lemma}\label{lem:uni} Let $\a_i,\lambda_i\in\mbb R$,
$i=1,2,3$ satisfy
$\a_1>0>\a_3$, and $\lambda_1>\lambda_2>\lambda_3$. Then the
function $f(x):=\sum_{i=1}^3\a_i\mathrm{e}^{\lambda_i x}$ has a unique
root $c^*$ of $f(c^*)=0$, and it holds $f'(c^*)>0$, and
\[
f(x) > 0 \qquad\mbox{for all $x> c^*$}.
\]
Furthermore, if $h\dvtx\mbb R_+\to\mbb\R$ is such that $h'(x)=k(x) f(x)$
for $x>0$, where $k\dvtx\mbb R_+\to
\mbb R_+\setminus\{0\}$,
then $h$ is
unimodal.
\end{lemma}

\begin{pf}
The function $g(x):= \mathrm{e}^{-\lambda_3 x}f(x)$ tends to
$+\infty$ and to $\alpha_3<0$ as $x\to\pm\infty$. If it holds $\a
_2\ge0$,
$g$ is strictly convex and strictly increasing. In the case $\a_2<0$, $g$
attains a minimum at the unique root of $g'$. In both cases
the equation $g(c)=0$ admits a unique root $c$, and it holds $g'(c)>0$.
Hence it holds
that $c$ is a unique root of $f(c)=0$, with $f'(c)>0$ and with
$f(x)>0$ for $x>c$. In particular, $h$ has a unique stationary
point where it attains a maximum, so that it is unimodal.
\end{pf}

The optimal level
$b^*$ is characterized as follows:

(i) For $K=0$ and in the case of an exponential penalty, $b^*_{v,+}=0$
if and only if
\[
G^{(v)\prime}(0)\leq0 \Eq(q+\lambda)^2 - \lambda\mu p \ge -c
\lambda q\frac{\mu^2}{v+\mu},
\]
as follows from the expression for $D^{(v)}(x)$. Similarly, in the case
of linear
penalty, it holds $b^*_{1,+}=0$ if and only if
\[
G_1'(0)\leq0 \Eq(q+\lambda)^2 - \lambda\mu
p \ge \lambda q (c - c_0 \mu),
\]
in view of the expression for $D_1(x)$. If $b^*_{+}$ is positive, it is
a stationary point, and hence
solves the equation
\[
G^{(v)\prime}(b) = 0 \Eq0 = D^{(v)}(b) = \a_+\mathrm{e}^{\zeta
^+(q) b}
- \a_- \mathrm{e}^{\zeta^-(q) b} + c\a_{v}\mathrm{e}^{(\zeta^+(q)+\zeta^-(q)) b},
\]
if the penalty $w$ is exponential and
\[
G_1'(b)=0 \Eq0 = D_1(b) = \a_+
\mathrm{e}^{\zeta^+(q) b} - \a_- \mathrm{e}^{\zeta^-(q) b} + \a_{1}
\mathrm{e}^{(\zeta^+(q)+\zeta^-(q)) b},
\]
if $w$ is an affine penalty.

(ii) Suppose next $K>0$. Then $b^*_{+}$ is strictly positive
as a consequence of the positive transaction cost $K$, and the
optimal levels $(b_{-}^*, b^*_{+})$ are given by $(b_{-}^*,
b^{*}_-+d^*)$ where $(b,d)$ maximizes over $(b,d)\in\mathbb
R_+\times\mbb R_+\setminus\{0\}$ the function
\[
\WT G^{(v)}\dvtx(b,d)\mapsto\frac{d - K - B_+\mathrm{e}^{\zeta
^+(q)b}(\mathrm{e}^{\zeta^+(q)d}-1)
+
B_-\mathrm{e}^{\zeta^-(q)b}(\mathrm{e}^{\zeta^-(q)d}-1)}{A_+\mathrm
{e}^{\zeta^+(q)b}(\mathrm{e}^{\zeta
^+(q)d}-1)-
A_-\mathrm{e}^{\zeta^-(q)b}(\mathrm{e}^{\zeta^-(q)d}-1)}
\]
if $w$ is an exponential penalty, and the function
\[
\WT G_1\dvtx(b,d)\mapsto\frac{d - K -
C_{+}\mathrm{e}^{\zeta^+(q)b}(\mathrm{e}^{\zeta^+(q)d}-1) +
C_{-}\mathrm{e}^{\zeta^-(q)b}(\mathrm{e}^{\zeta
^-(q)d}-1)}{A_+\mathrm{e}^{\zeta
^+(q)b}(\mathrm{e}^{\zeta
^+(q)d}-1)-
A_-\mathrm{e}^{\zeta^-(q)b}(\mathrm{e}^{\zeta^-(q)d}-1)}
\]
if $w$ is an affine penalty.

The following result sums up the form of the optimal dividend
policy:

\begin{lemma} Consider a Cram\'{e}r--Lundberg process
(\ref{cramermod}) with exponential jump sizes with mean $1/\mu$,
and fixed cost $ K\geq0$. The optimal dividend policy is given by
a \emph{single dividend-band strategy} $\pi_{b^*}$ for the
following Gerber--Shiu penalties $w$:
\begin{enumerate}[(a)]
\item[(a)] Exponential penalties: $w(x)=c \mathrm{e}^{x v}$, with $v,c< 0$
such that
the integrability condition \eqref{eq:cw2}
is satisfied.
\begin{enumerate}[(iii)]
\item[(i)] In the case $\{K=0$ and $(q+\lambda)^2 - \lambda\mu p \ge
-c\lambda q\frac{\mu^2}{v+\mu}\}$, then $b^*=0$.
\item[(ii)] In the case
$\{K=0$ and $(q+\lambda)^2 - \lambda\mu p < -c\lambda
q\frac{\mu^2}{v+\mu}\}$, then $b^*$ is the unique solution
$b\in\mbb R_+\setminus\{0\}$ of the equation $D^{(v)}(b)=0$.
\item[(iii)] In the case $K>0$,
we have $b_{+}^*=b_{-}^* + d^*$ where $b_{-}^*$ and $d^*$ maximize over
$b\ge0$, $d>0$, the function $\WT G^{(v)}$.
\end{enumerate}

\item[(b)] Affine penalties: $w(x)=cx + c_0$, with $c\ge0$ and $c_0\leq0$
such that \eqref{eq:cw2} is satisfied.
\begin{enumerate}[(iii)]
\item[(i)] In the case $\{K=0$ and $(q+\lambda)^2 - \lambda\mu p \ge
\lambda q (c - c_0 \mu)\}$, then we have $b^* = 0$.
\item[(ii)]
In the case $\{K=0$ and $(q+\lambda)^2 - \lambda
\mu p < \lambda q (c - c_0 \mu)\}$, then $b^*$ is the unique solution
$b\in\mbb R_+\setminus\{0\}$ of the equation $D_1(b)=0$.
\item[(iii)] In the case $K>0$,
we have $b_{+}^*=b_{-}^* + d^*$ where $b_{1,-}^*\ge0$ and $d^*>0$
maximize over $(b,d)$, the function $\WT G_1$.
\end{enumerate}
\end{enumerate}
\end{lemma}

\subsection{Cram\'{e}r--Lundberg model with Erlang jumps}
Suppose next that $X$ is given by the Cram\'{e}r--Lundberg model
(\ref{cramermod}) with the Erlang $(n,\mu)$ jump sizes.
The corresponding Laplace exponent is $\psi(s) = ps +
\frac{\lambda\mu^n}{(\mu+s)^n} -\lambda$, and by Laplace
inversion it follows that its $q$-scale function is given by
\[
W^{(q)}(x)=\sum_{j=0}^{n}
A_j\mathrm{e}^{\zeta_j(q) x},\qquad A_j =
\frac{(\zeta_j(q)+\mu
)^n}{p\prod_{k\neq j}(\zeta_j(q)-\zeta_k(q))}, \qquad x\ge0,
\]
where
$\zeta_0(q)>0>\zeta_1(q)>-\mu>\zeta_2(q) > \cdots$ are the $n+1$
roots of the Cram\'{e}r--Lundberg equation $\psi(\zeta)=q$.

Let $K=0$ and $w(x)=c\mathrm{e}^{vx}$ an exponential penalty ($c<0$), and
denote by $b$ the point where $G^{(v)}$ attains its maximum. In
general a single dividend-band strategy may not be optimal. A
necessary and sufficient criterion for optimality of $\pi_b$ is
the complete monotonicity of the function $\Xi_v\dvtx(\Phi(q),\infty
)\to\mbb R_+$ given by
\begin{eqnarray*}
 \Xi_v(s) &= &\frac{\psi(s) - q}{s} \cdot\mathrm{e}^{s b}
\int_b^\infty\mathrm {e}^{-s z}
\bigl(W^{(q)\prime}(z)G^*(b) - \bigl[1-F'(z)\bigr] \bigr) \,\td
c,
\\
%
 I(s) &=& s^{-1} \biggl[ p s + \frac{\lambda\mu^n}{(\mu+s)^n}
- \lambda- q \biggr],
\\
I_v(s) &=& I_0(s) - c\sum
_{j>i} k^{(v,q)}_{i,j}(s)
A_jA_i\mathrm{e}^{(\zeta_i(q) + \zeta_j(q))b},
\\
I_0(s) &=& \int_0^\infty
\mathrm{e}^{-s x}\bigl[ W^{(q)\prime
}(b+x)-W^{(q)\prime}(b)\bigr]
\,\td x = \sum_{j=0}^n A_j
k^{(q)}_{1,i,j}(s) \mathrm{e}^{\zeta_j(q)b},
\end{eqnarray*}
where $k^{(v,q)}_{i,j}(s) =
\frac{(\zeta_j(q)-\zeta_i(q)^2(v -\zeta_i(q)-\zeta_j(q))}%
{(s-\zeta_j(q))(s-\zeta_i(q))(v-\zeta_j(q))(v-\zeta_i(q))}$ and
$k^{(q)}_{1,i,j}(s) = \frac{\zeta_j(q)^2}{s(s-\zeta_j(q))}$.
If in addition there is no penalty ($w=0$), the expressions simplify.
If $b$ denotes the value where the minimum of $W^{(q)\prime}$ is
attained, $\pi_b$ is optimal
precisely if $\Xi_0\dvtx(\Phi(q),\infty)\to\mbb R_+$
is completely monotone, where
$
\Xi_0(s) = I(s) \cdot I_0 (s)$.

\textit{The Azcue--Muler example.} Consider next the example in Azcue and
Muller \cite{AM}, with pure Erlang claims of order $n=2$, with $\mu=1$,
$\lambda=10$, $p =\frac{107}{5}$, $q =\frac{1}{10}$,
$\th=\frac{7}{100}$ and Laplace exponent
$\psi(s) - q = ps + \lambda(\frac{\mu}{\mu+s})^2 -\lambda-q
=\frac{p}{(\mu+s)^2}(s +\zeta_1)(s +\zeta_2)(s -\zeta_0)$,
with $\zeta_0 \approx0.0396$, $\zeta_1\approx0.0794$, $\zeta_2
\approx1.4882$.
In addition we consider a linear penalty $w(x)=cx$, $c\in\mbb R_+$. We
will analyze below
four particular cases $c\in\{0,0.2,0.6,1.0\}$. In cases $c\in\{
0.6,1.0\}$
the optimal strategy is a single dividend
band strategy at level $b_1$, while in the cases $c\in\{0, 0.2\}$ it is
optimal to adopt
a two-band strategy with $b_1=0$ (in the case $c=0$ we thus recover the
form of the optimal strategy found in \cite{AM}). The parameters of the
optimal strategies are summarized in
Table~\ref{table} (with $v_2$ denoting the difference of the value
function and the identity $x\mapsto x$ at the end of the nonempty
continuation band).

\begin{table}
\caption{The values of the optimal band levels under a linear penalty $w(x)=c x$}
\label{table}
\begin{tabular*}{300pt}{@{\extracolsep{\fill}}lcccc@{}}
\hline
& \multicolumn{1}{c}{$\bolds{b_1}$} & \multicolumn{1}{c}{$\bolds{v_2}$} &
\multicolumn{1}{c}{$\bolds{a_2}$} &
\multicolumn{1}{c@{}}{$ \bolds{b_2}$} \\
\hline
$c=0$ & 0\phantom{0.} & $2.44$ &$ 1.83$ & $ 10.45$\\
$c=0.2$ & 0\phantom{0.} &$1.72$ &$ 1.90$& $ 10.47$\\
$c=0.6$ & $ 10.96$ & $1.71$ & $\infty$ & $\infty$\\
$c=1.0$ & $ 11.37$ & $1.30$ & $\infty$ & $\infty$\\
\hline
\end{tabular*}
\end{table}

In the cases $c\in\{0.6,1\}$ a plot of the function $G_1$ defined in
\eqref{eq:GvG1}
reveals that $G_1$ is monotone decreasing on the right of the level at which
attains its unique global maximum
which implies the optimal strategy is a single-dividend band strategy
at this level (Theorem~\ref{thm:cm2}).
In the cases $c\in\{0,0.2\}$
a plot of $G_1$ shows that this function attains its global maximum at $0$
but also attains a second local maximum at some strictly positive level,
so that the optimal value function is given by
\[
v(x) = \cases{ x + v_1, & \quad $b_1=0\leq x <
a_2$,\vspace*{2pt}
\cr
F_1 (x-a_1), &\quad $x
\in[a_2,b_2]$,\vspace*{2pt}
\cr
x + v_2, &\quad  $x >
b_2$. }
\]
Here $v_2 = - b_2 + F_1(b_2-a_2)$
and $v_1 = \frac{p - 20c}{q + \lambda} = \frac{214-200c}{101}$
is the value of the strategy (at zero) of paying
all premiums as dividends until the moment the first claim arrives,
which is also
the moment of ruin, and $F_1(x)$ is given by
\begin{eqnarray*}
F_1(x) &=& p (a_2 + v_1) W^{(q)}(x)
- \int_0^x W^{(q)}(x-y)
\bigl[f_{\nu,a_2}(y)\bigr]\,\td y,
\\
f_{\nu,a}(y) &=& \int_0^{a}(a-z +
v_0)k(y+z)\,\td z + c\int_{a}^\infty
(a-z)k(y+z)\,\td z,
\end{eqnarray*}
where $k(y)=\lambda\mu^2 y\mathrm{e}^{-\mu y}$ denotes the L\'{e}vy density
at $y$.

The function $v$ is the value function of a two-band strategy at levels
$(b_0, a_1, b_1)$ with $b_0=0$. The unknowns $a_1, b_1$ are determined
by the optimality equations $F'_1((b_1-a_1)-)=1$ and $F_1''((b_1-a_1)-)
= 0$ which yield the following system
of two nonlinear equations for $a_1$ and $b_1$:
\begin{eqnarray*}
1 &=& p(a_1+v_0)W^{(q)\prime}(b_1-a_1)
- p^{-1}f_{\nu,a_1}(b_1)
\\
&&{}- \int_0^{b_1-a_1} W^{(q)\prime}(b_1-a_1-y)f_{\nu, a_1}(y)\,\td y,
\\
0& =& p(a_1+v_0)qW^{(q)\prime\prime}(b_1-a_1)
- p^{-1}f'_{\nu, a_1}(b_1)
\\
&&{}- W^{(q)\prime}(0) f_{a_1,\nu}(b_1) - \int
_0^{b_1-a_1} W^{(q)\prime\prime}(b_1-a_1-y)f_{\nu, a_1}(y)\,\td y,
\end{eqnarray*}
with $W^{(q)\prime}(0) = \frac{101}{10} \cdot\frac{25}{107^2}$.
The two-band strategies at the levels $(a_1,b_1) = (1.83, 10.45)$ [$c=0$]
and $(a_1,b_1) = (1.90, 10.47)$ [$c=0.2$]
are indeed optimal since it holds
$({} _{b_1} \Gamma^w_{\infty} v - q v)(y)\leq0$ for all $y>b_1$
and $({} _{0} \Gamma^w_{\infty} v - q v)(y)\leq0$ for all $y\in(0,a_1)$.
%

\begin{appendix}\label{app}
\section{Proof of dynamic programming equation}\label{app:dp}
\begin{pf*}{Proof of Lemma~\ref{prop:mart}(ii)}
Fix arbitrary $\pi\in\Pi$, $x\in\mbb R_+$ and $s,t\in\mbb R_+$
with $s<t$. The process $V^\pi_t$ is $\mc F_t$-measurable, and is
UI on account of Lemma~\ref{lem:est}.
Fix arbitrary $\pi\in\Pi$, $x\in\mbb R_+$.
Define by $W^\pi=\{W^\pi_s, s\in\mbb R_+\}$ the value-process
$W_s^\pi= \operatorname{ess.sup}_{\tilde\pi\in\Pi_s}
J_s^{\tilde\pi}$ with
%
\begin{equation}
\label{eq:Wspi} J_s^{\tilde\pi} = \E \biggl[ \int
_{[0, \tau^{\tilde\pi})}\mathrm{e}^{-qu}\mu_K^{\tilde\pi
}(
\td u) + \mathrm{e}^{-q\tau^{\tilde\pi}}w\bigl(U^{\tilde\pi}_{\tau^{\tilde\pi
}}\bigr)
\Big|\mc F_s \biggr],
\end{equation}
where $\Pi_s =  \{\tilde\pi= (\pi,\ovl\pi) = \{D^{\pi,\ovl
\pi}_u,
u\in\mbb R_+\}\dvtx\ovl\pi\in\Pi  \}$, and $D^{\pi,\ovl
\pi}$
is given in terms of
the process $D^{\ovl\pi}(x)$ of cumulative dividends of the strategy
$\ovl\pi$ corresponding to initial capital $X_0=x$ by
\[
D^{\pi,\ovl\pi}_u = \cases{ D_u^\pi, &\quad
$u\in[0,s)$;\vspace*{2pt}
\cr
D_s^\pi+
D_{u-s}^{\ovl\pi}\bigl(U^\pi_s\bigr), &\quad
$u\ge s$. }
\]
It follows that $V^\pi$ is a supermartingale as direct
consequence of the following $\P$-a.s. relations:
\[
\mbox{(a)\quad $V_s^\pi= W_s^\pi$,\qquad
(b)\quad $W_s^\pi\ge\E\bigl[W^\pi_t|
\mc F_s\bigr]$},
\]
where $W^\pi$ is the process defined in \eqref{eq:Wspi}.

\emph{Proof of} (b): The identity follows by classical arguments.
Since the family of random variables $\{J_t^{\tilde\pi}, \tilde\pi
\in\Pi
_t\}$
is directed upwards, it follows from Neveu \cite{Neveu}
that there exists a sequence $\pi_n\in\Pi_t$ such that
$J_t^{\tilde\pi_n}\uparrow W_t^\pi$.
Since $\Pi_t\subset\Pi_s$ it follows that $W_s^\pi$ dominates
$J_s^{\pi_n}=\E[J_t^{\pi_n}|\mc F_s]$, so that monotone convergence
implies that we have
\[
W_s^\pi\ge\lim_n \E
\bigl[J_t^{\pi_n}|\mc F_s\bigr] = \E
\bigl[W_t^\pi|\mc F_s\bigr].
\]

\textit{Proof of} (a): The form of $D^{\tilde\pi}$ implies that,
conditional on $U^\pi_s$,
$\{D^{\tilde\pi}_u-D^{\tilde\pi}_s, u\ge s \}$
is independent of $\mc F_s$. On account of
the Markov property of $X$ it also follows that
conditional on $U^\pi_s$,
$\{U^{\tilde\pi}_u-U^{\tilde\pi}_s, u\ge s \}$
is independent of $\mc F_s$. As a consequence,
we have the following identity on the set $\{s < \tau^\pi\}$:
\begin{eqnarray*}
&&\E \biggl[\int_{[0, \tau^{\tilde\pi})}\mathrm{e}^{-qu}\mu
^{\tilde\pi
}_K(\td u) + \mathrm{e}^{-q\tau^{\tilde\pi}} w
\bigl(U^{\tilde\pi}_{\tau^{\tilde\pi}}\bigr) \Big|\mc F_s \biggr]
\\
&&\qquad= \mathrm{e}^{-qs}\E_{U^{\pi}_s} \biggl[ \int_{[0, \tau^{\ovl\pi})}
\mathrm{e}^{-qu}\mu^{\ovl\pi}_K(\td u) +
\mathrm{e}^{-q\tau^{\ovl\pi}}w\bigl(U^{\ovl\pi}_{\tau^{\ovl\pi
}}\bigr) \biggr]+ \int_{[0,s]}\mathrm{e}^{-qu}\mu^{\pi}_K(
\td u)\\
&&\qquad = \mathrm{e}^{-qs} v_{\ovl\pi}\bigl(U^\pi_s
\bigr) + \int_{[0,s]}\mathrm {e}^{-qu}
\mu^{\pi
}_K(\td u).
\end{eqnarray*}
In particular, $\P_x$-a.s.
the following representation holds true:
\[
J_s^{\tilde\pi} = \mathrm{e}^{-q(s\wedge\tau^{\pi})} v_{\ovl\pi}
\bigl(U^\pi_{s\wedge\tau^{\pi}}\bigr) + \int_{[0, s\wedge\tau^{\pi}]}
\mathrm{e}^{-qu}\mu^{\pi}_K(\td u),
\]
which yields the following $\P_x$-a.s. representation for
$W^\pi_s$:
%
\begin{eqnarray}
\label{eq:Wpi}&& W^\pi_s = \int_{[0, s\wedge\tau^{\pi}]}
\mathrm{e}^{-qu}\mu^{\pi
}_K(\td u)
\nonumber
\\[-8pt]
\\[-8pt]
\nonumber
&&\hspace*{45pt}\qquad{} +
\mathrm{e}^{-q(s\wedge\tau^{\pi})} \mathop{\operatorname{ess.sup}}_{\tilde\pi=(\pi,\ovl\pi)\in\Pi_s}
v_{\ovl\pi}\bigl(U^\pi_{s\wedge\tau^{\pi}}\bigr).
\end{eqnarray}
In view of the definitions of $\Pi_s$ and $v_*$,
the essential supremum in \eqref{eq:Wpi}
is $\P$-a.s. equal to $v_*(U^\pi_{s\wedge\tau^{\pi}})$,
which implies that, $\P$-a.s., $W^{\pi}_s = V^\pi_s$.
\end{pf*}

\section{Proof of properties of value function}\label{sec:lem:est}

\begin{pf*}{Proof of Lemma~\ref{lem:est}(i)}
Let $x>y$. Denote by
$\pi_\varepsilon(y)$ an $\varepsilon$-optimal strategy for the case $U_0=y$.
Then a possible strategy is to immediately pay out $x-y$ and
subsequently to adopt the strategy $\pi_\varepsilon(y)$, so that the
following holds:
\[
v_*(x) \ge x- y - K + v_{\pi_\varepsilon}(y) \ge v_*(y) - \varepsilon+ x- y - K.
\]
Since this inequality holds for any $\varepsilon>0$, the stated lower
bound follows.

To prove the stated continuity we first establish
an upper bound for the difference $v_*(x)-v_*(y)$ with $x>y$.
Let $\tilde\pi_\varepsilon(x)$ denote an $\varepsilon$-optimal
strategy for the case $U_0=x$ for a given $\varepsilon>0$.
Then a possible strategy is to refrain from
paying any dividends until the first time that the reserves
hit the level $x$, and to subsequently follow
the policy $\tilde\pi_\varepsilon$. Hence $v_*(y)$, $x\ge y$, is
bounded below by
\begin{eqnarray*}
&&\frac{W^{(q)}(y)}{W^{(q)}(x)}\bigl(v_{\tilde\pi_\varepsilon}(x) - F_w(x)\bigr) +
F_w(y) \\
&&\qquad\ge\frac{W^{(q)}(y)}{W^{(q)}(x)}\bigl(v^*(x) - \varepsilon -
F_w(x)\bigr) + F_w(y).
\end{eqnarray*}
Rearranging and letting $\varepsilon$ tend to zero yields
the upper-bound
%
\begin{equation}\qquad
\label{eq:upper} v_*(x) - v_*(y) \leq \biggl(1 -\frac{W^{(q)}(y)}{W^{(q)}(x)} \biggr)
\bigl[v_*(x) - F_w(x)\bigr] + F_w(x) -F_w(y).
\end{equation}

In the case $K=0$, continuity of $W^{(q)}|_{\mbb R_+\setminus\{0\}}$,
the lower bound from part (i) and~\eqref{eq:upper}
yield that $v_*$ is continuous on $\mbb R_+$. In the case $K>0$
continuity of $v_*$ on $\mbb R_+$ follows by combining the upper bound
in \eqref{eq:upper} with a different lower bound that is derived next.

For fixed $\varepsilon>0$ and given initial reserves
$U_0 = y$ for some $y>x$, a possible strategy is to adopt
$\tilde\pi_\varepsilon(x)$ until the first moment that the reserves
$U$ fall below
$\delta:=y-x$, and to follow then a waiting strategy $\pi_\varnothing$
(of not paying any dividends).
Taking $\pi= \tilde\pi_\varepsilon(x)$ it follows by the monotonicity
of $w$
that $v_*(y) - v_*(x)$
for $y\ge x$ is bounded below by
\begin{eqnarray*}
&&\E_{y} \biggl[\int_0^{\tau_\d^\pi}
\mathrm{e}^{-qt}\mu _K^\pi(\td t) +
\mathrm{e}^{-q\tau_\d^\pi}w \bigl(U^\pi_{\tau_\d^\pi} \bigr)
\mbf1_{\{
{\tau_\d
^\pi}=\tau_0^\pi\}} + \mathrm{e}^{-q\tau_\d^\pi}v_{\pi_\varnothing}
\bigl(U^\pi _{\tau_\d
^\pi
} \bigr)\mbf1_{\{{\tau_\d^\pi}<\tau_0^\pi\}} \biggr]
 \\
 &&\quad{}-v_*(x) \\
 &&\qquad= \E_{y} \bigl[\mathrm{e}^{-q\tau_\d^\pi} \bigl(w
\bigl(U^\d_{\tau_\d^\pi} \bigr) - w \bigl(U^\d_{\tau
_\d
^\pi
}-
\d \bigr) \bigr)\mbf1_{\{{\tau_\d^\pi}=\tau_0^\pi\}} \bigr] + f_{\varepsilon}(x,y)
\\
&&\qquad\quad{} + v_{\pi}(x)-v_*(x) \ge-\varepsilon+ f_{\varepsilon}(x,y),
\end{eqnarray*}
where $\tau_\d^\pi= \inf\{t\ge0\dvtx U^\pi_t < \d\}$ and
\[
f_{\varepsilon}(x,y) = \E_{y} \bigl[\mathrm{e}^{-q\tau_\d^\pi
}
\bigl(\mc V_w \bigl(U^\pi_{\tau_\d
^\pi} \bigr) - w
\bigl(U^\pi_{\tau_\d^\pi}-\d \bigr) \bigr) \mbf1_{\{{\tau_\d^\pi}<\tau_0^\pi\}}
\bigr].
\]
Assume for the moment that
$f_{\varepsilon}(x,y)$ tends to zero when $\d=y-x$ tends to $ 0$.
Given this assumption and the bound in \eqref{eq:upper}
it follows (since $\varepsilon$ was arbitrary)
%
\begin{equation}
\label{eq:liminf} \liminf_{|x-y|\to0}\bigl[v_*(y) - v_*(x)\bigr]
\geq0.
\end{equation}
Similarly, it can be shown
$\limsup_{|x-y|\to0}[v_*(y) - v_*(x)] \leq0$. Combining the two
limits yields that
$v_*(x)$ is continuous at each $x\in\mbb R_+$.

Finally, the claim that $f_{\varepsilon}(x,y)$ tends to zero is verified.
First, note the estimate
%
\begin{equation}
\label{fe} f_{\varepsilon}(x,y) \leq \Bigl(\sup_{x\in[0,\d]}
\mc V_w(x) - w(-\d ) \Bigr) \E_y\bigl[
\mathrm{e}^{-q\tau^\pi_\d}\mbf1_{\{\tau^\pi_\d
<\tau^\pi_0\}}\bigr].
\end{equation}
If $X$ has unbounded variation, then the
left-continuity of $w$ at zero and the fact $\mc V_w(0+)=w(0)$ combined
with the inequality in equation \eqref{fe}
imply $f_{\varepsilon}(x,y)\to0$ when $\d=y-x\to0$. If $X$ has bounded
variation, $v_{\pi_w}(0)$ is (in general) not equal to $w(0)$, and it
is next shown that the second factor
in equation \eqref{fe} tends to zero if $\delta\to0$.
Note that the policy $\tilde\pi_\varepsilon(x)$, being element of
$\Pi$,
consists of at most countably many dividends payments almost surely.
Denoting the times of the dividend payments
by $\tau_1, \tau_2,\ldots,$ and the values of $U^{\tilde\pi
_\varepsilon(x)}$ at
those times by $U_1, U_2,\ldots,$ the strong Markov property of $X$ implies
\begin{eqnarray*}
\E_y\bigl[\mathrm{e}^{-q\tau^\pi_\d}\mbf1_{\{\tau^\pi_\d<\tau^\pi
_0\}}\bigr]
&=& \sum_i \E_y\bigl[
\mathrm{e}^{-q\tau^\pi_\d}\mbf1_{\{\tau^\pi_\d
<\tau^\pi
_0, \tau^\pi
_\d\in[\tau_i, \tau_{i+1})\}}\bigr]
\\
&\leq& \sum_i \E_y\bigl[
\mathrm{e}^{-q\tau_i}\mbf1_{\{\tau_i<\tau
^\pi_0\}} \E _{U_i}\bigl[
\mathrm{e}^{-q T^-_\d}\mbf1_{\{T^-_\d< T^-_0\}}\bigr]\bigr].
\end{eqnarray*}
As $X$ has bounded variation, we have $\P_x(X(T^-_\d)<\d)=1$ for all
$x\in[\d,\infty)$
so that it follows that, for any $x\in[\d,\infty)$, the probability
$\P
_x(T_\d^- < T_0^-) = \P_x(0<X(T^-_\d)<\d)$
tends to zero as $\d$ tends to zero.
Lebesgue's dominated convergence theorem implies that the right-hand
side of the previous display
converges to zero when $\d$ tends to $0$.
This completes the proof of the claim in \eqref{eq:liminf}
\end{pf*}

\section{Proof of analytical optimality criterion}\label{pf:pr:key}

\begin{pf*}{Proof of Lemma~\ref{prop:key}}
(i) First consider the case
$K=0$. The proof is based on the following identity
that holds for any $c>0$ and any $x\leq b^*_++c$:
%
\begin{eqnarray}
\label{eq:vbb} &&\E_x \biggl[\mathrm{e}^{-q(t\wedge\tau)}
v_b(U_{t\wedge\tau}) + \int_{[0,t\wedge\tau]}
\mathrm{e}^{-qs}\,\td D_s \biggr] - v_b(x)
\nonumber
\\[-8pt]
\\[-8pt]
\nonumber
&&\qquad= \E_x \biggl[\int_0^{t\wedge\tau}
\mathrm{e}^{-qs} \bigl({}_{b_+} \Gamma^{\ovl w}_\infty
v_b\bigr) (U_{s-})\mbf1_{\{U_{s-}> b_+\}
}\,\td s \biggr],
\end{eqnarray}
with $b=b^*$, $b_+=b^*_+$ and $\tau=\tau^{\pi_{(b^*_-,b^*_++c)}}$,
$\ovl w=v_{b^*}$, $\mu_K=\mu_K^{{\pi_{(b^*_-,b^*_++c)}}}$,
$D=D^{\pi_{(b^*_-,b^*_++c)}}$, $U=U^{\pi_{(b^*_-,b^*_++c)}}$.
The proof of \eqref{eq:vbb} is similar to the proof of Lemma~\ref
{lem:smp}(ii) and is omitted.

Letting $t\to\infty$ in \eqref{eq:vbb} Lebesgue's dominated
convergence theorem
implies for $x\in[0,b^*_++c]$
\begin{eqnarray*}
v_{b^*+c}(x) - v_{b^*}(x) &=&\E_x \biggl[\int
_0^{\tau_{b^*+c}} \mathrm{e}^{-qs}\bigl[
{}_{b^*_+}\Gamma^{\ovl w}_\infty v_{b^*}\bigr]
\bigl(U^{b^*+c}_{s-}\bigr)\mbf1_{\{
U^{b^*+c}_{s-}> b^*_+\}}\,\td s \biggr]
\\
&= &\int_{(b^*_+,b^*_++c]} \bigl[{}_{b^*_+}\Gamma^{\ovl w}_\infty
v_b\bigr](y) R^q_{0,b^*_++c}(x,\,\td y) \qquad\mbox{with}
\\
R^q_{0,b^*_++c}(x,\td y)&=& \int_0^\infty
\mathrm{e}^{-qt}\P_x\bigl(Y^{b^*_++c}_{t}
\in\td y, t<\tau_0\bigr)\,\td t.
\end{eqnarray*}
Inserting the
explicit expressions from \eqref{eq:wapap} and Pistorius \cite{P}, Theorem~1 (see also proof of Proposition~\ref{prop:two}) for
$v_b^*$, $v_{b^*+c}$ and $R^q_{0,b^*_++c}(x,\td y)$
yields for $x\in
x \in[0,b^*_+]$
\begin{eqnarray*}
&&W^{(q)}(x)\bigl[G\bigl(b^*_++c\bigr) - G\bigl(b^*_+\bigr)\bigr]\\
&&\qquad =
W^{(q)}(x) \int\bigl[{}_{b^*_+}\Gamma ^{\ovl w}_\infty
v_{b^*}\bigr](y)\frac{W^{(q)}(b^*_++c-\td y)}{W^{(q)\prime
}(b^*_++c)},
\end{eqnarray*}
where the integral is over the interval $(b^*_+,b^*_++c]$
with $G=G_{b^*_-}$ and using that $W^{(q)}(x)$ is equal to 0 for $x<0$.
Changing coordinates
in the integral and using that $W^{(q)}(x)$ is strictly positive at any
$x>0$ yields the
first equality in \eqref{eq:key}. The second equality in
\eqref{eq:key} follows
by the representation in \eqref{eq:Fwd0}.
The second statement is a direct consequence of \eqref{eq:key}
and the fact $\{G(b^*_-, b^*_++c) < G(b^*_-, b^*_+)\ \forall c>0\}$
(from the definition of $d^*$ as \emph{last} supremum).
The proof of the case $K>0$ is similar and omitted.

The ultimate monotonicity of $G(b^-,y)$ and $G^\#(y)$ follows from the fact
that ${}_{b_+} \Gamma^w_\infty v_b(x)$ tends to minus infinity when
$x\to\infty$ (by Lemma~\ref{lem:smp}).

(ii) Taking the Laplace transform
in $c$ in \eqref{eq:key} and using the form of the Laplace
transform of $W^{(q)}$ yields that, for $\theta>\Phi(q)$ and with $G=G_{b_-}$,
\begin{eqnarray*}
\mc L g(\theta) \cdot\frac{\theta}{\psi(\theta) - q} &=& \int_{[0,\infty)}
\mathrm{e}^{-\theta c} W^{(q)\prime
}(b_++c)\bigl[G(b_++c) - G(b_+)\bigr]\,
\td c
\\
&=& \int_{[0,\infty)} \int_{[z,\infty)}
\mathrm{e}^{-\theta c} W^{(q)\prime}(b_++c)\,\td c\, G(b_++\td z)
\\
&=& \mathrm{e}^{\theta b_+} \int_{[b_+,\infty)} \int
_{[z,\infty)} \mathrm{e}^{-\theta c} W^{(q)\prime}(c)\,\td c\,
G(\td z)
\\
&=& \frac
{\mathrm{e}^{\theta
b_+}}{\psi(\theta)-q}\int_{[b_+,\infty)} \mathrm{e}^{-\theta
z}Z^{(q,\theta)\prime}(z)
G(\td z),
\end{eqnarray*}
by a change of the order of integration, which is justified by Fubini's
theorem, and the form \eqref{eq:Zd} of $Z^{(q,\theta)\prime}(z)$.
The second assertion follows since a function $f\dvtx(c,\infty)\to
\mbb R$
with $c>0$
is completely monotone if and only if it is the Laplace transform of a
nonnegative
measure supported on $\mbb R_+$.
\end{pf*}

\section{On optimality of single band strategies}\label{ssec:12band}

\begin{pf*}{Proof of Corollary~\ref{thm:expl}}
In view of verification Theorem~\ref{cor:repg}, it suffices to verify that
it holds $J(x)\leq0$ for any $x>0$ with $J(x):=
({}_{b^*_+} \Gamma^{\T w}_\infty v_{b^*})(b^*_+ + x)$.
This assertion follows once the following three facts are verified:
\begin{longlist}[(iii)]
\item[(i)] $J$ is concave on $\mbb R_+\setminus\{0\}$,
\item[(ii)] $J(0+)=0$ and
\item[(iii)] $J'(0+)\leq0$.
\end{longlist}

To show (i) note that under the stated assumptions, for $y\in(0,b)$,
$[v(b-y) - v(b) + y]\leq0 \Eq v(b)-v(b-y)\ge y$ (as $K=0$), and
for $y\ge b$ it holds $w(b-y)-v(0) -b+y\leq0$ and $v(0)-v(b)
+ b\leq0$ which yields that $w(b-y)-v(b)\leq y$ for $y\ge b$. As
$\nu'$ is convex, and a mixture of convex functions with positive
weights is again convex, it follows that $J$ is concave on
$\mbb R_+\setminus\{0\}$.

Given (ii), statement (iii) follows since
if $J'(0+)$ were positive,
$(J(x)-J(0+))/x=J(x)/x$
would be positive for all $x$ sufficiently small which would be in contradiction
with \eqref{eq:key}.

To see that (ii) holds, note that, from
\eqref{eq:key}, $\int_{[0,c]} J(c-y) W^{(q)}(\td y)\leq0$
for all $c>0$ sufficiently small.
Thus since $J$ is continuous on $\mbb R_+\setminus\{0\}$ (as it is concave)
it follows $J(0+)\leq0$. To complete the proof it is next shown that
also $J(0+)\ge0$.

First consider the case that $\sigma^2$ is strictly positive:
The observations that, for any $b>0$,
$\mathrm{e}^{-q(t\wedge T_{0,b})}v_b(X_{t\wedge T_{0,b}})$
is a martingale with $v_b\in C^2$ together with It\^{o}'s lemma
yield that $({}_{0} \Gamma^w_\infty v_{b})(x) = 0$ for all $x\in(0,b_+)$
which in turn implies that $J(0+) = {}_{0} \Gamma^w_\infty v_{b}(b_+)
= 0$
on account of the continuity of $x\mapsto({}_{0} \Gamma^w_\infty v_{b^*})(x)$
at \mbox{$x=0$}.

Consider next the case $\sigma^2=0$, which follows by approximation.
By adding a small Brownian component with variance $\s^2>0$ to $X$
and subsequently letting $\s^2\to0$, it can be shown that
in this case $J(0+)\ge0$: If $\sigma\searrow
0$, the continuity theorem implies that the scale functions
$W^{(q)(\sigma)}$ and $F^{(\sigma)}_w$ of the perturbed process
$X^{(\sigma)}:= X+\sigma B$ (where $B$ is a Brownian motion
independent of $X$) and the corresponding derivatives $W^{(q)(\sigma
)\prime}$ and $F^{(\sigma)\prime}_w$ converge
pointwise to the corresponding (derivatives of) scale functions of $X$
at any point
of continuity. Denote by $J^{(\sigma)}$ the function $J$ with the
function $v$ replaced by the function
$v^{(\sigma)}$ corresponding to the perturbed process $X^{(\sigma)}$.
An application of Fatou's lemma, which is justified on account of the
bounds in Lemma~\ref{lem:est}, then yields
that
\[
0 = \lim_{\sigma\searrow0} J^{(\sigma)}(x)\leq J(x)\qquad \mbox{for any
$x>0$}.
\]
The proof is complete.
\end{pf*}
\end{appendix}
\section*{Acknowledgments}
We are grateful to the anonymous referees for their many helpful
suggestions and careful reading, which led to improvements of the paper.

%





\printaddresses
\end{document}